\documentclass[a4paper,11pt,makeidx]{amsart}
\def\mapright#1{\smash{\mathop{\longrightarrow}\limits^{#1}}}
\usepackage{amscd}
\usepackage{xypic}  
\usepackage{amssymb}
\usepackage{amsthm}
\usepackage{epsf}
\makeindex
\newtheoremstyle{break}
  {9pt}
  {9pt}
  {\itshape}
  {}
  {\bfseries}
  {.}
  {\newline}
  {}
\theoremstyle{break}

\theoremstyle{plain}
\newtheorem{thm}{Theorem}[section]
\newtheorem{cor}[thm]{Corollary}

\newtheorem{lemma}[thm]{Lemma}

\newtheorem{prop}[thm]{Proposition}

\newtheorem{defn}[thm]{Definition} 
\newtheorem{claim}[thm]{Claim}
\newtheorem{rem}[thm]{Remark}

\renewcommand{\proofname}{Proof}
\newtheorem{notation}[thm]{Notation}
\font\goth=eusm10
\def\fake{\hbox{\goth I}}

\def\Sing{\operatorname{Sing}}
\def\Ker{\operatorname{Ker}}

\def\min{\operatorname{min}}
\def\Im{\operatorname{Im}}

\def\cork{\operatorname{cork}}
\def\max{\operatorname{max}}

\def\c1{\operatorname{c_1}}
\def\c2{\operatorname{c_2}}
\def\Cliff{\operatorname{Cliff}}
\def\gon{\operatorname{gon}}
\def\Id{\operatorname{Id}}

\def\QQ{{\mathbb Q}}
\def\PP{{\mathbb P}}

\def\L{{\mathcal L}}

\def\O{{\mathcal O}}
\def\I{{\mathcal J}}

\def\H{{\mathcal H}}

\def\X{{\mathcal X}}
\def\Ss{{\mathcal S}}

\def\eqv{\equiv}
\def\sub{\subseteq}

\def\+{\oplus}                   
\def\*{\otimes}                  
\def\hpil{\longrightarrow}       
\def\khpil{\rightarrow}

\def\Pic{\operatorname{Pic}}
\def\Num{\operatorname{Num}}

\def\hs{\hspace{.1in}}

\begin{document}

\title[Extendability of surfaces and a genus bound for Enriques-Fano threefolds]{On the
extendability of projective surfaces and \\ a genus bound for Enriques-Fano threefolds}

\dedicatory{\normalsize \dag \ Dedicated to the memory of Giulia Cerutti, Olindo Ado
Lopez and Sa\'ul S\'anchez} 

\author[A.L. Knutsen, A.F. Lopez and R. Mu\~{n}oz]{Andreas Leopold Knutsen*, Angelo Felice Lopez**
and Roberto Mu\~{n}oz***}

\address{\hskip -.43cm Andreas Leopold Knutsen, Dipartimento di Matematica, 
Universit\`a di Roma Tre, Largo San Leonardo Murialdo 1, 00146, Roma, Italy. e-mail {\tt
knutsen@mat.uniroma3.it.}}

\address{\hskip -.43cm Angelo Felice Lopez, Dipartimento di Matematica, 
Universit\`a di Roma Tre, Largo San Leonardo Murialdo 1, 00146, Roma, Italy. e-mail {\tt
lopez@mat.uniroma3.it.}}

\address{\hskip -.43cm Roberto Mu\~{n}oz, ESCET, Universidad Rey Juan Carlos, 28933 M\'ostoles 
(Madrid), Spain. \newline e-mail {\tt roberto.munoz@urjc.es}} 

\thanks{* Research partially supported by a Marie Curie Intra-European Fellowship within the 6th European
Community Framework Programme}

\thanks{** Research partially supported by the MIUR national project ``Geometria delle variet\`a algebriche"
COFIN 2002-2004.}

\thanks{*** Research partially supported by the MCYT project BFM2003-03971.}


\begin{abstract}
We introduce a new technique, based on Gaussian maps, to study the possibility, for a given
surface, to lie on a threefold as a very ample divisor with given normal bundle. We give several applications,
among which one to surfaces of general type and another one to Enriques surfaces. For the latter we prove that
any  threefold (with no assumption on its singularities) having as hyperplane section a smooth Enriques surface
(by definition an Enriques-Fano threefold)  has genus $g \leq 17$ (where $g$ is the genus of its smooth curve
sections). Moreover we find a new Enriques-Fano threefold of genus $9$ whose normalization has canonical 
but not terminal singularities and does not admit $\QQ$-smoothings.
\end{abstract}

\maketitle
 
\section{Introduction}
\label{intro}

One of the most important contributions given in algebraic geometry in the
last century, is the scheme of classification of higher dimensional varieties proposed by Mori
theory. Despite the fact that many statements still remain conjectural, several beautiful
theorems have been proved and the goal, at least in the birational realm, is particularly
clear in dimension three: starting with a threefold $X_0$ with terminal singularities and
using contractions of extremal rays, the Minimal Model Program (see e.g. \cite{km}) predicts to
arrive either at a threefold $X$ with $K_X$ nef or at a Mori fiber space, that is (\cite{re2})
there is an elementary contraction $X \to B$ with $\dim B < \dim X$. Arguably the
simplest case of such spaces is when $B$ is a point, that is $X$ is a Fano threefold. As is well
known, {\it smooth} Fano threefolds have been classified (\cite{isk1, isk2, mm}), while, in the
singular case, a classification, or at least a search for the numerical invariants, is still underway 
\cite{muk1, pr1, jpr}. 

Both the old and recent works on the classification of smooth Fano threefolds use
the important fact \cite{sh1} that a general anticanonical divisor is a smooth K3 surface. In
\cite{clm1, clm2} the authors studied varieties with canonical curve section and recovered,
in a very simple way, using the point of view of Gaussian maps, a good part of the classification
\cite{muk2}. The starting step of the latter method is Zak's theorem \cite[page 277]{za} (see also 
\cite[Thm.0.1]{lv}): If $Y \subset \PP^r$ is a smooth variety of codimension at least two with normal
bundle $N_{Y/\PP^r}$ and $h^0(N_{Y/ \PP^r}(-1)) \leq r + 1$, then the only variety $X \subset
\PP^{r + 1}$ that has $Y$ as hyperplane section is a cone over $Y$ (when this happens $Y \subset
\PP^r$ is said to be {\it nonextendable}). Now the key point in the application of this theorem
is to be able to calculate the cohomology of the normal bundle. This is of course an often
difficult task, especially when the codimension of $Y$ grows. It is here that Gaussian maps
enter the picture by giving a big help in the case of curves \cite[Prop.1.10]{wa}: if $Y$ is a curve then 
\begin{equation} 
\label{eq:gm}   
h^0(N_{Y/ \PP^r}(-1)) = r + 1 + \cork \Phi_{H_Y, \omega_Y}
\end{equation}   
where $\Phi_{H_Y, \omega_Y}$ is the Gaussian map associated to the canonical and hyperplane
bundle $H_Y$ of $Y$. For example when $X \subset \PP^{r+1}$ is a smooth anticanonically
embedded Fano threefold with general hyperplane section $Y$, in \cite[Thm.4 and Prop.3]{clm1}, $h^0(N_{Y/
\PP^r}(-1))$ was computed by calculating these coranks for the general curve section
$C$ of $Y$. 

In the case above the proof was strongly based on the fact that $C$ is a general
curve on a general K3 surface and that the Hilbert scheme of K3 surfaces is essentially
irreducible. On the other hand the latter fact is quite peculiar of K3 surfaces and we immediately
realized that if one imposes different hyperplane sections to a threefold, for example Enriques
surfaces, it becomes quite difficult to usefully rely on the curve section.

To study this and other cases it became apparent that it would be an important help
to have an analogue of the formula (\ref{eq:gm}) in higher dimension. We accomplish this in
Section \ref{key} by proving the following general result in the case of surfaces (a similar result
holds in higher dimension):

\begin{prop} 
\label{anysurface}
Let $Y \subset \PP^r$ be a smooth irreducible linearly normal surface and let $H$ be its
hyperplane bundle. Assume there is a base-point free and big line bundle $D_0$ on $Y$
with $H^1(H-D_0)=0$ and such that the general element $D \in |D_0|$ is not rational
and satisfies 
\begin{itemize}
\item[(i)] the Gaussian map $\Phi_{H_D, \omega_{D}(D_0)}$ is surjective; 
\item[(ii)] the multiplication maps $\mu_{V_D, \omega_D}$ and $\mu_{V_D, \omega_{D}(D_0)}$ are
surjective, where 

$V_D := \Im \{H^0(Y, H-D_0) \khpil H^0(D, (H-D_0)_{|D})\}$.
\end{itemize}
Then
\[ h^0(N_{Y/ \PP^r}(-1)) \leq r+1 + \cork \Phi_{H_D, \omega_D}. \]
\end{prop}
Despite the apparent complexity of the above hypotheses, it should be mentioned
that as soon as both $D_0$ and $H-D_0$ are base-point free and the degree of $D$ is large with
respect to its genus, the hypotheses are fulfilled unless $D$ is hyperelliptic. Proposition
\ref{anysurface} is therefore a flexible instrument to study threefolds whose hyperplane sections
have large Picard group. This aspect complements in a nice way the recent work of Mukai
\cite{muk1}, where a classification of Gorenstein {\it indecomposable} Fano threefolds has been
achieved: In fact indecomposable implies that a decomposition $H = D_0 + (H-D_0)$ with both
$D_0$ and $H-D_0$ moving essentially does not exist.

As we will see in Section \ref{absver}, Proposition \ref{anysurface} has several
applications. A nice sample of this is the following consequence: a pluricanonical embedding of
a surface of general type, and even some projection of it, is not, in many cases, hyperplane
section of a threefold (different from a cone) (see also Remark \ref{examples} for sharpness).
We recall that if $Y$ is a minimal surface of general type containing no (-2)-curves,
then $mK_Y$ is very ample for $m \geq 5$ \cite[Main Thm.]{bo}.

\begin{cor} 
\label{pluricanonical}
Let $Y \subset \PP V_m$ be a minimal surface of general type whose canonical bundle is base-point free and
nonhyperelliptic and $V_m \subseteq H^0(mK_Y+\Delta)$ where $\Delta \geq 0$ and either $\Delta$ is nef or
$\Delta$ is reduced and $K_Y$ is ample. Suppose that either $Y$ is regular or linearly normal and that 
\[ m \geq \begin{cases} 9 & {\rm if} \ K_Y^2 = 2; \\ 7 & {\rm if} \ K_Y^2 = 3; \\ 6 & {\rm if}
\ K_Y^2 = 4 \ {\rm and \ the \ general \ curve \ in} \ |K_Y| \ {\rm is \ trigonal \ or \ if } \ K_Y^2 = 5
\ {\rm and} \\ & {\rm the \ general \ curve \ in} \ |K_Y| \ {\rm is \ a \ plane \ quintic;}
\\ 5 & {\rm if \ either \ the \ general \ curve \ in} \ |K_Y| \ {\rm has \ Clifford \ index} \ 2 \
{\rm or} \\ & 5 \leq K_Y^2 \leq 9 \ {\rm and \ the \ general \ curve \ in} \ |K_Y| \ {\rm is \
trigonal;}
\\ 4 & {\rm otherwise}. 
\end{cases} \]  
Then $Y$ is nonextendable.
\end{cor}

Besides the mentioned applications, in the present article we will concentrate most of
our attention on the case of Enriques-Fano threefolds: in analogy with Fano threefolds where an
anticanonical divisor is a K3 surface, we define

\begin{defn} 
An {\bf Enriques-Fano threefold} is an irreducible three-dimensional variety $X \subset \PP^N$ having a 
hyperplane section $S$ that is a smooth Enriques surface, and such that $X$ is not a cone over $S$. We will say
that $X$ has genus $g$ if $g$ is the genus of its general curve section.
\end{defn}

Fano himself, in a 1938 article \cite{fa}, claimed a classification of such
threefolds, but his proof contains several gaps. Conte and Murre \cite{cm} were the
first to remark that an Enriques-Fano threefold must have some isolated singularities,
typical examples of which are quadruple points with tangent cone the cone over the
Veronese surface. Filling out some of the gaps in \cite{fa} and making some special assumptions on the
singularities, Conte and Murre recovered some of the results of Fano, but not enough to give
a classification, nor to bound the numerical invariants. On the opposite extreme, with the
strong assumption that the Enriques-Fano threefold is a quotient of a {\it smooth} Fano
threefold by an involution (this corresponds to having only cyclic quotient terminal
singularities), a list was given by Bayle \cite[Thm.A]{ba} and Sano \cite[Thm.1.1]{sa}, by using the 
classification of smooth Fano threefolds and studying which of them have such involutions.

Moreover, by the results of Minagawa \cite[MainThm2]{mi}, any Enriques-Fano threefold
with at most terminal singularities admits a $\QQ$-smoothing, that is \cite{mi, re1}, it appears as 
central fiber of a small deformation over the 1-parameter unit disk, such that a general fiber
has only cyclic quotient terminal singularities. This, together with the results of Bayle and Sano, gives
then the bound $g \leq 13$ for Enriques-Fano threefolds with at most terminal singularities.

Bayle and Sano recovered all of the examples of Enriques-Fano threefolds given by
Fano and Conte-Murre. As these were the only known examples, it has been conjectured for
some time now that this list is complete or, at least, that the genus is bounded, in analogy with the celebrated
genus bound for smooth Fano threefolds \cite{isk1, isk2, sh2}.

In Section \ref{new3-fold}, we will show that the list of known Enriques-Fano 
threefolds of Fano, Conte-Murre, Bayle and Sano {\it is not complete} (in fact not even after
specialization), by finding a new Enriques-Fano threefold enjoying several peculiar properties (for a
more precise statement and related questions, see Proposition \ref{prop:newthreefold2} and Remark \ref{remko}):

\begin{prop} 
\label{newthreefold}
There exists an Enriques-Fano threefold $X \subset \PP^9$ of genus $9$ such that neither 
$X$ nor its polarized normalization belongs to the list of Fano-Conte-Murre-Bayle-Sano. 

Moreover, $X$ does not have a $\QQ$-smoothing and in particular $X$ is not in the closure
of the component of the Hilbert scheme made of Fano-Conte-Murre-Bayle-Sano's examples. 
Its normalization $\widetilde{X}$ has canonical but not terminal singularities and does not admit
$\QQ$-smoothings.
\end{prop}

Observe that $\widetilde{X}$ is a $\QQ$-Fano threefold of Fano index $1$ with canonical singularities not having
a $\QQ$-smoothing, thus showing that Minagawa's theorem \cite[MainThm.2]{mi} cannot be extended to the canonical
case.

In the present article we apply Proposition \ref{anysurface} to get a genus bound on 
Enriques-Fano threefolds, with no assumption on their singularities:

\begin{thm}
\label{main}
Let $X \subset \PP^r$ be an Enriques-Fano threefold of genus $g$. Then $g \leq 17$.
\end{thm}

A more precise result for $g = 15$ and $17$ is proved in Proposition \ref{precisa}.

We remark that very recently Prokhorov \cite{pr1, pr2} proved the same genus bound $g \leq 17$
for Enriques-Fano threefolds and at the same time constructed an example of a new Enriques-Fano 
threefold of genus $17$ \cite[Prop.3.2]{pr2}, thus showing that the bound $g \leq 17$ is in fact
optimal. His methods are completely different from ours, in that he uses the log minimal model program
in the category of $G$-varieties and results about singularities of Enriques-Fano threefolds of Cheltsov
\cite{ch}. On the other hand, our procedure relies only on the geometry of curves on Enriques surfaces (see also
Remark \ref{rem2}). In any case, both our example in Proposition \ref{newthreefold} and Prokhorov's new examples
(in fact, he also gives a new example in genus $13$), shows that new methods were required in the classification
of Enriques-Fano threefolds with arbitrary singularities.

Now a few words on our method of proof. In Section \ref{basic} we review some basic results that will be needed 
in our study of Enriques surfaces. In Section \ref{extend} we apply  Proposition \ref{anysurface} to Enriques
surfaces and obtain the main results on nonextendability needed in the rest of the article (Propositions
\ref{mainextenr}, \ref{ramextenr}, \ref{ramextenr6} and \ref{ramextenr4}). In Section \ref{class} we prove
Theorem \ref{main} for all Enriques-Fano threefolds except for some concrete embedding line bundles on the
Enriques surface section. These are divided into different groups and then handled one by one in Sections
\ref{caseD}-\ref{caseS} by finding suitable divisors satisfying the conditions of Proposition \ref{anysurface},
thus allowing us to prove our main theorem and a more precise statement for $g = 15$ and $17$ in Section
\ref{precisazioni}. 

To prove our results it turns out that one needs effective criteria to ensure the surjectivity
of Gaussian maps on curves on Enriques surfaces and of multiplication maps of (not always complete) linear 
systems on such curves. To handle the first problem a good  knowledge of the Brill-Noether theory of a
curve lying on an Enriques surface and general in its linear system must be available. We studied this 
independent problem in another article (\cite{kl1}) and consequently obtained results ensuring the surjectivity
of Gaussian maps in \cite{klgm}  (see also Theorem \ref{tendian}). To handle the multiplication maps, we find an
effective criterion in Lemma \ref{multhelp} (which holds on any surface) in the present article.

\section{Normal bundle estimates}
\label{key}

We devise in this section a general method to give an upper bound on the cohomology of
the normal bundle of an embedded variety. We state it here only in the case of surfaces to avoid a
lengthy list of conditions.

\begin{notation}
{\rm Let $L, M$ be two line bundles on a smooth projective variety. Given $V \subseteq H^0(L)$ we
will denote by $\mu_{V, M} : V \otimes H^0(M) \hpil H^0(L \otimes M)$ the multiplication map of
sections, $\mu_{L, M}$ when $V = H^0(L)$, by $R(L, M)$ the kernel of $\mu_{L, M}$ and by
$\Phi_{L,M} : R(L, M) \hpil  H^0(\Omega^1_X \otimes L \otimes M)$ the Gaussian map (that can be
defined locally by $\Phi_{L,M}(s \otimes t) = sdt - tds$, see
\cite[1.1]{wa}).}
\end{notation}

\renewcommand{\proofname}{Proof of Proposition {\rm \ref{anysurface}}}  
\begin{proof}
To estimate $h^0(N_{Y/ \PP^r}(-1))$ we will use the exact sequence
\[0 \hpil N_{Y/\PP^r}(-D_0-H) \hpil N_{Y/\PP^r}(-H) \hpil N_{Y/\PP^r}(-H)_{|D} \hpil 0\]
and prove that
\begin{equation} 
\label{-D-H}
h^0(N_{Y/\PP^r}(-D_0-H)) = 0
\end{equation}
and 
\begin{equation}
\label{|D}
h^0(N_{Y/\PP^r}(-H)_{|D}) \leq r+1 + \cork \Phi_{H_D, \omega_D}.
\end{equation}
To prove (\ref{-D-H}), let us see first that it is enough to have
\begin{equation} 
\label{-D-HonD}
h^0(N_{Y/\PP^r}(-D_0-H)_{|D}) = 0 \hs \mbox{for a general} \ D \in |D_0|.
\end{equation}
In fact, by hypothesis there is a nonempty open subset $U \subseteq
|D_0|$ such that every $D \in U$ is smooth irreducible and satisfies (i) and (ii). Now if
$H^0(N_{Y/\PP^r}(-D_0-H))$ has a nonzero section $\sigma$, then, as $h^0(D_0) \geq 2$,
there is a nonempty open subset $U_{\sigma} \subseteq U$ such that, for every $D \in U_{\sigma}$,
we can find a point $x \in D$ with $\sigma(x) \neq 0$. The latter, of course, contradicts
(\ref{-D-HonD}).

Now (\ref{-D-HonD}) follows from the exact sequence
\[0 \hpil N_{D/Y}(-D_0-H) \hpil N_{D/\PP^r}(-D_0-H) \hpil N_{Y/\PP^r}(-D_0-H)_{|D} \hpil 0\] 
and the two conditions
\begin{equation} 
\label{D}
h^0(N_{D/\PP^r}(-D_0-H)) = 0,
\end{equation}
\begin{equation} 
\label{phiH+D}
\varphi_{H+D}: H^1(N_{D/Y}(-D_0-H)) \hpil H^1(N_{D/\PP^r}(-D_0-H)) \hs \mbox{is injective.}
\end{equation}
To see (\ref{D}), we note that the multiplication map $\mu_{H_D, \omega_{D}(D_0)}$ is
surjective by the $H^0$-lemma \cite[Thm.4.e.1]{gr}, since $|{D_0}_{|D}|$ is base-point free,
whence $D_0^2 \geq 2$, therefore $h^1(\omega_D(D_0-H)) = h^0((H-D_0)_{|D}) \leq h^0(H_D) - 2$, as
$H_D$ is very ample. Now let $\PP^k \sub \PP^r$ be the linear span of $D$. The exact sequence
\begin{equation} 
\label{normalD}
0 \hpil N_{D/\PP^k}(-D_0-H) \hpil N_{D/\PP^r}(-D_0-H) \hpil \O_D(-D_0)^{\oplus (r-k)} \hpil 0
\end{equation}
and the hypothesis $D_0^2 >0$ imply that $h^0(N_{D/\PP^r}(-D_0-H)) = h^0(N_{D/\PP^k}(-D_0-H))$. Since 
$Y$ is linearly normal and $H^1(H-D_0) = 0$, we have that also $D$ is linearly normal. As $\mu_{H_D,
\omega_{D}(D_0)}$ is surjective, by \cite[Prop.1.10]{wa}, we get that
$h^0(N_{D/\PP^k}(-D_0-H)) = \cork \Phi_{H_D, \omega_{D}(D_0)} = 0$ because of (i), and this proves \eqref{D}.
As for (\ref{phiH+D}), we prove the surjectivity of $\varphi_{H+D}^{\ast}$ with the help of the
commutative diagram 
\begin{equation}  
\label{diagram} 
\xymatrix{
H^0(\fake_{D/\PP^r}(H)) \otimes H^0(\omega_{D}(D_0)) \ar[r] \ar[d]^{f} &
H^0(N^*_{D/\PP^r}\otimes \omega_{D}(D_0+H))
\ar[d]^{\varphi_{H+D}^{\ast}} 
\\ H^0(\fake_{D/Y}(H)) \otimes H^0(\omega_{D}(D_0)) \ar[r]^{h} & H^0(N^*_{D/Y} \otimes
\omega_{D}(D_0+H)).   
}
\end{equation}
Here $f$ is surjective by the linear normality of $Y$, while $h$ factorizes as
\[ H^0(\I_{D/Y}(H)) \otimes H^0(\omega_{D}(D_0)) \twoheadrightarrow V_{D} \otimes
H^0(\omega_{D}(D_0)) \mapright{\mu_{V_D, \omega_{D}(D_0)}} H^0(N^*_{D/Y} \otimes
\omega_{D}(D_0+H)), \]
whence also $h$ is surjective by (ii). 

Finally, to prove (\ref{|D}), recall that the multiplication map $\mu_{H_D,
\omega_D}$ is surjective by \cite[Thm.1.6]{as} since $D$ is not rational, whence
$h^0(N_{D/\PP^k}(-H))  = k + 1 + \cork \Phi_{H_D, \omega_D}$ by \cite[Prop.1.10]{wa}. Therefore,
twisting (\ref{normalD}) by $\O_D(D_0)$, we get $h^0(N_{D/\PP^r}(-H)) \leq r + 1 + \cork \Phi_{H_D,
\omega_D}$ and (\ref{|D}) will now follow by the exact sequence
\[ 0 \hpil N_{D/Y}(-H) \hpil N_{D/\PP^r}(-H) \hpil N_{Y/\PP^r}(-H)_{|D} \hpil 0 \] 
and the injectivity of $\varphi_{H} : H^1(N_{D/Y}(-H)) \hpil  H^1(N_{D/\PP^r}(-H))$. The latter
follows, as in (\ref{diagram}), from the commutative diagram
\begin{equation} 
\nonumber  
\xymatrix{
H^0(\fake_{D/\PP^r}(H)) \otimes H^0(\omega_{D}) \ar[r] \ar[d] & H^0(N^*_{D/\PP^r}
\otimes \omega_{D}(H)) \ar[d]^{\varphi_{H}^{\ast}} 
\\ H^0(\fake_{D/Y}(H)) \otimes H^0(\omega_{D}) \ar[r] & H^0(N^*_{D/Y}\otimes \omega_{D}(H)) 
}
\end{equation}
by the linear normality of $Y$ and the surjectivity of $\mu_{V_D, \omega_D}$.
\end{proof}
\renewcommand{\proofname}{Proof}

\begin{rem} 
\label{keyrem1}
{\rm In the above proposition and also in Corollary \ref{cor:anysurface} below, the surjectivity of
$\mu_{V_D, \omega_{D}(D_0)}$ can be replaced by either one of the following
\begin{itemize}
\item[(i)] the multiplication map $\mu_{\omega_{D}(H-D_0), {D_0}_{|D}}$ is surjective; 
\item[(ii)] $h^0((2D_0-H)_{|D}) \leq h^0({D_0}_{|D}) - 2$;
\item[(iii)] $H.D_0 > 2D_0^2$.
\end{itemize}
{\it Proof.} The commutative diagram 
\begin{equation}  \nonumber 
\xymatrix{
V_D \otimes H^0(\omega_{D}) \otimes H^0({D_0}_{|D}) \hskip.3cm \ar[r]^{\mu_{V_D, \omega_{D}} \otimes
\Id} \ar[d] & \hskip.3cm H^0(\omega_{D}(H-D_0)) \otimes H^0({D_0}_{|D})
\ar[d]^{\mu_{\omega_{D}(H-D_0), {D_0}_{|D}}} 
\\ V_D \otimes H^0(\omega_{D}(D_0)) \ar[r]^{\mu_{V_D, \omega_{D}(D_0)}} & H^0(\omega_{D}(H))   }
\end{equation}
and the assumed surjectivity of $\mu_{V_D, \omega_{D}}$, show that (i) is enough. Now (ii) implies
(i), by the $H^0$-lemma \cite[Thm.4.e.1]{gr}, while, under hypothesis (iii), we have that
$h^0((2D_0-H)_{|D}) = 0$, whence (ii) holds.}
\end{rem}

\begin{rem} 
\label{keyrem2} 
{\rm The important conditions, in Proposition \ref{anysurface}, are the surjectivity of $\mu_{V_D, \omega_{D}}$ 
and the control of the corank of $\Phi_{H_D, \omega_D}$. On the Gaussian map side it must be said that all the
known results imply surjectivity under some conditions (more or less of the type $H.D >> g(D)$), but no good
bound on the corank is in general known. On the other hand, in many applications, one of the most important
advantages is that one can reduce to numerical conditions involving $H$ and $D_0$ (see for example Proposition
\ref{ramextenr}).} 
\end{rem}

The upper bound provided by Proposition \ref{anysurface} can be applied in many
instances to control how many times $Y$ can be extended to higher dimensional varieties.
However the most interesting and useful application will be to one simple extension. 

\begin{cor} 
\label{cor:anysurface}
Let $Y \subset \PP^r$ be a smooth irreducible surface which is either linearly normal or regular
(that is, $h^1(\O_Y) = 0$) and let $H$ be its hyperplane bundle. Assume there is a base-point
free and big line bundle $D_0$ on $Y$ with $H^1(H-D_0) = 0$ and such that the general element $D
\in |D_0|$ is not rational and satisfies
\begin{itemize}
\item[(i)] the Gaussian map $\Phi_{H_D, \omega_D}$ is surjective; 
\item[(ii)] the multiplication maps $\mu_{V_D, \omega_D}$ and $\mu_{V_D, \omega_{D}(D_0)}$ are
surjective, where 

$V_D := \Im \{H^0(Y, H-D_0) \khpil H^0(D, (H-D_0)_{|D})\}$.
\end{itemize}
Then $Y$ is nonextendable. 
\end{cor}

\begin{proof}
Note that $g(D) \geq 2$, else $\Phi_{H_D, \omega_D}$ is not surjective. Also since $\mu_{V_D,
\omega_{D}(D_0)}$ is surjective, we must have that $V_D$ (whence also $|(H-D_0)_{|D}|$) is base-point
free, as $|\omega_{D}(H)|$ is such. Therefore $2g(D) - 2 + (H-D_0).D > 0$, whence
$h^1(\omega_{D}^2(H-D_0)) = 0$ and the $H^0$-lemma \cite[Thm.4.e.1]{gr} implies that the
multiplication map $\mu_{\omega_{D}^2(H), {D_0}_{|D}}$ is surjective. Now (i) and the commutative diagram
\begin{equation} \nonumber    
\xymatrix{
R(H_D, \omega_{D}) \otimes H^0({D_0}_{|D}) \hskip.3cm \ar@{->>}[r]^{\Phi_{H_D, \omega_D}
\otimes \Id} 
\ar[d] & \hskip.3cm H^0(\omega_{D}^2(H)) \otimes H^0({D_0}_{|D}) \ar@{->>}[d]^{\mu_{\omega_{D}^2(H),
{D_0}_{|D}}} \\ R(H_D, \omega_{D}(D_0)) \ar[r]^{\Phi_{H_D, \omega_{D}(D_0)}} & H^0(\omega_{D}^2(H+D_0))  
}
\end{equation}
give that also $\Phi_{H_D, \omega_{D}(D_0)}$ is surjective.

If $Y$ is linearly normal the result therefore follows by Zak's theorem \cite[page 277]{za}, \cite[Thm.0.1]{lv},
and Proposition \ref{anysurface}. 

Assume now that $h^1(\O_Y) = 0$ and that $Y \subset \PP^r$ is extendable, that is, there
exists a nondegenerate threefold $X \subset \PP^{r+1}$ which is not a cone over $Y$ and such that
$Y = X \cap \PP^r$ is a hyperplane section. Let $\pi: \widetilde{X} \rightarrow X$ be a
resolution of singularities and let $L = \pi^*\O_X(1)$ and $\widetilde{Y} = \pi^{-1}(Y)$. 
Since $Y$ is smooth we have $Y \cap \Sing X = \emptyset$, whence there is an isomorphism
$(\widetilde{Y}, L_{|\widetilde{Y}}) \cong (Y, \O_Y(1))$. Now $L$ is nef and birational,
whence $H^1(\widetilde{X}, - L) = 0$ by Kawamata-Viehweg vanishing. Moreover, as
$\widetilde{Y} \in |L|$, we have, for all $k$, an exact sequence
\[0 \rightarrow  k L \rightarrow  (k+1)L \rightarrow (k+1)L_{|\widetilde{Y}} \rightarrow 0. \]

Setting $k = -1$ we get that $H^1(\O_{\widetilde{X}}) \subseteq H^1(\O_{\widetilde{Y}}) =
H^1(\O_Y) = 0$, therefore, setting $k=0$, we deduce the surjectivity of the restriction map 
$H^0(\widetilde{X},L)  \rightarrow  H^0(\widetilde{Y}, L_{|\widetilde{Y}})$. 

Consider the birational map $\varphi_L : \widetilde{X} \rightarrow \PP^N$ where $N
\geq r+1$, let $\overline{X} = \varphi_L(\widetilde{X})$ and let $\overline{Y}$ be the
hyperplane section of $\overline{X}$ corresponding to $\widetilde{Y} \in |L|$. Now
$\overline{Y} \cong Y$ and the embedding $\overline{Y} \subset \PP^{N-1}$ is given by the complete
linear series $|\O_Y(1)|$. Note also that, by construction, $\overline{X} \subset \PP^N$
projects to $X \subset \PP^{r+1}$, whence $\overline{X}$ is not a cone over $\overline{Y}$.
Therefore $\overline{Y} \subset \PP^{N-1}$ is linearly normal and extendable. But also on
$\overline{Y}$ we have a line bundle $\overline{D_0}$ satisfying the same properties as $D_0$,
whence, by the proof in the linearly normal case, $\overline{Y}$ is nonextendable, a contradiction.
\end{proof}

\section{Absence of Veronese embeddings on threefolds} 
\label{absver}

It was already known to Scorza in 1909 \cite{sc} that the Veronese varieties
$v_m(\PP^n)$ are nonextendable for $m > 1$ and $n > 1$. For a Veronese embedding of any variety we
can use Zak's theorem to deduce nonextendability, as follows (we omit the case of
curves that can be done, as is well-known, via Gaussian maps)

\begin{rem} 
Let $X \subset \PP^r$ be a smooth irreducible nondegenerate n-dimensional variety, $n \geq 2$, $L = \O_X(1)$ 
and let $\varphi_{mL}(X) \subset \PP^N$ be the m-th Veronese embedding of $X$. 

If $H^1(T_X(-mL)) = 0$ then $\varphi_{mL}(X)$ is nonextendable. In particular the latter
holds if
\[ m > \max \{2, n + 2 + \frac{K_X.L^{n-1} - 2r + 2n + 2}{L^n}\}. \]
\begin{proof} {\rm Set $Y = \varphi_{mL}(X)$. From the exact sequences
\[0 \hpil \O_Y(-1) \hpil \O_Y^{\oplus (N+1)} \hpil T_{\PP^N}(-1)_{|Y} \hpil 0\]
\[0 \hpil T_Y(-1) \hpil T_{\PP^N}(-1)_{|Y} \hpil N_{Y/\PP^N}(-1) \hpil 0\] 
and Kodaira vanishing we deduce that $h^0(N_{Y/\PP^N}(-1)) \leq h^0(T_{\PP^N}(-1)_{|Y}) +
h^1(T_Y(-1)) = N + 1 + h^1(T_X(-mL)) = N + 1$, and it just remains to apply Zak's theorem
\cite[page 277]{za}, \cite[Thm.0.1]{lv}.

To see the last assertion observe that since $n \geq 2$ and $m \geq 3$ we have, as is well-known,
$h^1(T_X(-mL)) = h^0(N_{X/\PP^r}(-mL))$.
Now if the latter were not zero, the same would hold for a general hyperplane section $X \cap
H$ of $X$ and so on until the curve section $C \subset \PP^{r-n+1}$. Now taking $r-n-1$ general
points $P_{j} \in C$, we have an exact sequence \cite[2.7]{bel}
\[0 \hpil \bigoplus\limits_{j=1}^{r-n-1} \O_C(1-m)(2 P_{j}) \hpil N_{C/\PP^{r-n+1}} (-m)
\hpil  \omega_C (3-m)(-2 \sum\limits_{j=1}^{r-n-1} P_{j}) \hpil 0\] 
whence $h^0(N_{C/\PP^{r-n+1}}(-m)) = 0$ for reasons of degree.} 
\end{proof}
\end{rem}

In the case of surfaces, as an application of Corollary \ref{cor:anysurface}, we can give
an extension of the above remark to multiples of big and nef line bundles. 

\begin{defn} 
\label{m(L)}
Let $Y$ be a smooth surface and let $L$ be an effective line bundle on $Y$ such that the general
divisor $D \in |L|$ is smooth and irreducible. We say that $L$ is {\bf hyperelliptic, trigonal},
etc., if $D$ is such. We denote by $\Cliff(L)$ the Clifford index of $D$.
Moreover, when
$L^2 > 0$, we set 
\[ \varepsilon(L) =  \begin{cases} 3 & {\rm if} \ L \ {\rm is \ trigonal}; \\ 
5 & {\rm if} \ \Cliff(L) \geq 3; \\ 
0 & {\rm if} \ \Cliff(L) = 2.
\end{cases} \]
and
\[ m(L) =  \begin{cases} \frac{16}{L^2} & {\rm if} \  L.(L+K_Y) = 4; \\ \frac{25}{L^2} &
{\rm if} \ L.(L+K_Y) = 10 \ {\rm and \ the \ general \ divisor \ in} \ |L| \ {\rm is \ a \ plane \
quintic}; \\ \frac{3L.K_Y + 18}{2L^2} + \frac{3}{2} & {\rm if} \ 6 \leq L.(L+K_Y) \leq 22 \
{\rm and} \ L \ {\rm is \ trigonal}; \\ \frac{2 L.K_Y - \varepsilon(L)}{L^2} + 2 & {\rm otherwise}.
\end{cases} \]
\end{defn}

\begin{cor} 
\label{veronese} 
Let $Y \subset \PP V$ be a smooth surface with $V \subseteq H^0(mL+\Delta)$, where $L$ is a base-point
free, big, nonhyperelliptic line bundle on $Y$ with $L.(L+K_Y) \geq 4$ and $\Delta \geq 0$ is a divisor. Suppose
that either $Y$ is regular or linearly normal and that $m$ is such that $H^1((m-2)L+\Delta) = 0$ and 
$m > \max \{m(L) - \frac{L.\Delta}{L^2}, \lceil \frac{L.K_Y + 2 - L.\Delta}{L^2} \rceil + 1\}$. Then $Y$ is
nonextendable. 
\end{cor}
\begin{proof} 
We apply Corollary \ref{cor:anysurface} with $D_0 = L$ and $H = mL+\Delta$. By hypothesis the general $D \in
|L|$ is smooth and irreducible of genus $g(D) = \frac{1}{2}L.(L+K_Y) + 1$. Since $H^1(H-2L) = 0$,
we have $V_D = H^0((H-L)_{|D})$. Also $(H-L).D = (m-1)L^2 + L.\Delta \geq L.(L+K_Y) + 2 = 2g(D)$ by
hypothesis, whence $|(H-L)_{|D}|$ is base-point free and birational (as $D$ is not hyperelliptic) and the
multiplication map $\mu_{V_D, \omega_D}$ is surjective by \cite[Thm.1.6]{as}. Moreover
$H^1((H-L)_{|D}) = 0$, whence also $H^1(H-L) = 0$ by the exact sequence
\[0 \hpil H-2L \hpil H-L \hpil (H-L)_{|D} \hpil 0. \]
 
The surjectivity of $\mu_{V_D, \omega_{D}(L)}$ now follows by \cite[Cor.4.e.4]{gr} since $\deg \omega_{D}(L)
\geq 2g(D) + 1$ because $L^2 \geq 3$: If $L^2 \leq 2$ we have that $h^0(L_{|D}) \leq 1$ as $D$ is not 
hyperelliptic, whence $h^0(L) \leq 2$, contradicting the hypotheses on $L$. Finally the surjectivity of
$\Phi_{H_D, \omega_D}$ follows by the inequality $m > m(L) - \frac{L.\Delta}{L^2}$ and well-known results about
Gaussian maps (\cite[Prop.1.10]{wa}, \cite[Prop.2.9, Prop.2.11 and Cor.2.10]{klgm}, \cite[Thm.2]{bel}).
\end{proof}

\begin{rem} 
{\rm The above result does apply, in some instances, already for $m = 1$ or $2$. Also observe that the
base-point free ample and hyperelliptic line bundles are essentially classified by several results in
adjunction theory (see \cite{bs} and references therein).} 
\end{rem}

We can be a little bit more precise in the interesting case of pluricanonical
embeddings.

\renewcommand{\proofname}{Proof of Corollary {\rm \ref{pluricanonical}}}
\begin{proof} 
We apply Corollary \ref{veronese} with $L = K_Y$ and $H = mK_Y+\Delta$, and we just need to check that
$H^1((m-2)K_Y+\Delta) = 0$. If $\Delta$ is nef this follows by Kawamata-Viehweg vanishing. Now
suppose that $\Delta$ is reduced and $K_Y$ is ample. Again by Kawamata-Viehweg vanishing we have that
$H^1((m-2)K_Y) = 0$, whence  the exact sequence
\[ 0 \hpil (m-2)K_Y \hpil (m-2)K_Y+\Delta \hpil \O_{\Delta} ((m-2)K_Y+\Delta) \hpil 0 \] 
shows that $H^1((m-2)K_Y+\Delta) = 0$ since $h^1(\O_{\Delta} ((m-2)K_Y+\Delta)) = 
h^0(\O_{\Delta}(-(m-3)K_Y)) = 0$. 
\end{proof}
\renewcommand{\proofname}{Proof}

\begin{rem}
\label{examples} 
{\rm Consider the $5$-uple embedding $X$ of $\PP^3$ into $\PP^{55}$. A general hyperplane 
section of $X$ is a smooth surface $Y$ embedded with $5K_Y$ and satisfying $K_Y^2 = 5$. Also
consider the $4$-uple embedding of a smooth quadric hypersurface in $\PP^4$ into $\PP^{54}$. Its
general hyperplane section is a smooth surface $Y$ embedded with $4K_Y$ and satisfying
$K_Y^2 = 8$. Therefore, in general, the conditions on $K_Y^2$ and $m$ cannot be weakened.}
\end{rem}

\begin{rem} 
{\rm If $K_Y$ is hyperelliptic, then $2K_Y$ is not birational and these surfaces have been
classified by the work of several authors (see \cite{bcp} and references therein).} 
\end{rem}

We can be even more precise in the interesting case of adjoint embeddings.

\begin{cor} 
\label{adjoint}
Let $Y \subset \PP V$ be a minimal surface of general type with base-point free and nonhyperelliptic canonical
bundle and $V \subseteq H^0(K_Y + L + \Delta)$, where $L$ is a line bundle on $Y$ and $\Delta \geq 0$ is a
divisor. Suppose that $Y$ is either regular or linearly normal, that $H^1(L + \Delta - K_Y) = 0$ and that 
\[ L.K_Y + K_Y.\Delta > \begin{cases} 14 & {\rm if} \  K_Y^2 = 2; \\ 20 & {\rm if} \ K_Y^2 = 5 \ {\rm and \ the \
general \ divisor \ in} \ |K_Y| \ {\rm is \ a \ plane \ quintic}; \\ 2K_Y^2 + 9 & {\rm if} \ 3 \leq
K_Y^2 \leq 11 \ {\rm and} \ K_Y \ {\rm is \ trigonal}; \\ 3K_Y^2 - \varepsilon(K_Y) & {\rm otherwise}.
\end{cases}. \]
Then $Y$ is nonextendable.
\end{cor}

\begin{proof} Similar to the proof of Corollary \ref{veronese} with $D_0 = K_Y$ and $H = K_Y + L + \Delta$.
\end{proof}

To state the pluriadjoint case, given a big line bundle $L$ on a smooth surface $Y$ we
define the function
\[ \nu(L) = \begin{cases} \frac{12}{L^2} + 1 & {\rm if} \ L.(L+K_Y) = 4; \\ \frac{15}{L^2} + 1 &
{\rm if} \ L.(L+K_Y) = 10 \ {\rm and \ the \ general \ divisor \ in} \ |L| \ {\rm is \ a \ plane \
quintic}; \\ \frac{L.K_Y + 18}{2L^2} + \frac{3}{2} & {\rm if} \ 6 \leq L.(L+K_Y) \leq 22 \ {\rm
and} \ L \ {\rm is \ trigonal}; \\ \frac{L.K_Y - \varepsilon(L)}{L^2} + 2 & {\rm otherwise}.
\end{cases} \]
\begin{cor} 
\label{pluriadjoint}
Let $Y \subset \PP V$ be a smooth surface with $V \subseteq H^0(K_Y + mL + \Delta)$ where $L$ is a
base-point free, big and nonhyperelliptic line bundle on $Y$ with $L.(L+K_Y) \geq 4$ and $\Delta \geq 0$
is a divisor such that $H^1(K_Y + (m-2)L + \Delta) = 0$. Suppose
that $Y$ is either regular or linearly normal and that $m > \max\{2 + \frac{1}{L^2}, \nu(L) \} -
\frac{L.\Delta}{L^2}$. Then $Y$ is nonextendable.
\end{cor}

\begin{proof} Similar to the proof of Corollary \ref{veronese} with $D_0 = L$ and $H = K_Y + mL + \Delta$.
\end{proof}

\section{Basic results on line bundles on Enriques surfaces}
\label{basic}

\begin{defn}
Let $S$ be an Enriques surface. We denote by $\sim$ (respectively $\eqv$) the linear
(respectively numerical) equivalence of divisors (or line bundles) on $S$. A line bundle $L$ is {\bf primitive} 
if $L \eqv hL'$ for some line bundle $L'$ and some integer $h$, implies $h = \pm 1$. A {\bf nodal} curve on $S$
is a smooth rational curve. A {\bf nodal cycle} is a divisor $R>0$ such that, for any $0 < R' \leq R$ we have
$(R')^2 \leq -2$. An {\bf isotropic divisor} $F$ on $S$ is a divisor such that $F^2 = 0$ and $F \not\eqv 0$. An
{\bf isotropic $k$-sequence} is a set $\{f_1, \ldots, f_k\}$ of isotropic divisors such that $f_i.f_j=1$ for $i
\neq j$.
\end{defn} 

We will often use the fact that if $R$ is a nodal cycle, then $h^0(\O_S(R)) = 1$ and $h^0(\O_S(R+K_S)) = 0$. 

Let $L$ be a line bundle on $S$ with $L^2 >0$. Following \cite{cd} we define
\[ \phi(L) =\inf \{|F.L| \; : \; F \in \Pic S, F^2=0, F \not\eqv 0\}. \]

Two important properties of this function, which will be used throughout the article, are
that $\phi(L)^2 \leq L^2$ \cite[Cor.2.7.1]{cd} and that, if $L$ is nef, then there exists a
genus one pencil $|2E|$ such that $E.L = \phi(L)$ (\cite[2.11]{cos} or by \cite[Cor.2.7.1,
Prop.2.7.1 and Thm.3.2.1]{cd}). Moreover we will extensively use (often without further mentioning) the fact
that a nef line bundle $L$ with $L^2 \geq 4$ is base-point free if and only if
$\phi(L) \geq 2$ \cite[Prop.3.1.6, 3.1.4 and Thm.4.4.1]{cd}.

\begin{lemma}
\label{nefred} 
Let $S$ be an Enriques surface, let $L$ be a line bundle on $S$ such that $L > 0$ and $L^2>0$ and let $F$ be an
effective divisor on $S$ such that $F^2=0$ and $\phi(L)=|F.L|$. Moreover let $A,B$ be two effective divisors on 
$S$ such that $A^2 \geq 0$ and $B^2 \geq 0$. Then
\begin{itemize}
\item[(a)] $F.L>0$;
\item[(b)] if $\alpha$ is a positive integer such that $(L-\alpha F)^2 \geq 0$, then 
$L-\alpha F > 0$;
\item[(c)] $A.B \geq 0$ with equality if and only if there exists a primitive divisor $D > 0$ and
integers $a \geq 1, b \geq 1$ such that $D^2 = 0$ and $A \eqv aD, B \eqv bD$.
\end{itemize}
\end{lemma}  
\begin{proof}
For parts (a) and (b) see \cite[Lemma2.5]{kl1}. Part (c) is proved in \cite[Lemma2.1]{klvan}.
\end{proof}

We will often use the ensuing

\begin{lemma} 
\label{fi}
For $1 \leq i \leq 4$ let $F_i>0$ be four isotropic divisors such that $F_1.F_2 = F_3.F_4 = 1$ and $F_1.F_3 =
F_2.F_3 = 2$. If $F_4.(F_1+F_2) = 4$ then $F_1.F_4 = F_2.F_4 = 2$.
\end{lemma}

\begin{proof}
By symmetry and Lemma \ref{nefred} we can assume, to get a contradiction, that $F_1.F_4 = 1$ and $F_2.F_4 = 3$. 
Then $(F_2 + F_4)^2 = 6$ and $\phi(F_2 + F_4) = 2$ whence, by Lemma \ref{nefred}, we can write $F_2 + F_4 \sim
A_1 + A_2 + A_3$ with $A_i > 0$, $A_i^2 = 0$ and $A_i.A_j = 1$ for $i \neq j$. But this gives the contradiction
$8 = (F_2 + F_4).(F_1 + F_2 + F_3) \geq 3\phi(F_1 + F_2 + F_3) = 9$. 
\end{proof}
 
\begin{lemma} 
\label{lemma0}
Let $L > 0$ be a line bundle on an Enriques surface with $L^ 2 \geq 0$. Then there exist 
(not necessarily distinct) divisors $F_i > 0, 1 \leq i \leq m$, such that $F_i^2=0$ and $L
\sim F_1 + \ldots + F_m$.
\end{lemma}

\begin{defn} 
\label{ar0} 
We call such a decomposition of $L$ an {\bf arithmetic genus $1$ decomposition}.
\end{defn}

\renewcommand{\proofname}{Proof of Lemma {\rm \ref{lemma0}}} 
\begin{proof}
The assertion being clear for $L^ 2 = 0$ we suppose $L^ 2 > 0$. By Lemma \ref{nefred} there is an $F_1 > 0$ such
that $F_1.L = \phi(L)$. Since $\phi(L) \leq \lfloor \sqrt{L^ 2} \rfloor$, we have $(L-F_1)^2 \geq
0$. Again by Lemma \ref{nefred} we have $L-F_1 > 0$ and $(L-F_1)^2 < L^2$ so that we can proceed by induction.
\end{proof}
\renewcommand{\proofname}{Proof}

\begin{defn} 
\label{ST}
Let $L \geq 0$ be an effective line bundle on an Enriques surface with $L^2 \geq 0$. Then 
$L$ is said to be of {\bf small type} if either $L=0$ or for every decomposition of $L$ with
\[ L \eqv a_1E_1 + \ldots + a_rE_r, \hs E_i >0, \hs E_i^2=0, \hs E_i.E_j >0 \mbox{ for } i
\neq j, \] 
and $a_i>0$, we have that $a_i=1$ for all $1 \leq i \leq r$. 
\end{defn}

The next two results are immediate consequences of the previous ones.

\begin{lemma} 
\label{STlemma1}
Let $L \geq 0$ be an effective line bundle on an Enriques surface with $L^2 \geq 0$. Then $L$ is of small type 
if and only if either {\rm (i)} $L^ 2=0$ and $L$ is either trivial or primitive, or {\rm (ii)}
$L^2 >0$ and $(L-2F)^2 < 0$ for any $F > 0$ with $F^2=0$ and $F.L=\phi(L)$.
\end{lemma}

\begin{lemma} 
\label{STlemma2}
Let $L \geq 0$ be an effective line bundle on an Enriques surface with $L^2 \geq 0$. Then $L$ is of small type 
if and only if it is of one of the following types (where $E_i >0$, $E_i^2=0$ and $E_i$ primitive): {\rm (a)}
$L=0$; {\rm (b)} $L^2=0$, $L \sim E_1$; {\rm (c)} $L^2=2$, $L \sim E_1+E_2$, $E_1.E_2=1$;  {\rm (d)} $L^2=4$,
$\phi(L)=2$, $L \sim E_1+E_2$, $E_1.E_2=2$; {\rm (e)} $L^2=6$, $\phi(L)=2$, $L \sim E_1+E_2+E_3$,
$E_1.E_2=E_1.E_3=E_2.E_3=1$; {\rm (f)} $L^2=10$, $\phi(L)=3$, $L \sim E_1+E_2+E_3$, $E_1.E_2=1$,
$E_1.E_3=E_2.E_3=2$.
\end{lemma}

Given an effective line bundle $L$ with $L^2 >0$, among all arithmetic genus $1$
decompositions of $L$ we want to choose the most convenient for us (in a sense that will be
clear in the following sections). To this end let us first record the following

\begin{lemma} 
\label{defalpha}
Let $L > 0$ be a line bundle on an Enriques surface such that $L^2 >0$ and suppose there exists an $F > 0$ with 
$F^2 = 0, \phi(L) = F.L$ and $(L-2F)^2 > 0$. Then there exist an integer $k \geq 2$ and an $F' > 0$ with 
$(F')^2 = 0, F'.F > 0, (L-kF)^2 > 0$ and $F'.(L-kF) = \phi(L-kF)$.
\end{lemma}

\begin{proof} As $(L-2F)^2 > 0$ we can choose an integer $k \geq 2$ such that $(L-kF)^2 > 0$ and $(L-(k+1)F)^2 
\leq 0$. Set $L' = L - k F$, so that $L' > 0$ by Lemma \ref{nefred} and pick any $F'>0$ such that $(F')^2 = 0$
and $F'.L' = \phi(L')$. If $(L'-F')^2 > 0$ then $F' \not\eqv F$, whence $F'.F > 0$ by Lemma \ref{nefred} and we
are done.

If $(L'-F')^2 \leq 0$ one easily sees that $((L')^2, \phi(L')) = (2, 1)$ or $(4, 2)$ and $L' \sim F' + F''$ with
$F''>0$, $(F'')^2 = 0$ and $F'.F'' = 1, 2$. Hence either $F'.F > 0$ or $F''.F > 0$ and we are again done.
\end{proof}

Now for any line bundle $L > 0$ which is {\it not of small type} with $L^2 > 0$ and $\phi(L) =
F.L$ for some $F > 0$ with $F^2 = 0$, define
\begin{eqnarray}
\label{alpha1}
\alpha_F(L) = \min \{ k \geq 2 & | & (L-kF)^2\geq 0 \mbox{ and if } (L-kF)^2 > 0 \mbox{ there
exists } F' > 0 \mbox{ with } \\ \nonumber & & (F')^2=0, F'.F>0  \mbox{ and } F'.(L-kF) \leq \phi(L) \}.
\end{eqnarray}

By Lemma \ref{defalpha}, $\alpha_F(L)$ exists and it is easily seen that an equivalent
definition is 
\begin{eqnarray}
\label{alpha2}
\alpha_F(L) = \min \{ k \geq 2 & | & (L-kF)^2\geq 0 \mbox{ and if } (L-kF)^2 > 0 \mbox{ there
exists } F' > 0 \mbox{ with } \\ \nonumber & & (F')^2 = 0, F'.F>0  \mbox{ and } F'.(L-kF) = \phi(L-kF)
\}.
\end{eqnarray}

If $L^2 = 0$ and $L$ is {\it not of small type}, then let $k \geq 2$ be the maximal integer
such that there there exists an $F > 0$ with $F^2=0$ and $L \eqv kF$. In this case we define 
$\alpha_F(L) = k$.

\begin{lemma} 
\label{STlemma3}
Let $L > 0$ be a line bundle not of small type with $L^2 > 0$ and $(L^2, \phi(L)) \neq (16,4)$,
$(12,3)$, $(8,2)$, $(4,1)$. Then $(L-\alpha_F(L)F)^2 >0$. 
\end{lemma}

\begin{proof}
Set $\alpha=\alpha_F(L)$. Assume that $(L-\alpha F)^2=0$. Then, since $L^2 > 0$, we have $L
\sim \alpha F+F'$ for some $F' > 0$ with $(F')^2=0$ and $F.F'= \phi(L)$. Now 
$(L-(\alpha-1) F)^2 = (F+F')^2 = 2\phi(L) > 0$ and $F'. (L-(\alpha-1) F) = \phi(L)$, whence
$\alpha = 2$. Therefore $L \sim 2F+F'$, whence $L^2 = 4F.F'= 4\phi(L)$, which gives $\phi(L)^2
\leq 4\phi(L)$, in other words $\phi(L) \leq 4$ and we are done.
\end{proof}

Finally we recall a definition and some results, proved in \cite{kl1} and \cite{klvan}, that 
will be used throughout the article.

\begin{lemma} \cite[Lemma2.4]{kl1}
\label{A}
Let $L > 0$ and  $\Delta > 0$ be divisors on an Enriques surface with $L^2 \geq 0$, $\Delta^2 =
-2$ and $k := - \Delta.L > 0$. Then there exists an $A > 0$ such that $A^2 = L^2$, $A.\Delta=k$ and 
$L \sim A + k \Delta$. Moreover, if $L$ is primitive, then so is $A$. 
\end{lemma}

\begin{defn} 
\label{def:qnef}
An effective line bundle $L$ on a K3 or Enriques surface is said to be {\bf quasi-nef} if $L^2
\geq 0$ and $L.\Delta \geq -1$ for every $\Delta$ such that $\Delta > 0$ and $\Delta^2 = -2$. 
\end{defn}

\begin{thm} \cite[Cor.2.5]{klvan} 
\label{lemma:qnef}
An effective line bundle $L$ on a K3 or Enriques surface is quasi-nef if and only if $L^2 \geq 0$
and either $h^1(L) = 0$ or $L \eqv nE$ for some $n \geq 2$ and some primitive and nef divisor $E > 0$
with $E^2 = 0$.
\end{thm}

We will often make use of the following simple

\begin{lemma} 
\label{qnef0}
Let $L$ be a nef and big line bundle on an Enriques surface and let $F$ be a divisor satisfying $F.L < 2\phi(L)$
(respectively $F.L=\phi(L)$ and $L$ is ample). Then $h^0(F) \leq 1$ and if $F>0$ and $F^2 \geq 0$ we have
$F^2=0$, $h^0(F)=1$, $h^1(F)=0$ and $F$ is primitive and quasi-nef (resp. nef).
\end{lemma} 

\begin{proof}
If $h^0(F) \geq 2$ we can write $|F| = |M| + G$, with $M$ the moving part and $G \geq 0$ the fixed part of $|F|$. By
\cite[Prop.3.1.4]{cd} we get $F.L \geq 2\phi(L)$, a contradiction. Then $h^0(F) \leq 1$ and if $F>0$ and 
$F^2 \geq 0$ it follows that  $F^2=0$ and $h^1(F)=0$ by Riemann-Roch. Hence $F$ is quasi-nef and primitive by
Theorem \ref{lemma:qnef}. If $F.L=\phi(L)$, $L$ is ample and $F$ is not nef, by Lemma \ref{A} we can write $F
\sim F_0 + \Gamma$ with $F_0>0$, $F_0^2=0$ and $\Gamma$ a nodal curve. But then $F_0.L < \phi(L)$.
\end{proof}

\section{Main results on extendability of Enriques surfaces}
\label{extend}

It is well-known that abelian and hyperelliptic surfaces are nonextendable (see for
example \cite[Rmk.3.12]{glm}). In the case of K3 surfaces the extendability problem is open,
but beautiful answers are known for general K3's (even with assigned Picard lattice) (\cite{clm1,
clm2, be}). Let us deal now with Enriques surfaces.

We will state here a simplification of Corollary \ref{cor:anysurface} that will be a
central ingredient for us. An analogous result can be written for K3 surfaces.

\begin{prop} 
\label{mainextenr}
Let $S \subset \PP^r$ be an Enriques surface and denote by $H$ its hyperplane section. Suppose we
can find a nef and big (whence $>0$) line bundle $D_0$ on $S$ with $\phi(D_0) \geq 2,
H^1(H-D_0)=0$ and such that the following conditions are satisfied by the general element $D \in
|D_0|$:
\begin{itemize}
\item[(i)] the Gaussian map $\Phi_{H_D, \omega_D}$ is surjective; 
\item[(ii)] the multiplication map $\mu_{V_D, \omega_D}$ is surjective, where

$V_D := \Im \{H^0(S, \O_S(H-D_0)) \khpil H^0(D, \O_D(H-D_0))\}$;
\item[(iii)] $h^0 (\O_D(2D_0 - H)) \leq \frac{1}{2} D_0^2 - 2$.
\end{itemize}
Then $S$ is nonextendable.
\end{prop}

\begin{proof}
Note that $D_0^2 \geq \phi(D_0)^2 \geq 4$. Now the line bundle $D_0$ is base-point free 
since $\phi(D_0) \geq 2$ by \cite[Prop.3.1.6, 3.1.4 and Thm.4.4.1]{cd}.
Therefore we just apply Corollary \ref{cor:anysurface} and Remark \ref{keyrem1}.
\end{proof} 

Our first observation will be that, for many line bundles $H$, a line bundle $D_0$
satisfying the conditions of Proposition \ref{mainextenr} can be found with the help of Ramanujam's
vanishing theorem.

\begin{prop} 
\label{ramextenr}
Let $S \subset \PP^r$ be an Enriques surface, denote by $H$ its hyperplane section and assume that
$H$ is not $2$-divisible in $\Num S$. Suppose there exists an effective line bundle
$B$ on $S$ with the following properties:
\begin{itemize}
\item[(i)] $B^2 \geq 4$ and $\phi(B) \geq 2$, 
\item[(ii)] $(H-2B)^2 \geq 0$ and $H-2B \geq 0$,
\item[(iii)] $H^2 \geq 64$ if $B^2 =4$ and $H^2 \geq 54$ if $B^2 =6$.
\end{itemize} 
Then $S$ is nonextendable.
\end{prop}

\begin{proof}
We first claim that we can find a {\it nef} divisor $D'>0$ still satisfying (i)-(iii) with $D' \leq B$,
$(D')^2=B^2$ and $\phi(D')=\phi(B)$ by using {\it Picard-Lefschetz reflections}. 

Recall that if $\Gamma$ is a nodal curve on an Enriques surface, then the Picard-Lefschetz reflection 
with respect to $\Gamma$ acting on $\Pic S$ is defined as $\pi_{\Gamma}(L):= L + (L.\Gamma) \Gamma$. It
is straightforward to check that $\pi_{\Gamma}(\pi_{\Gamma}(L)) = L$ and that
$\pi_{\Gamma}(L).\pi_{\Gamma}(L') = L.L'$ for any $L, L' \in \Pic S$. Moreover $\pi_{\Gamma}$ preserves
effectiveness \cite[Prop.VIII.16.3]{bpv} and the function $\phi$, when $L^2 >0$.

Now if $B$ is not nef, then there is a nodal curve $\Gamma$ such that $\Gamma.B <0$.
By the properties of $\pi_{\Gamma}$ just mentioned and the fact that
clearly $0 < \pi_{\Gamma}(B) < B$, 
it follows that $\pi_{\Gamma}(B)^2=B^2$ and $\phi(\pi_{\Gamma}(B))=\phi(B)$, whence 
$\pi_{\Gamma}(B)$ still satisfies (i)-(iii). If $\pi_{\Gamma}(B)$ is not nef, we 
repeat the process, which must eventually end, as $\pi_{\Gamma}(B) < B$.

We have therefore found the desired nef divisor $D'$. 

Since $H-D' \geq H-B > H-2B \geq 0$ and $(D')^2 > 0$, we have $D'.(H-D') >0$. 

Now define the set
\[ \Omega (D') = \{M \in \Pic S :  M \geq D', M \mbox{ is nef, satisfies
(i)-(ii) and } M.(H-M) \leq D'.(H-D') \}. \]
 
We have just seen that this set is nonempty.

Note that for any $M \in \Omega (D')$ we have $H-2M>0$, whence $H.M < \frac{1}{2}H^2$ 
is bounded. Let then $D_0$ be a {\it maximal} divisor in $\Omega (D')$, that is a divisor in
$\Omega (D')$ such that $H.D_0 \geq H.M$ for any $M \in \Omega (D')$. We want to show that
$h^1(H-2D_0)=0$.

Set $R=H-2D_0$. Assume, to get a contradiction that $h^1(H-2D_0)>0$. Since
$R^2 \geq 0$ it follows from Ramanujam vanishing \cite[Cor.II.12.3]{bpv} that 
$R +K_S$ is not 1-connected, whence $R+K_S \sim R_1+R_2$, for $R_1 >0$ and $R_2 >0$ with $R_1.R_2 \leq 0$.

We can assume that $R_1.H \leq R_2.H$. Define $D_1=D_0+R_1$. If $D_1$ is nef,
$\phi(D_1)$ is calculated by a nef divisor, whence $\phi(D_1) \geq \phi(D') \geq 2$ and 
$D_1^2 \geq D_0^2 \geq (D')^2 \geq 4$ (since $D_1 \geq D_0 \geq D'$). Moreover
\[ (H-2D_1)^2 = (R-2R_1)^2 = (R_2-R_1)^2 = R^2 -4R_1.R_2 \geq R^2 \geq 0, \]
and since 
\[ (H-2D_1).H=(R-2R_1).H = (R_2-R_1).H \geq 0, \]
we get by Riemann-Roch and the fact that $H$ is not $2$-divisible in $\Num S$, that
$H-2D_1 >0$.

Furthermore
\[ D_1.(H-D_1)= (D_0+R_1).(H-D_0-R_1) =  D_0.(H-D_0) + R_1.R_2 \leq D_0.(H-D_0), \]
whence $D_1$ is an element of $\Omega (D')$ with $H.D_1 > H.D_0$, contrary to our
assumption that $D_0$ is maximal.

Hence $D_1$ cannot be nef and there exists a nodal curve $\Gamma$ with $\Gamma.D_1
<0$ (whence $\Gamma.R_1 <0$). Since $H$ is ample we must have $\Gamma.(H-D_1) \geq -\Gamma.D_1 +1 \geq 2$. 
Let now $D_2 = D_1-\Gamma$. Since $\Gamma.R_1 <0$ we have $D_2 \geq D_0$, whence, if $D_2$ is
nef, we have as above that $\phi(D_2) \geq \phi(D') \geq 2$ and $D_2^2 \geq D_0^2 \geq
(D')^2 \geq 4$. Moreover $H-2D_2 > H-2D_1 >0$ and 
\[ (H-2D_2)^2 = (H-2D_1+2\Gamma)^2 =  (H-2D_1)^2 -8 +4(H-2D_1).\Gamma \geq (H-2D_1)^2+ 4 
>0. \]

Furthermore we also have
\begin{eqnarray*}
D_2.(H-D_2) & = &  (D_1-\Gamma).(H-D_1+\Gamma)=  D_1.(H-D_1) -\Gamma.(H-D_1)+\Gamma.D_1
+2 \leq \\ & \leq & D_1.(H-D_1) -1 < D_1.(H-D_1) \leq D_0.(H-D_0),
\end{eqnarray*}
whence $D_2 = D_0 + (R_1 - \Gamma) \not= D_0$. Now if $D_2$ is nef, then it is
an element of $\Omega (D')$ with $H.D_2 > H.D_0$, contrary to our assumption that
$D_0$ is maximal.

Hence $D_2$ cannot be nef, and we repeat the process by finding a nodal
curve $\Gamma_1$ such that $\Gamma_1 . (R_1 - \Gamma) < 0$ and so on. However, since
$R_1$ can only contain finitely many nodal curves, the process must end, that is
$h^1(H-2D_0)=0$, as we claimed.

Note that since $D_0^2 \geq (D')^2 = B^2$, then $D_0$ also satisfies (iii)
above, that is $H^2 \geq 64$ if $D_0^2 =4$ and $H^2 \geq 54$ if $D_0^2=6$. Furthermore
$D_0$ is base-point free since it is nef with $\phi(D_0) \geq \phi(D') \geq 2$ 
\cite[Prop.3.1.6, 3.1.4 and Thm.4.4.1]{cd}.

Now let $D$ be a general smooth curve in $|D_0|$. We have
\[ \deg (H-D_0)_{|D} = (H-D_0).D_0= D_0^2 + (H-2D_0).D_0 \geq 2g(D), \]
where we have used that $(H-2D_0).D_0 \geq \phi(D_0) \geq 2$ by Lemma \ref{lemma0}. 
Since $D$ is not hyperelliptic, it follows
that $(H-D_0)_{|D}$ is base-point free and birational, whence the map $\mu_{(H-D_0)_{|D}, \omega_D}$ is
surjective by \cite[Thm.1.6]{as}.

From the short exact sequence 
\[ 0 \hpil \O_S(H-2D_0) \hpil \O_S(H-D_0) \hpil \O_D(H-D_0) \hpil 0, \]
and the fact that $h^1(\O_D(H-D_0)) = 0$ for reasons of degree, we find $h^1(H-D_0)=0$. 

Now to show that $S$ is nonextendable, we only have left to show, by
Proposition \ref{mainextenr}, that the map 
$\Phi_{H_D, \omega_D}$ is surjective.

From $(H-2D_0).D_0 \geq 2$ again, we get $\deg H_D \geq 4g(D)-2$, whence 
by \cite[Thm.2]{bel}, the map $\Phi_{H_D, \omega_D}$ is surjective provided
that $\Cliff(D) \geq 2$. This is satisfied if $D_0^2 \geq 8$ by \cite[Cor.1 and Prop.4.15]{kl1}.

If $D_0^2=6$, then $g(D)=4$, whence by \cite[Prop.1.10]{wa}, the
map $\Phi_{H_D, \omega_D}$ is surjective if $h^0(\O_D(3D_0+K_S-H))=0$ (see also
Theorem \ref{tendian}(b) below). Since $H^2 \geq 54$, we get by the Hodge index
theorem that $H.D \geq 18$ with equality if and only if $H \eqv 3D_0$. If 
$H.D_0 >18$, we get $\deg \O_D(3D_0+K_S-H) <0$ and $\Phi_{H_D, \omega_D}$ is
surjective. If $H \eqv 3D_0$, then either $H \sim 3D_0$ and 
$h^0(\O_D(3D_0+K_S-H))=h^0(\O_D(K_S))=0$ or $H \sim 3D_0+K_S$ and then we
exchange $D_0$ with $D_0+K_S$ and we are done.

If $D_0^2=4$, then $g(D)=3$, whence by \cite[Prop.1.10]{wa}, the map
$\Phi_{H_D, \omega_D}$ is surjective if $h^0(\O_D(4D_0-H))=0$ (see also
Theorem \ref{tendian}(a) below). Since $H^2 \geq 64$, we get by the Hodge index
theorem that $H.D \geq 17$, whence $\deg \O_D(4D_0-H) <0$ and the map
$\Phi_{H_D, \omega_D}$ is surjective. 
\end{proof}

We recall the following result on Gaussian maps on curves on Enriques surfaces.

\begin{thm} \cite{klgm}
\label{tendian}
Let $S$ be an Enriques surface, let $L$ be a very ample line bundle on $S$
and let $D_0$ be a line bundle such that $D_0$ is nef, $D_0^2 \geq 4,
\phi(D_0) \geq 2$ and $H^1(D_0-L)=0$. Let $D$ be a general divisor in
$|D_0|$. Then the Gaussian map $\Phi_{L_{|D}, \omega_D}$ is surjective if
one of the hypotheses below is satisfied:
\begin{itemize}
\item[(a)] $D_0^2 = 4$ and $h^0(\O_D(4D_0-L)) = 0$;
\item[(b)] $D_0^2 =6$ and $h^0(\O_D(3D_0+K_S-L)) = 0$; 
\item[(c)] $D_0^2 \geq 8$ and $h^0(2D_0-L)=0$;
\item[(d)] $D_0^2 \geq 12$ and $h^0(2D_0-L) = 1$;
\item[(e)] $H^1(L_{|D})=0$, $L.D_0 \geq \frac{1}{2} D_0^2 + 2 \geq 6$ and $h^0(2D_0-L) \leq \Cliff(D) - 2$.
\end{itemize}
\end{thm}

We now get an improvement of Proposition \ref{ramextenr} in the cases
$B^2=6$ and $B^2=4$.

\begin{prop} 
\label{ramextenr6}
Let $S \subset \PP^r$ be an Enriques surface, denote by $H$ its hyperplane section and assume that
$H$ is not $2$-divisible in $\Num S$. Suppose there exists a line bundle $B$ on $S$
with the following properties:
\begin{itemize}
\item[(i)] $B^2 =6$ and $\phi(B) =2$, 
\item[(ii)] $(H-2B)^2 \geq 0$ and $H-2B \geq 0$,
\item[(iii)] $h^0(3B-H)=0$ or $h^0(3B+K_S-H)=0$.
\end{itemize} 
Then $S$ is nonextendable.
\end{prop}

\begin{proof}
Argue exactly as in the proof of Proposition \ref{ramextenr} and let $D'$, $D_0$ and $D$
be as in that proof, so that, in particular, $D_0^2 \geq (D')^2 = 6$. 
If $D_0^2 \geq 8$, we are done by Proposition \ref{ramextenr}. If
$D_0^2 = 6$ write $D_0 = D' + M$ with $M \geq 0$. Since both $D_0$ and $D'$ are nef we find
$6 = D_0^2 = (D')^2 + D'.M + D_0.M \geq 6$, whence $D'.M = D_0.M = 0$, so that $M^2 = 0$. Therefore 
$M = 0$ and $D_0=D'$, whence $3D_0-H \sim 3D'-H \leq 3B-H$. It follows that either
$h^0(3D_0-H)=0$ or $h^0(3D_0+K_S-H)=0$. Possibly after exchanging $D_0$ with $D_0+K_S$,
we can assume that $h^0(3D_0+K_S-H)=0$. From the exact sequence 
\[ 0 \hpil \O_S(2D_0+K_S-H) \hpil \O_S(3D_0+K_S-H) \hpil \O_D(3D_0+K_S-H) \hpil 0, \]
and the fact that $h^1(2D_0+K_S-H)=h^1(H-2D_0)=0$, we get
$h^0(\O_D(3D_0+K_S-H))=0$, whence the map $\Phi_{H_D, \omega_D}$ is surjective by Theorem \ref{tendian}(b).

The multiplication map $\mu_{V_D, \omega_D}$ is surjective as in the proof of 
Proposition \ref{ramextenr}, whence $S$ is nonextendable by Proposition \ref{mainextenr}.
\end{proof}

\begin{prop} 
\label{ramextenr4}
Let $S \subset \PP^r$ be an Enriques surface, denote by $H$ its hyperplane section and assume that
$H$ is not $2$-divisible in $\Num S$. Suppose there exists a line bundle $B$ on $S$
with the following properties:
\begin{itemize}
\item[(i)] $B$ is nef, $B^2 =4$ and $\phi(B) =2$, 
\item[(ii)] $(H-2B)^2 \geq 0$ and $H-2B \geq 0$,
\item[(iii)] $H.B > 16$.
\end{itemize}  
Then $S$ is nonextendable.
\end{prop}

\begin{proof}
Argue exactly as in the proof of Proposition \ref{ramextenr} and let $D'$,
$D_0$ and $D$ be as in that proof. Since $B$ is assumed to be nef, we have
$D'=B$, and since $D_0 \geq D'$, we get $H.D_0 > 16$. If $D_0 ^2 \geq 8$, we are done by
Proposition \ref{ramextenr}. If $D_0 ^2 =6$, then we must have $D_0 > D'= B$, so that $H.D_0 \geq
18$ whence $(3D_0-H).D_0 \leq 0$. This gives that if $3D_0-H>0$, then it is a nodal cycle, whence either 
$h^0(3D_0-H)=0$ or $h^0(3D_0+K_S-H)=0$. Now we are done by Proposition \ref{ramextenr6}. 

If $D_0 ^2 =4$, then, as in the proof of Proposition \ref{ramextenr6}, $D_0=D'=B$, whence the map 
$\Phi_{H_D, \omega_D}$ is surjective by Theorem \ref{tendian}(a), since $\deg \O_D(4D_0-H) <0$.
The multiplication map $\mu_{V_D, \omega_D}$ is surjective as in the proof of 
Proposition \ref{ramextenr}, whence $S$ is nonextendable by Proposition \ref{mainextenr}.
\end{proof}

In several cases the following will be very useful:

\begin{lemma} 
\label{subram}
Let $S \subset \PP^r$ be an Enriques surface with hyperplane section $H \sim 2B+A$, for $B$ nef, $B^2
\geq 2$, $A^2 = 0$, $A > 0$ primitive, $H^2 \geq 28$ and satisfying one of the following conditions:
\begin{itemize}
\item[(i)] $A$ is quasi-nef and $(B^2, A.B) \not \in \{ (4,3), (6,2) \}$;
\item[(ii)] $\phi(B) \geq 2$ and $(B^2, A.B) \not \in \{ (4,3), (6,2) \}$;
\item[(iii)] $\phi(B) = 1$, $B^2 =2l$, $B \sim lF_1 + F_2$, $l \geq 1$, $F_i > 0$, $F_i^2 = 0$, 
$i = 1, 2$, $F_1.F_2 = 1$, and either
\begin{itemize}
\item[(a)] $l \geq 2$, $F_i.A \leq 3$ for $i=1,2$ and $(l,F_1.A,F_2.A) \neq (2,1,1)$; or 
\item[(b)] $l=1$, $5 \leq B.A \leq 8$, $F_i.A \geq 2$ for $i = 1, 2$
and $(\phi(H),F_1.A,F_2.A) \neq (6,4,4)$.
\end{itemize}
\end{itemize}
Then $S$ is nonextendable.
\end{lemma}

\begin{proof}
Note that possibly after replacing $B$ with $B + K_S$ if $B^2 = 2$ we can, without
loss of generality, assume that $B$ is base-component free.

We first prove the lemma under hypothesis (i).

We have that $B + A$ is nef, since any nodal curve $\Gamma$ with $\Gamma.(B + A) < 0$ would
have to satisfy $\Gamma.A = -1$ and $\Gamma.B = 0$, whence $\Gamma.H = -1$, a contradiction.

Now let $D_0 = B + A$. Then $D_0^2=B^2+2B.A \geq 6$, since $A.B \geq 2$ for $2A.B = A.H 
\geq \phi(H) \geq 3$, and $\phi(D_0) \geq \phi(B) \geq 1$.

If $\phi(D_0)=1=F.D_0$ for some $F>0$ with $F^2=0$ we get $F.B=1$, $F.A=0$ giving the
contradiction $F.H=2$. Therefore $\phi(D_0) \geq 2$.

One easily checks that (i) implies $D_0^2 \geq 12$.  Since $h^0(2D_0-H) =h^0(A)=1$
by Theorem \ref{lemma:qnef}, we have that $\Phi_{H_D, \omega_D}$ is surjective by Theorem
\ref{tendian}(d).

Also $h^1(H-2D_0)=h^1(-A)=0$, again by Theorem \ref{lemma:qnef}, so that we have $V_D=H^0(\O_D(H-D_0))$.
As $H - D_0 = B$ is base-component free and $|D_0|$ is base-point free and birational by
\cite[Lemma4.6.2, Thm.4.6.3 and Prop.4.7.1]{cd}, also $V_D$ is base-point free and is either a complete pencil
or birational, and then $\mu_{V_D, \omega_D}$ is surjective by the base-point free pencil trick and by 
\cite[Thm.1.6]{as} (see also \eqref{biraz}). Then $S$ is nonextendable by
Proposition \ref{mainextenr}.

Therefore the lemma is proved under the assumption (i) and, in particular, the whole lemma
is proved with the additional assumption that $A$ is quasi-nef.

Now assume that $A$ is not quasi-nef. Then there is a $\Delta >0$ with $\Delta^2=-2$ and $\Delta.A \leq
-2$. We have $\Delta.B \geq 2$ by the ampleness of $H$. Furthermore, among all such $\Delta$'s we will
choose a minimal one, that is such that no $0 < \Delta' < \Delta$ satisfies
$(\Delta')^2=-2$ and $\Delta'.A \leq -2$.

We now claim that $B_0:=B+\Delta$ is nef. Indeed, if there is a nodal curve
$\Gamma$ with $\Gamma.(B+\Delta) <0$ then $\Gamma.\Delta <0$ and we must have
$\Delta_1 := \Delta-\Gamma >0$ with $\Delta_1^2=-4-2\Delta.\Gamma$.

If $\Delta.\Gamma \leq -2$ then $\Delta_1^2 \geq 0$ whence $\Delta_1.A \geq 0$ by Lemma \ref{nefred}
and $\Gamma.A \leq \Delta.A \leq -2$, contradicting the minimality of $\Delta$. Therefore
$\Delta.\Gamma=-1$, $\Gamma.B=0$ and $\Delta_1^2=-2$. The ampleness of $H$ implies $\Gamma.A >0$,
whence $\Delta_1.A <\Delta.A \leq -2$, again contradicting the minimality of $\Delta$.

Therefore $B_0:=B+\Delta$ is nef with $B_0^2 \geq 2 + B^2$, and, as $\phi(B_0)$ is computed by a nef
isotropic divisor, we have that $\phi(B_0) \geq \phi(B)$.

We also note that $H-2B_0 \sim A-2\Delta >0$ and primitive by Lemma \ref{A} with $(H-2B_0)^2 \geq 0$.

Under the assumptions (ii), we have $\phi(B_0) \geq 2$. Then $S$ is nonextendable by Proposition
\ref{ramextenr} if $B_0^2 \geq 8$. If $B_0^2=6$, we have $B^2=4$ and $\Delta.B=2$, so that
$\Delta.A =-2$ or $-3$ by the ampleness of $H$. Hence $H \sim 2B_0+ A'$, with $B_0^2=6$ and $A' \sim
A-2\Delta$ satisfies $(A')^2=0$ or $4$. In the first case we are done by conditions (i) if $A'$ is
quasi-nef (because $B_0.A' = (B+\Delta).A' \geq 4$), and if not we can just repeat the process and find
that $S$ is nonextendable by Proposition \ref{ramextenr} (because we find a divisor $B_0'$ with
$(B_0')^2 \geq 8$). In the case $(A')^2=4$ we have $A'.B_0 \geq 5$ by the Hodge index theorem.
Therefore $(3B_0-H).B_0 = (B_0-A').B_0 \leq 1 < \phi(B_0)$, so that if $3B_0-H>0$ it is a nodal cycle.
Hence either $h^0(3B_0-H)=0$ or $h^0(3B_0+K_S-H)=0$ and $S$ is nonextendable by Proposition
\ref{ramextenr6}.  

We have therefore shown that $S$ is nonextendable under conditions (ii).

Now assume (iii) and, using Lemma \ref{A}, write $A \sim A_0 + k\Delta$
with $A_0 > 0$ primitive, $A_0^2 = 0$ and $k:=-\Delta.A = \Delta.A_0 \geq 2$. 

As $F_i.A = F_i.A_0 + kF_i.\Delta$, the primitivity of $F_i$, $\Delta.B \geq 2$ and the hypotheses in
case (iii-a) yield the only possibility $k=2$, $F_1.\Delta=F_1.A_0=1$. Then $H
\sim 2B_0 +A_0$ with $B_0^2 \geq 6$, $\phi(B_0) \geq 2$, $B_0$ nef, $A_0^2=0$ and
$B_0.A_0=(B+\Delta).A_0 \geq 3$, so that conditions (ii) are satisfied and $S$ is nonextendable.

Finally we assume we are in case (iii-b), so that $F_i.A \leq 6$ for $i=1,2$ by hypothesis.

Suppose $\Delta.F_1 \leq 0$. Then $F_2.\Delta \geq 2$. As $6 \geq F_2.A = F_2.A_0 + kF_2.\Delta$, we
must have $k=F_2.\Delta=2$, so that  $\Delta.F_1=0$ and $4 \leq F_2.A \leq 6$. In particular, 
$F_1.B_0=1$, so that $B_0 \sim 2F_1+ F_2'$, where $F_2' \sim F_2+\Delta-F_1>0$ satisfies
$(F_2')^2=0$. We have $F_1.A_0=F_1.A \leq 4$, and equality implies $F_2.A=4$, whence $F_2 \eqv A_0$,
so that $F_1.A_0=F_1.F_2=1$, a contradiction. Hence $F_1.A_0 \leq 3$. Moreover $F_2'.A_0=
(F_2+\Delta-F_1).A_0 =(F_2-F_1).A -2 \leq 2$, as $F_2.A \leq 6$. Also it cannot be
$(F_1.A_0,F_2'.A_0)=(1,1)$, for then $F_1.A=1$. Therefore $H \sim 2B_0 +A_0$ satisfies the conditions
in (iii-a), so that $S$ is nonextendable.

We can therefore assume $\Delta.F_1 >0$, and by symmetry, also $\Delta.F_2 >0$. Hence $\phi(B_0) \geq 2$.

If $k \geq 3$, then $F_i.A = F_i.A_0 + k F_i.\Delta \geq 4$ for $i=1,2$, so that by our assumptions 
we can only have $k=3$, $F_i.A=4$ and $F_i.\Delta = F_i.A_0=1$. Then $B.A=8$ and $H^2=40$, so that
$\phi(H) \leq 5$ by hypothesis. Pick any isotropic divisor $F>0$ satisfying $F.H=\phi(H)$. Since
$(A')^2=4$ we have $5 \geq F.H = 2F.B_0 + F.A' \geq 5$, so that $F.H = 5$, $F.A'=1$, $(A'-2F)^2=0$, $A'-2F >0$,
and $(A'-2F).H =(A-2\Delta-2F).H=4$, a contradiction.

Therefore $k=2$, so that $A_0^2=0$. As $B_0.A_0=(B+\Delta).A_0 =B.A_0+2 \geq 3$, we see that the
conditions (ii) are satisfied, unless possibly if $B_0^2=4$ and $B.A_0=1$. In this case
$B.\Delta=2$ and $A_0 \eqv F_i$, for $i=1$ or $2$. Hence $\Delta.B = \Delta.(F_1+F_2) = 3$, a
contradiction. Therefore the conditions (ii) are satisfied and $S$ is nonextendable.
\end{proof}

We also have the following helpful tools to check the surjectivity of 
$\mu_{V_D, \omega_D}$ in the cases where $h^1(H-2D_0) \neq 0$. The first lemma holds on any
smooth surface.

\begin{lemma} 
\label{multhelp}
Let $S$ be a smooth surface, let $L$ be a line bundle on $S$ and let $D_1$ and 
$D_2$ be two effective nonzero divisors on $S$ not intersecting the base locus of $|L|$ and such that
$h^0(\O_{D_1})=1$ and $h^0(\O_{D_1}(-L)) = h^0(\O_{D_2}(-D_1)) = 0$.
For any divisor $B>0$ on $S$ define $V_B = \Im \{H^0 (S,L) \khpil
H^0(B,\O_B(L))\}$. 
If the multiplication maps $\mu_{V_{D_1}, \omega_{D_1}}$ and $\mu_{V_{D_2}, \omega_{D_2}(D_1)}$ are
surjective then $\mu_{V_D, \omega_D}$ is surjective for general $D \in |D_1 + D_2|$.
\end{lemma}

\begin{proof} Let $D'= D_1 + D_2$. By hypothesis we have the exact sequence
\[ 0 \hpil H^0(\omega_{D_1}) \hpil H^0(\omega_{D'}) \mapright{\psi}
H^0(\omega_{D_2}(D_1)) \hpil 0. \] 
 
Moreover by definition we have two surjective maps $\pi_i : V_{D'} \to V_{D_i}, i = 1, 2$, whence a commutative 
diagram
\[  
\xymatrix{
0  \ar[r] & W \ar[r] \ar[d]^{\varphi} &
V_{D'} \* H^0(\omega_{D'})  \ar[r]^{ \hskip -.7cm \pi_2 \* \psi} \ar[d]^{\mu_{V_{D'},
\omega_{D'}}}   &  V_{D_2} \* H^0(\omega_{D_2}(D_1)) \ar[r] \ar[d]^{\mu_{V_{D_2},
\omega_{D_2}(D_1)}}  & 0 
\\
0 \ar[r] & H^0(\omega_{D_1}(L)) \ar[r]^{\chi} & H^0(\omega_{D'}(L)) \ar[r] &
H^0(\omega_{D_2}(D_1+L))  \ar[r] & 0
}
\]
where $W := \Ker \pi_2 \* H^0(\omega_{D'}) + V_{D'} \* \Ker \psi$ and $\varphi$ is just
the restriction of $\mu_{V_{D'}, \omega_{D'}}$ to this subspace. Since $\mu_{V_{D_2},
\omega_{D_2}(D_1)}$ is surjective, to conclude the surjectivity of $\mu_{V_{D'},
\omega_{D'}}$ we just show the  surjectivity of
$\varphi$. Now the commutative diagram
\[  
\xymatrix{
V_{D'} \*  H^0(\omega_{D_1}) \ar[r]^{\cong} \ar[rd]^{\pi_1 \* \Id} & 
V_{D'} \* \Ker \psi \ar[r] & H^0(\omega_{D'}(L))
 \\
& V_{D_1} \* H^0(\omega_{D_1}) \ \ \ar[r]^{\mu_{V_{D_1}, \omega_{D_1}}} 
\ar[r]^{\mu_{V_{D_1}, \omega_{D_1}}}  & \ H^0(\omega_{D_1}(L)) \ar[u]^{\chi}
}
\]
and the injectivity of $\chi$ show that $H^0(\omega_{D_1}(L))=\Im \mu_{V_{D_1}, \omega_{D_1}} = \Im
\varphi_{|_{V_{D'} \* \Ker \psi}}$, as required.

Therefore $\mu_{V_{D'}, \omega_{D'}}$ is surjective. Now for any divisor $B \in |D'|$, let
$M_B = \Ker \{ V_B \* \O_B \to L_{|B} \}$. By hypothesis we have that
$V_{D'}$ globally generates $L_{|D'}$, whence we have an exact sequence
\[ 0 \hpil M_{D'} \* \omega_{D'} \hpil V_{D'} \* \omega_{D'} \to L_{|D'} \* 
\omega_{D'} \hpil 0. \] 

Since $h^1(\omega_{D_1}(L))=h^1(\omega_{D_2}(D_1+L))=0$ we get
$h^1(L_{|D'} \* \omega_{D'})=0$ and the surjectivity of $\mu_{V_{D'},\omega_{D'}}$ implies 
that $h^1(M_{D'} \* \omega_{D'}) \leq \dim V_{D'} = \dim V_D$. By semicontinuity the same holds
for a general $D \in  |D'|$, whence $\mu_{V_D, \omega_D}$ is surjective. 
\end{proof}

\begin{lemma} 
\label{2E}
Let $S$ be an Enriques surface, let $L$ be a very ample line bundle on $S$ and let $D_0$ be a nef and big 
divisor on $S$ such that $\phi(D_0) \geq 2$. Let $E > 0$ be such that $E^2 = 0$, $E.L = \phi(L)$ and define, on
a general $D \in |D_0|$, $V_D = \Im \{H^0(\O_S(L-D_0)) \khpil H^0(\O_D(L-D_0))\}$.

If $|L-D_0-2E|$ is base-component free, $h^1(D_0+K_S-2E) = h^2(D_0+K_S-4E) = 0$ and 
\begin{equation}
\label{green}
h^0(L-2D_0-2E) + h^0(\O_D(L-D_0-4E)) \leq \frac{1}{2}(L-D_0-2E)^2 -1
\end{equation}
then $\mu_{V_D, \omega_D}$ is surjective.
\end{lemma}

\begin{proof}
Consider the natural restriction maps $p_D: H^0(\O_S(L-D_0)) \khpil H^0(\O_{D}(L-D_0))$, 
$p'_D: H^0(\O_S(D_0+K_S)) \khpil H^0(\omega_D)$, $r_D: H^0(\O_S(L-D_0-2E)) \khpil H^0(\O_D(L-D_0-2E))$, $r'_D:
H^0(\O_S(2E+D_0+K_S)) \khpil H^0(\omega_D(2E))$. Then $V_D = \Im p_D$, $W_D := \Im r_D$ and let $\mu =
\mu_{2E,D_0+K_S}$, $\mu' = \mu_{2E,L-D_0-2E}$ be the multiplication maps of line bundles on $S$. We have a
commutative diagram
\begin{equation*}
\xymatrix{
H^0(2E) \* H^0(L-D_0-2E) \*  H^0(D_0+K_S) \ar[r]^{ \Id \* \mu}
 \ar[d]_{\mu' \* \Id} & H^0(L-D_0-2E) \* H^0(2E+D_0+K_S) 
\ar[d]^{r_D \* r'_D}
\\
H^0(L-D_0)  \* H^0(D_0+K_S) \ar[d]_{p_D \* p'_D}  & W_D \*
H^0(\omega_D(2E)) \ar[d]^{\mu_{W_D, \omega_D(2E)}}
\\
 V_D  \* H^0(\omega_D)  \ar[r]^{\mu_{V_D,\omega_D}}  & H^0(\O_D(L+K_S)).  
}
\end{equation*} 
Since $H^1(D_0+K_S-2E) = H^2(D_0+K_S-4E) = 0$ we have that $\mu$ is surjective by Castelnuovo-Mumford's lemma. 
At the same time, since $h^1(2E+K_S)=0$, we have that $r'_D$ is surjective. To conclude the surjectivity of
$\mu_{V_D, \omega_D}$, by the above diagram, we just need to prove that $\mu_{W_D, \omega_D(2E)}$ is surjective.
To see the latter note that, as $D$ is general, $W_D$ is base-point free and $\dim W_D - 2 = h^0(L-D_0-2E) -
h^0(L-2D_0-2E) - 2$, whence $\mu_{W_D, \omega_D(2E)}$ is surjective by the $H^0$-lemma \cite[Thm.4.e.1]{gr} as 
soon as $h^1(\omega_D(2E-(L-D_0-2E)) \leq h^0(L-D_0-2E) - h^0(L-2D_0-2E) - 2$, which is equivalent to
(\ref{green}) by Riemann-Roch on $S$ and Serre duality on $D$.
\end{proof}

\section{Strategy of the proof of Theorem \ref{main}}
\label{class}

In this section we prove Theorem \ref{main} for all very ample line bundles on an Enriques surface
except for some concrete cases, and then we give the main strategy of the proof in these remaining 
cases, which will then be carried out in Sections \ref{caseD}-\ref{caseS}. We also set some notation 
and conventions that will be used throughout the proofs, often without further mentioning. 

Let $S \subset \PP^r$ be an Enriques surface of sectional genus $g$ and let $H$ be its hyperplane bundle.
As we will prove a result also for $g = 15$ and $17$ (Proposition \ref{precisa}) we will henceforth assume
$g \geq 17$ or $g = 15$, so that $H^2=2g-2 \geq 32$ or $H^2 = 28$, and, as $H$ is very ample, $\phi(H) \geq 3$.
We choose a genus one pencil $|2E|$ such that $E.H = \phi(H)$ and, as $H$ is not of small type by Lemma
\ref{STlemma2}, we define $\alpha:= \alpha_E (H)$ and $L_1 := H - \alpha E$, where $\alpha_E (H)$ is as in
(\ref{alpha1}) and  (\ref{alpha2}). By Lemmas \ref{nefred} and \ref{STlemma3} we have that $L_1 > 0$ and 
$L_1^2>0$. Now suppose that $L_1$ is not of small type. Starting with $L_0:=H$ and $E_0:=E$ we continue the
process inductively until we reach a line bundle of small type, as follows. Suppose given, for $i \geq 1$,
$L_i>0$ not of small type with $L_i^2 > 0$. We choose $E_i > 0$ such that $E_i^2 = 0$, $E_i.E_{i-1} > 0$,
$E_i.L_i = \phi(L_i)$ and define $\alpha_i = \alpha_{E_i}(L_i)$ and $L_{i+1} = L_i-\alpha_iE_i$. Note that
$L_{i+1}>0$ by Lemma \ref{nefred}. Now if $L_{i+1}^2 = 0$ we write $L_{i+1} \eqv \alpha_{i+1}E_{i+1}$ and define
$L_{i+2} = 0$, which is of small type by definition and we also have $E_{i+1}.E_i > 0$ because $L_i^2 > 0$. If
$L_{i+1}^2 > 0$ then either $L_{i+1}$ is of small type or we can continue.

We then get, for some integer $n \geq 1$:
\begin{equation} 
\label{decomposition}  
H = \alpha E + \alpha_1 E_1 + \ldots + \alpha_{n-1} E_{n-1} + L_n,
\end{equation}
with $\alpha \geq 2$, $\alpha_i \geq 2$ for $1 \leq i \leq n-1$ and $L_n$ is of small type.

Moreover $E.E_1 \geq 1$, $E_i.E_{i+1} \geq 1$, $E$ and $E_i$ are primitive for all $i$, $L_i^2 > 0$ and $E_i.L_i
= \phi(L_i)$ for $0 \leq i \leq n-2$ and $L_{n-1}^2 \geq 0$.

We record for later the following fact, which follows immediately from the definitions:
\begin{equation}
\label{eq:class.1}   
E_1.(H - \alpha E) \leq \phi(H) \mbox{ and if } \alpha \geq 3 \mbox{ then } 
E_1.(H - \alpha E) \geq \phi(H) + 1 - E.E_1. 
\end{equation}  

Furthermore we claim that $\alpha_i = 2$ for $1 \leq i \leq n-1$. To see this we proceed by induction on $i$.
If $(L_1 - 2E_1)^2 = 0$ then $\alpha_1 = 2$ by definition. On the other hand if $(L_1 - 2E_1)^2 > 0$ to see that
$\alpha_1 = 2$ we just need to prove that $E_0.(L_1 - 2E_1) \leq \phi(L_1)$, or, equivalently, that
$\phi(L_0) \leq E_1.L_0 + (2-\alpha_0) E_1.E_0$. Now the latter holds both when $\alpha_0 = 2$ and, by
\eqref{eq:class.1}, when $\alpha_0 \geq 3$. By induction and the same proof for $i=1$ we can prove that 
$\alpha_i = 2$ for $1 \leq i \leq n-2$ and also for $i=n-1$ if $L_{n-1}^2 > 0$. Finally if $L_{n-1}^2 = 0$ we
have $L_{n-2} \eqv 2E_{n-2} + \alpha_{n-1}E_{n-1}$, whence $(\alpha_{n-1}E_{n-2}.E_{n-1})^2 = \phi(L_{n-2})^2
\leq L_{n-2}^2 = 4\alpha_{n-1}E_{n-2}.E_{n-1}$. Therefore $\alpha_{n-1}E_{n-2}.E_{n-1} \leq 4$ and if
$\alpha_{n-1} \geq 3$ we get $E_{n-2}.E_{n-1}=1$, giving the contradiction $\alpha_{n-1} = \phi(L_{n-2}) \leq
E_{n-1}.L_{n-2}=2$ and the claim is proved. 

We now search for a divisor $B$ as in Proposition \ref{ramextenr} to show
that $S \subset \PP^r$ is nonextendable.

Assume for the moment that $H$ is not 2-divisible in $\Num S$ and that
$n \geq 2$ (that is $L_1$ is not of small type).

If $n \geq 4$, then set $B=E+E_1+E_2+E_3$. We have $B^2 \geq 6$ with equality
if and only if $E.E_2 = E.E_3 = E_1.E_3 = 0$. But the latter implies the contradiction
$E_2 \eqv E \eqv E_3$. Hence $B^2 \geq 8$ and clearly $\phi(B) \geq 2$. Now
\[H - 2B = (\alpha -2)E + 2 \sum_{i=4}^{n-1} E_i + L_n \geq 0, \] 
where the sum is empty if $n = 4$. Hence $(H-2B)^2 \geq 0$,
therefore $B$ satisfies the conditions in Proposition \ref{ramextenr} and $S$ is nonextendable.

If $n =3$, then $H = \alpha E + 2 E_1 + 2 E_2 + L_3$. Set 
$B=\lfloor \frac{\alpha}{2} \rfloor E + E_1 +  E_2$.  Then $B$ satisfies the conditions in
Proposition \ref{ramextenr}, whence $S$ is nonextendable, unless
\begin{itemize}
\item[(I-A)] $n = 3, E_2 \eqv E$, $E.E_1=1$.
\item[(II)] $n = 3, E.E_1=E.E_2=E_1.E_2=1$, $\alpha \in \{2,3 \}$, $H^2 \leq 52$.
\end{itemize}

If $n = 2$, then $H = \alpha E + 2 E_1 + L_2$. Set $B=\lfloor
\frac{\alpha}{2} \rfloor E + E_1$.  Then $B$ satisfies the conditions in Proposition \ref{ramextenr}, 
whence $S$ is nonextendable, unless
\begin{itemize}
\item[(I-B)] $n = 2, E.E_1=1$.
\item[(III)] $n = 2, E.E_1=2$, $\alpha \in \{2,3 \}$, $H^2 \leq 62$,
\end{itemize}
or $n = 2, E.E_1=3$, $\alpha \in \{2,3 \}$ and $H^2 \leq 52$. But the latter case does not occur. Indeed then
$6 \leq 6 + E.L_2 = E.H = \phi(H) \leq 6$ by \cite[Prop.1]{kl1}, whence $E.L_2 = 0$, therefore either $L_2=0$ or
$L_2 \eqv E$. Now since $(E + E_1)^2 = 6$ and $\phi(E + E_1)=2$, as $E$ and $E_1$ are primitive, we can write 
$E + E_1 \sim A_1+A_2+A_3$ with $A_i>0$, $A_i^2=0$. Therefore $18 \geq 6 + 3 \alpha + E_1.L_2 = (E + E_1).H \geq
3\phi(H) = 18$, whence $\alpha = 3$ and $E_1.H=12$. But then $E_1.(H-2E) = 6$ so that $\alpha = 2$, a
contradiction.

Now $L_n \geq 0$ and  $L_n^2 \geq 0$ so that, if $L_n>0$, it has (several) arithmetic genus $1$
decompositions. We want to extract from them any divisors numerically equivalent
to $E$ or to $E_1$, if possible. If, for example, we give priority to $E$, we will 
write $L_n \eqv E + L_n'$ and then, if $L_n'$ has an arithmetic genus $1$ decomposition with
$E_1$ present, we write $L_n' \eqv E_1 + M_n$. In case the priority is
given to $E_1$ we do it first with $E_1$ and then with $E$. Finally, for a reason that will 
be clear below, in the case (I-A), where only $M_3$ is defined, we will set $M_2 = M_3$.

To avoid treating the same case more times we make the following choice of
{\bf ``removing conventions''}:
\label{remconv}

\begin{itemize}
\item[(I-A)] Remove $E$ and $E_1$ from $L_3$, the one with lowest intersection
number with $L_3$ first, giving priority to $E_1$ in case $E.L_3 = E_1.L_3$.
\item[(I-B)] Remove $E$ and $E_1$ from $L_2$, the one with lowest intersection
number with $L_2$ first, giving priority to $E$ in case $E.L_2 = E_1.L_2$.
\item[(II)] Remove $E$, $E_1$ and $E_2$ from $L_3$, the one with lowest intersection
number with $L_3$ first, giving priority to $E$ first and then to $E_2$.
\item[(III)] Remove $E$ and $E_1$ from $L_2$, the one with lowest intersection
number with $L_2$ first, giving priority to $E$ in case $E.L_2 = E_1.L_2$.
\end{itemize}

At the end we get the following cases where the extendability of $S$ still has to be
checked, where $\gamma, \delta \in \{2,3 \}$:

\begin{itemize}
\item[(I)] $H \eqv \beta E + \gamma E_1 + M_2$, $E.E_1=1$, $H^2 \geq 32$ or $H^2 = 28$, 
\item[(II)] $H \eqv \beta E + \gamma E_1 + \delta E_2 + M_3$,
$E.E_1=E.E_2=E_1.E_2=1$, $\beta \in \{2,3 \}$,
$32 \leq H^2 \leq 52$ or $H^2 = 28$.
\item[(III)] $H \eqv \beta E + \gamma E_1 + M_2$, $E.E_1=2$,
$\beta \in \{2,3 \}$, $32 \leq H^2 \leq 62$ or $H^2 = 28$
\end{itemize}
(where the limitations on $\beta$ are obtained using the same $B$'s as above), in addition to:
\begin{itemize}
\item[(D)] $H \eqv 2H_1$ for some $H_1 >0$, $H_1^2 \geq 8$,
\item[(S)] $L_1$ is of small type and $H^2 \geq 32$ or $H^2 = 28$.  
\end{itemize}

\begin{defn}
We call such decompositions as in {\rm (I), (II)} and {\rm (III)}, obtained by the 
inductive process and removing conventions above, a {\bf ladder decomposition} of $H$. 
\end{defn}

Note that $M_n, n = 2, 3$, satisfies: $M_n \geq 0$, $M_n^2 \geq 0$ and $M_n$ is
of small type. Moreover, when $M_n > 0$, we will replace $M_n$ with $M_n + K_S$ that
has the same properties, to avoid to study the two different numerical equivalence
classes for $H$. Also note that $\beta \geq \alpha \geq 2$ and $\beta \geq \alpha + 2$ in (I-A).

We will treat all these cases separately in the next sections.

A useful fact will be the following

\begin{lemma} 
\label{nef:main}
If $E.E_1 \leq 2$, then $E+E_1$ is nef.
\end {lemma}

\begin{proof}
Let $\Gamma$ be a nodal curve such that
$\Gamma.(E+E_1) <0$. Since $E$ is nef, we must have $\Gamma.E_1 <0$. By Lemma 
\ref{A} we can then write $E_1=A+k\Gamma$, for $A >0$ primitive with $A^2=0$,
$k=-\Gamma.E_1 \geq 1$ and $\Gamma.A =k$. Since $A.L_1 \geq \phi(L_1) = E_1.L_1$, we get $k\Gamma.L_1 = (E_1 -
A).L_1 \leq 0$, whence $\Gamma.E >0$, because $H$ is ample. This yields $k \geq \Gamma.E+1 \geq 2$.
Hence $E.E_1= E.A + k\Gamma.E \geq 2\Gamma.E$, and we get $k=2$,
$\Gamma.E=1$ and $E.A=0$. Hence $A \eqv E$ by Lemma \ref{nefred}, contradicting the fact that $\Gamma.A =2$.
\end{proof}

From \cite[Prop.3.1.6, 3.1.4 and Thm.4.4.1]{cd} and the above Lemma we now know that $E+E_1$
is base-point free when $E.E_1 =2$, and that $E+E_1$ is base-component free when
$E.E_1=1$, unless $E_1 \sim E+ R$, for a nodal curve $R$ such that $E.R = 1$. But
since we are free to choose between $E_1$ and $E_1+K_S$ (they both calculate
$\phi(L_1)$), {\bf we adopt the convention of always choosing $E_1$ such that $E+E_1$ is
base-component free}. We therefore have

\begin{lemma} 
\label{nef:main2}
If $E.E_1 =2$, then $E+E_1$ is base-point free.

If $E.E_1 =1$, then $E+E_1$ is base-component free. Furthermore
if there exists $\Delta > 0$ such that $\Delta^2 = -2$ and $\Delta.E_1 < 0$, then
$\Delta$ is a nodal curve and $E_1 \sim E+\Delta+K_S$.

Moreover in both cases we have $H^1(E_1) = H^1(E_1 + K_S) = 0$.
\end {lemma}

\begin{proof} We need to prove the last two assertions. Suppose there exists a $\Delta > 0$ such that
$\Delta^2 = -2$ and $\Delta.E_1 < 0$. By Lemma \ref{A} we can write $E_1 = A + k \Delta$, for $A >0$
primitive with $A^2 = 0$ and $k = -\Delta.E_1 = \Delta.A \geq 1$. Now $0 \leq (E+E_1).\Delta =
E.\Delta -k$ gives $E.\Delta \geq k$. From $2 \geq E.E_1 = E.A + kE.\Delta \geq k^2$
we get $k = 1$, whence $E_1$ is quasi-nef and primitive, so that the desired vanishing
follows by Theorem \ref{lemma:qnef}. Now if $E.E_1 = 1$ we get $E.A = 0$, whence $A \eqv E$ by 
Lemma \ref{nefred} and then $E_1 \eqv E + \Delta$. Since $E_1$ is not nef, by
\cite[Prop.3.1.4, Prop.3.6.1 and Cor.3.1.4]{cd}  there is a nodal curve $R$ such that $E_1 \sim E + R +
K_S$, whence $\Delta = R$.
\end{proof}

Another useful nefness lemma is the following. 

\begin{lemma} 
\label{lemmone}
Let $H \sim \beta E + \gamma E_1 +M_2$ be of type {\rm (I)} or {\rm (III)}, with $M_2 >0$ and $M_2^2 \leq
4$. Let $i = 2$ and $M_2 \sim E_2$ or $i = 2, 3$ and $M_2 \sim E_2+E_3$ be genus $1$ decompositions
of $M_2$ (note that, by construction, $E.E_j \geq 1$ for $j = 1, 2$).

Assume that $E_i$ is quasi-nef. Then:
\begin{itemize}
\item[(a)] $|2E+E_1+E_i|$ is base-point free.
\item[(b)] $|E+E_1+E_i|$ is base-point free if $\beta=2$ or if $E.E_1=1$ and $E_1.E_i \neq E.E_i-1$.
\item[(c)] Assume $\gamma=2$ and $E.E_1=E_1.E_i=1$. Then $E+E_i$ is nef if either $E.E_i \geq 2$ or if 
$M_2^2 \geq 2$ and $E_1.M_2 \geq 4$. 
\item[(d)] Assume $\gamma=2$, $M_2^2=2$, $E.E_1=E_1.E_2=E_1.E_3=1$ and that both $E_2$ and $E_3$ are
quasi-nef. Then either $E+E_2$ or $E+E_3$ is nef.
\item[(e)] If $E.E_1=E.E_i=1$ and $E_1.E_i \neq 1$ then $E_1+E_i$ is nef.
\end{itemize}
\end{lemma}

\begin{proof}
Assume there is a nodal curve $R$ with $R.(E+E_1+E_i) < 0$. Then Lemma \ref{nef:main} and the
quasi-nefness of $E_i$ yield $R.(E+E_1)=0$ and $R.E_i=-1$. Moreover $R.E_1 \leq 0$ by the nefness of
$E$.

If $R.E=0$ or if $\beta=2$, then $R.H=R.(\beta E+\gamma E_1+M_2) \leq R.M_2$ implies $M_2^2 \geq 2$ 
and $R.E_j \geq 2$ for $j \in \{2,3\}-\{i\}$. By Lemma \ref{A} we can write $E_i \sim A + R$ with $A > 0$
primitive, $A^2=0$ and $A.R=1$. From $2 \geq E_j.E_i=A.E_j+R.E_j \geq A.E_j+2$, we get $A \eqv E_j$ by Lemma
\ref{nefred} and the contradiction $1=R.A=R.E_j=2$.

Therefore $R.E >0$ and $\beta \geq 3$, so that $R.E_1 < 0$. By Lemma \ref{A} we can write $E_1 \sim A + kR$ 
with $A > 0$ primitive, $A^2=0$ and $k:=-R.E_1=A.R \geq 1$. Now $2 \geq E.E_1=A.E+kE.R \geq A.E+k$ gives
$k=1$ by Lemma \ref{nefred}. Hence $2E+E_1+E_i$ is nef, whence base-point free, as $\phi(2E+E_1+E_i)
\geq 2$, and (a) is proved. Moreover, if $E.E_1=1$, then $E_1 \eqv E+R$ by Lemma \ref{nef:main2},
whence $E_1.E_i = E.E_i-1$, and (b) is proved, again since $\phi(E+E_1+E_i) \geq 2$.

To prove (c), assume $\gamma=2$ and $E.E_1=E_1.E_i=1$ and that there is a nodal curve $R$ with $R.(E+E_i) < 0$.
Then $R.E=0$ and $R.E_i=-1$ by quasi-nefness, and moreover $R.E_1 >0$ by (a). Therefore $E_1.E_i=1$ yields 
$E_i \eqv E_1+R$ with $R.E_1=1$, so that $E.E_i=E.E_1=1$. Moreover, if $M^2_2 \geq 2$ and $j \in
\{2,3\}-\{i\}$, then  $R.H=R.(\beta E + 2E_1 + E_i +E_j) = 1 + R.E_j$, so that $R.E_j \geq 0$. It follows that
$E_j.E_1=E_j.(E_i-R) \leq E_j.E_i \leq 2$, whence $E_1.M_2=E_1.(E_i+E_j) \leq 3$ and (c) is proved.

Moreover, by what we have just seen, under the assumptions in (d), if neither $E+E_2$ nor $E+E_3$ is nef, then
$E_2 \eqv E_1+R_2$ and $E_3 \eqv E_1+R_3$ with $R_2$ and $R_3$ nodal curves such that $R_2.E_1=R_3.E_1=1$. Then
$R_2.R_3=(E_2-E_1).(E_3-E_1)=-1$, a contradiction. This proves (d). 

Finally, to prove (e), assume that there is a nodal curve $R$ such that $R.(E_1+E_i) < 0$. By Lemma \ref{A} 
we can write $E_1+E_i \sim B + kR$ with $B > 0$, $B^2=2E_1.E_i>0$ and $k:=-R.(E_1+E_i) \geq 1$. By
(a) we have $E.R \geq 1$. From $2 = E.(E_1+E_i)=E.B+kE.R \geq 1 + k \geq 2$, we get $k = E.R = 1$, so that
$R.(E_1+E_i) = -1$. Now if $E_1.R \leq -1$ we get, by Lemma \ref{nef:main2}, that $E_1 \eqv E + R$, $E_1.R
= -1$ and $E_i.R = 0$, giving the contradiction $E_1.E_i = 1$. If $E_1.R \geq 0$ the quasi-nefness of $E_i$
implies that $E_1.R = 0$ and $E_i.R = -1$, whence $E_i \eqv E + R$ by Lemma \ref{A}, giving again the
contradiction $E_1.E_i = 1$.
\end{proof}

{\bf The general strategy} to prove the nonextendability of $S$ in the remaining cases (I), (II), (III), (D) 
and (S), will be as follows.

We will first use the ladder decomposition and Propositions \ref{ramextenr}-\ref{ramextenr4}
to reduce to particular genus one decompositions of $M_2$ or $M_3$ where we know all the intersections involved. 

Then we will find a big and nef line bundle $D_0$ on $S$ such that $\phi(D_0) \geq 2$ and $H-D_0$ is 
base-component free with $(H-D_0)^2 > 0$. In particular, this implies that $H^1(H-D_0)=H^1(D_0-H) =0$. In many
cases this $D_0$ will  satisfy the numerical conditions in Lemma \ref{subram},
so that $S$ will be nonextendable. 

In the remaining cases we will apply  Proposition \ref{mainextenr} in the following way.

We will let $D \in |D_0|$ be a general smooth curve. (This will be used repeatedly without further 
mentioning.) 

The surjectivity of the Gaussian map $\Phi_{H_D, \omega_D}$ will be handled by means of 
Theorem \ref{tendian}, to which we will refer. Moreover in all of the cases where  we will apply 
Theorem \ref{tendian} (with the exception of (e)) we will have that $h^0 (\O_D(2D_0 - H)) \leq 1$ 
if $D_0^2 \geq 6$ and $h^0 (\O_D(2D_0 - H)) = 0$ if $D_0^2 = 4$. Therefore the hypothesis (iii) of 
Proposition \ref{mainextenr} will always be satisfied and we will skip its verification.

To study the surjectivity of the multiplication map $\mu_{V_D, \omega_D}$ we will use several tools, 
outlined below. In several cases we will find an effective decomposition $D \sim D_1+D_2$ and use Lemma
\ref{multhelp}. We remark that {\bf except possibly for the one case in \eqref{mu1isot} below where $D_1$
is primitive of canonical type, both $D_1$ and $D_2$ will always be smooth curves}. The reason for this
is that we will always have that $|D_1|$ and $|D_2|$ are base-component free and not multiple of elliptic 
pencils, whence their general members will be smooth and irreducible \cite[Prop.3.1.4 and Thm.4.10.2]{cd}. In
most cases this will again be used without further mentioning.

Furthermore the spaces $V_D$, $V_{D_1}$ and $V_{D_2}$ will always be base-point free. This is 
immediately clear for $V_D$, as $|D_0|$ is base-point free. As for $V_{D_1}$ and $V_{D_2}$, one only has to
make sure that, in the cases where $|H-D_0|$ has base points (that is, $\phi(H-D_0)=1$), in which case it
has precisely two distinct base points \cite[Prop.3.1.4 and Thm.4.4.1]{cd}, they do not intersect the
possible base points of $|D_1|$ and $|D_2|$. This will always be satisfied, and again, we will not
repeatedly mention this.

Here are the criteria we will use to verify that the desired multiplication maps are surjective:
\begin{eqnarray} 
\nonumber & \mu_{V_D, \omega_D} \; \mbox{is surjective in any of the following cases:} &   \\
\label{biraz} & H^1(H-2D_0) = 0 \; \mbox{and either $|D_0|$ or $|H-D_0|$ is birational (see Rmk.
\ref{rem:bir} below)}.  & \\
\label{pencil}  & H^1(H-2D_0) = 0 \; \mbox{and $|H-D_0|$ is a pencil}. &
\end{eqnarray}
\begin{eqnarray} 
\nonumber & \mbox{If $V_{D_1}$ is base-point free, then } \mu_{V_{D_1}, \omega_{D_1}} \; \mbox{is surjective in 
any of the following cases:} 
& \\
\label{mu1}  & H^1(H-D_0-D_1) = 0, \; \mbox{$D_1$ is smooth and} \; (H-D_0).D_1 \geq D_1^2 + 3;  & \\
\label{mu1isot} &   H^1(H-D_0-D_1) = 0 \; \mbox{and $D_1$ is nef and isotropic.}  & 
\end{eqnarray}
\begin{eqnarray} 
\nonumber & \mbox{If $D_2$ is smooth and $V_{D_2}$ is base-point free, then } \mu_{V_{D_2}, \omega_{D_2}(D_1)} 
\; \mbox{is surjective if}  & \\ 
\label{mu2}  & h^0(H-D_0-D_2) + h^0(\O_{D_2}(H-D_0-D_1)) \leq
\frac{1}{2} (H-D_0)^2-1 \; \mbox{(see Rmk. \ref{rem:rrcliff})}.& 
\end{eqnarray}
To see \eqref{biraz}-\eqref{pencil} note that if $H^1(H-2D_0) = 0$, we have that $V_D = H^0((H - D_0)_{|D})$, 
whence \eqref{pencil} is just the base-point free pencil trick, while \eqref{biraz} 
follows from the base-point free pencil trick and \cite[Thm.1.6]{as} since $(H-D_0)_{|D}$ is base-point free 
and is either a pencil or birational by any of the hypotheses. The same proves \eqref{mu1}. As for
\eqref{mu1isot} the hypotheses imply $V_{D_1} = H^0((H - D_0)_{|D_1})$ and $\omega_{D_1} \cong \O_{D_1}$, as
$D_1$ is either a smooth elliptic curve or indecomposable of canonical type \cite[III, \S 1]{cd}, and the
surjectivity is immediate.

Finally for \eqref{mu2} we use the $H^0$-lemma \cite[Thm.4.e.1]{gr}, which states that
$\mu_{V_{D_2}, \omega_{D_2}(D_1)}$ is surjective if
\begin{eqnarray*} \dim V_{D_2} - 2 & = & h^0(H-D_0) - h^0(H-D_0-D_2)-2 \geq
\\ & \geq &
h^1(\omega_{D_2}(D_1-(H-D_0))=h^0(\O_{D_2}(H-D_0-D_1)). 
\end{eqnarray*}
Using Riemann-Roch on $S$, this is equivalent to \eqref{mu2}.

\begin{rem} 
\label{rem:bir}
{\rm A complete linear system $|B|$ is birational if it defines a birational map.
By \cite[Prop.3.1.4, Lemma4.6.2, Thm.4.6.3, Prop.4.7.1 and Thm.4.7.1]{cd} a nef divisor $B$ with 
$B^2 \geq 8$ defines a birational morphism if $\phi(B) \geq 2$ and $B$ is not $2$-divisible in $\Pic S$ when 
$B^2 = 8$.} 
\end{rem}

\begin{rem} 
\label{rem:rrcliff}
{\rm The inequality in \eqref{mu2} will be verified by giving an upper bound on
$h^0(H-D_0-D_2)$ and using Riemann-Roch and Clifford's theorem on $D_2$ to bound $h^0(\O_{D_2}(H-D_0-D_1))$. 
We will often not mention this.  } 
\end{rem}

\section{Case (D)}
\label{caseD}

We have $H \eqv 2H_1$ whence $H_1$ is ample with $H_1^2 \geq 8$ and $\phi(H) = 2 \phi(H_1) \geq
3$ gives $\phi(H_1) \geq 2$.

\subsection{$H \sim 2H_1+K_S$} 

We will prove that $S \subset \PP^r$ is nonextendable. We set $D_0 = H_1$ and check that the
hypotheses of Proposition \ref{mainextenr} are satisfied by $D_0$. Of course $D_0$ is ample and $\phi(D_0) 
\geq 2$. The Gaussian map $\Phi_{H_D, \omega_D}$ is surjective by Theorem \ref{tendian}(c) since $H^0(2D_0-H) =
H^0(K_S) = 0$. Moreover also $H^1(H-2D_0) = 0$, whence and $(H-D_0)_{|D} = \omega_D$ so that the multiplication
map $\mu_{V_D, \omega_D}$ is just $\mu_{\omega_D, \omega_D}$ which is surjective since $D$ is not hyperelliptic.

\subsection{$H \sim 2H_1$}

\subsubsection{$\phi(H_1) \geq 3$}

In the course of the proof of Corollary \ref{cor:anysurface} we have seen that, given an
extendable Enriques surface $S \subset \PP^r$, we can reembed it in such a way that it is
linearly normal and still extendable. Therefore we can assume that $S \subset \PP H^0(2H_1)$. By
\cite[Cor.2, page 283]{cd} it follows that $H_1$ is very ample, whence $S \subset \PP
H^0(2H_1)$ is $2$-Veronese of $S_1 = \varphi_{H_1}(S) \subset \PP H^0(H_1)$ and therefore $S
\subset \PP^r$ is nonextendable by \cite[Thm.1.2]{glm}.

\subsubsection{$\phi(H_1) = 2$ and $H_1^2 = 8$}

Since $E.H_1 = 2$ we have $(H_1 - 2E)^2 = 0$ and we can write $H_1 = 2E +
F$ with $F > 0, F^2 = 0, E.F = 2$. According to whether $F$ is primitive or not, we get cases (a1) \label{a1}
and (a2) \label{a2} in the proof of Proposition \ref{precisa}.

\subsubsection{$\phi(H_1) = 2$ and $H_1^2 = 10$}

Since $E.H_1 = 2$ we have $(H_1 - 2E)^2 = 2$ and we can write $H_1 = 2E +
E_1 + E_2$ with $E_i > 0, E_i^2 = 0, E.E_i = E_1.E_2 = 1$. Then 
$H \sim 4E + 2E_1 + 2E_2$ with $\phi(H)=4$. Note that $\alpha = 2$ and $E_1.L_1 = E_2.L_1 = \phi(L_1)$.

First we show that either $E_1$ or $E_2$ is nef. If not, by Lemma \ref{nef:main2}, we have $E_1 \eqv
E+\Gamma_1$, $E_2 \eqv E+\Gamma_2$ for two nodal curves $\Gamma_i$ with $\Gamma_i.E=1$. But then
$\Gamma_1.\Gamma_2 = (E_1-E).(E_2-E)=-1$, a contradiction. Therefore, replacing $E_1$ with $E_2$ if necessary, 
we can assume that $E_1$ is nef. Moreover, possibly after substituting $E_2$ with $E_2+K_S$, we can 
assume that either $E_2$ is nef or there exists a nodal curve $\Gamma$ such that 
$E_2\sim E+\Gamma+K_S$. In particular $E+E_2$ is base-component free.

We set $D_0= E+2E_1+E_2$, which is nef with $\phi(D_0)=2$, $D_0^2=10$ and $H-D_0 \sim 3E+E_2$ is base-component
free. Moreover $2D_0-H \sim -2E+2E_1$ and $h^0(2D_0-H)=0$ by the nefness of $E_1$, whence $\Phi_{H_D, \omega_D}$ is
surjective by Theorem \ref{tendian}(c).

Now to see the surjectivity of $\mu_{V_D, \omega_D}$ note that since
$h^1(-2E_1)=h^1(K_S)=0$ the two restriction maps $H^0(E+E_2) \khpil H^0(\O_D(E+E_2))$ and
$H^0(E+2E_1+E_2+K_S) \khpil H^0(\omega_D)$ are surjective and $|\O_D(E+E_2)|$ is a base-point 
free pencil.

Consider the commutative diagram
\[  
\xymatrix{
H^0(2E) \* H^0(E+E_2) \* H^0(E+2E_1+E_2+K_S)  \ar[d]^{\mu_{2E, E+E_2}} 
\ar[r]^{\hskip 1.2cm r_D} &  W_D \* H^0(\O_D(E+E_2)) \* H^0(\omega_D) 
\ar[d]^{\Id \* \mu_{\O_D(E+E_2), \omega_D}}
\\
H^0(H-D_0) \* H^0(D_0+K_S) \ar[d]^{p_D} &  W_D \* H^0(\omega_D(E+E_2)) 
\ar[d]^{\mu_{W_D, \omega_D(E+E_2)}}
\\
V_{D} \* H^0(\omega_D) 
\ar[r]^{\mu_{V_D, \omega_D}} &  H^0(\O_D(H+K_S)), 
}
\]
where $p_D$ and $r_D$ are the natural restriction maps and $W_D := \Im \{ H^0(2E) \khpil
H^0(\O_D(2E)) \}$. The map $r_D$ is surjective and the map $\mu_{\O_D(E+E_2), \omega_D}$ is
surjective by the base-point free pencil trick. Now to prove that $\mu_{V_D, \omega_D}$
is surjective it suffices to show that $\mu_{W_D, \omega_D(E+E_2)}$ is surjective. 

Since $\dim W_D=2$ and $W_D$ is base-point free, the surjectivity of 
$\mu_{W_D, \omega_D(E+E_2)}$ follows by the $H^0$-lemma \cite[Thm.4.e.1]{gr} as soon as we 
prove that $h^1(\omega_D(E_2-E)) = 0$. Now $h^1(\omega_D(E_2-E))=h^0(\O_D(E-E_2))$ and, since
$\deg \O_D(E-E_2)=0$, the required vanishing follows unless $\O_D(E-E_2)$ is trivial. But if the
latter holds, then we have a short exact sequence
\[ 0 \hpil -2E_1-2E_2 \hpil E-E_2 \hpil \O_D \hpil 0. \]

If $E_2$ is nef then $h^0(E-E_2)=h^1(-2E_1-2E_2)=0$ and we get the contradiction $h^0(\O_D)=0$.

If $E_2$ is not nef then $E_2 \sim E+\Gamma+K_S$, and $\Gamma.D=0$ whence
$h^1(\omega_D(E_2-E)) = h^1(\O_D(D_0+\Gamma)) = h^1(\O_D(D_0)) = 0$.  

Therefore in both cases $\mu_{V_D, \omega_D}$ is surjective and $S$ is
nonextendable by Proposition \ref{mainextenr}.

\subsubsection{$\phi(H_1) = 2$ and $H_1^2 \geq 12$}

We set $D_0 = H_1$ and check that the hypotheses of Proposition \ref{mainextenr} are satisfied. Of course $D_0$
is ample and $\phi(D_0) = 2$. The Gaussian map $\Phi_{H_D, \omega_D}$ is surjective by Theorem \ref{tendian}(d)
since $h^0(2D_0-H) = h^0(\O_S) = 1$. Moreover also $H^1(H-2D_0) = 0$, whence $\mu_{V_D, \omega_D}$ is
surjective by (\ref{biraz}). Hence $S \subset \PP^r$ is nonextendable.

\section{Case (I) with $M_2 =0$}
\label{caseI-0}

We have $H \eqv \beta E + \gamma E_1$ with $E.E_1 = 1, \beta \geq 2,
\gamma \in \{2,3\}$ and $H^2 \geq 32$ or $H^2 = 28$. 

Now $\gamma = E.H = \phi(H) \geq 3$, whence $\gamma = 3$ and we will deal with
$H \eqv \beta E + 3 E_1, E.E_1 = 1, \beta \geq 6.$ 

We set $D_0 = \lfloor \frac{\beta}{2} \rfloor E + 2 E_1$ and check the
conditions of Proposition \ref{mainextenr}. We have that $D_0$ is nef by Lemma \ref{nef:main}, $\phi(D_0) = 2$ 
and $D_0^2 = 4 \lfloor \frac{\beta}{2} \rfloor \geq 12$. Now $H - D_0 \eqv \lfloor \frac{\beta+1}{2} \rfloor E +
E_1$ and $2D_0 - H \eqv E_1$(respectively $E_1 - E$) if $\beta$ is even (respectively if $\beta$ is odd). 
Therefore, by Lemma \ref{nef:main2}, we have that $h^0(2D_0 - H) \leq 1$ and the Gaussian map
$\Phi_{H_D, \omega_D}$ is surjective by Theorem \ref{tendian}(c)-(d). To prove the
surjectivity of the multiplication map $\mu_{V_D, \omega_D}$ we apply Lemma \ref{2E}. We have
$D_0+K_S-2E = (\lfloor \frac{\beta}{2} \rfloor - 2) E + 2 E_1 + K_S$ whence $H^1(D_0+K_S-2E)
= 0$ by Lemma \ref{nef:main2}.  Also $H^2(D_0+K_S-4E) = 0$ since $E.(D_0+K_S-4E)
= 2$. Now $H-D_0-2E \eqv \lfloor \frac{\beta-3}{2} \rfloor E + E_1$ is base-component free by Lemma
\ref{nef:main2} and $H-2D_0-2E \eqv (\beta - 2 \lfloor \frac{\beta}{2} \rfloor - 2) E - E_1$, whence
$h^0(H-2D_0-2E)=0$ by the nefness of $E$. Also $H-2D_0-4E \eqv - \varepsilon E - E_1, \varepsilon = 3, 4$, 
whence $H^1(H-2D_0-4E) = 0$ and the exact sequence
\[ 0 \hpil H-2D_0-4E \hpil H-D_0-4E  \hpil \O_D(H-D_0-4E)  \hpil  0 \]
shows that $h^0(\O_D(H-D_0-4E)) \leq h^0(H-D_0-4E)$. Since $H-D_0-4E \eqv \lfloor \frac{\beta-7}{2} \rfloor 
E + E_1$ we have that $h^0(H-D_0-4E) = \lfloor \frac{\beta-5}{2} \rfloor$ if $\beta \geq 7$, while if $\beta
= 6$, we have $H-D_0-4E \eqv - E + E_1$ and replacing $D_0$ with $D_0 + K_S$ if necessary,
we can assume that $h^0(H-D_0-4E) = 0$ by Lemma \ref{nef:main2}. In both cases we get that (\ref{green}) is 
satisfied. 

\section{Case (I) with $\gamma =3$ and $M_2 >0$}
\label{caseI-3}

We first note that, since $\gamma=3$, we must have $\beta \geq 3$. Indeed,
if $\beta=2$ we have $L_2 \sim E_1+M_2$ and $E.L_2=1+E.M_2= \phi(H)-2
\leq E_1.H - 2 = E_1.M_2 = E_1.L_2$, contradicting the removing conventions of Section \ref{class}, page
\pageref{remconv} (because then $(L_2 - E)^2 \geq (L_2 - E_1)^2 \geq 0$, therefore we could find $E$ in a genus 
$1$ decomposition of $L_2$, but then $\beta \geq 3$). Then we have, in this section,
\begin{equation} 
\label{eq:I-3.1}
H \sim \beta E + 3 E_1 + M_2, \hs E.E_1=1, \hs \beta \geq 3, \; 
M_2 >0, \hs M_2^2 \geq 0 \; \mbox{ and } \; H^2 \geq 32 \; \mbox{ or } \; H^2 = 28.
\end{equation}  

We will use the following:

\begin{lemma} 
\label{nef:I-3}
Let $N = E+E_1$. Then either $|H-2N|$ is base-point free and $h^1(H-3N)=0$, or $S$ is
nonextendable.
\end{lemma}

\begin{proof}
We have $H-3N \sim (\beta-3) E + M_2>0$ and $(H-3N)^2 \geq 0$ with equality only if 
$H-3N = M_2$. Since $M_2$ is of small type it follows by Theorem \ref{lemma:qnef} that
$h^1(H-3N)=0$ if and only if $H-3N$ is quasi-nef. Moreover, as $(H-2N)^2 = 2(\beta-2)(E.M_2+1)
+ 2E_1.M_2 +M_2^2 \geq 6$ and $\phi(H-2N) \geq 2$, we have that $H-2N$ is base-point free if
and only if it is nef. Now if $H-2N$ is not nef, there is a nodal curve $\Gamma$  with
$\Gamma.(H-2N) <0$. Then $\Gamma.N >0$, since $H$ is ample, whence $\Gamma.(H-3N) \leq -2$. 

Therefore, to prove the lemma, it suffices to show that $S$ is nonextendable if $H-3N$ is not
quasi-nef.

Let $\Delta >0$ be such that $\Delta^2=-2$ and $\Delta.(H-3N) \leq -2$. We have $\Delta.N >0$
since $H$ is ample. Also note that $\Delta.E_1 \geq 0$, for if not, we would have $\Delta.E
\geq 2$, whence the contradiction $(E+\Delta)^2 \geq 2$ and $E_1.(E+\Delta) \leq 0$.
Hence $M_2.\Delta \leq -2$ and by Lemma \ref{A} we can write $M_2 \sim A + k\Delta$, with
$A>0$ primitive, $A^2 =M_2^2$ and $k:=-\Delta.M_2=\Delta.A \geq 2$. Now
if $E.\Delta > 0$ we find $E.M_2 \geq k$ and if equality holds, then $E.A=0$ and $E.\Delta=1$, whence $E \eqv A$
by Lemma \ref{nefred}, a contradiction. We get the same contradiction if $E_1.\Delta > 0$. Therefore
\begin{equation} 
\label{eq:star}
E.M_2 \geq -\Delta.M_2+1 \geq 3 \; \mbox{if} \; E.\Delta >0 \; \mbox{and} \; 
E_1.M_2 \geq -\Delta.M_2+1 \geq 3 \; \mbox{if} \; E_1.\Delta > 0. 
\end{equation}

We first consider the case $E.\Delta>0$.

Note that we cannot have that $\beta=3$, for otherwise $H$ is of type (I-B) in Section
\ref{class} and $L_2 \sim (3-\alpha)E+E_1+M_2$ is of small type, whence $E_1.M_2 \leq 5$ by Lemma 
\ref{STlemma2}, so that $E_1.(H-2E)=E_1.(E+3E_1+M_2) \leq 6$. Since $\phi(H)=E.H=3+E.M_2 \geq 6$ by
\eqref{eq:star}, we get $\alpha=2$ and $E_1.H = 3 + E_1.M_2 \geq 6$, so that $E_1.M_2 \geq 3$. Hence $L_2 \sim 
E+E_1+M_2$ and $L_2^2 \geq 14$, a contradiction.

Therefore $\beta \geq 4$, whence $\Delta.M_2 \leq -2 - (\beta-3) \Delta.E \leq -3$,
so that $E.M_2 \geq 4$ by \eqref{eq:star} and $\phi(H)\geq 7$, whence $H^2 \geq 54$ by
\cite[Prop.1]{kl1}. Now one easily verifies that $B:=2E+E_1+\Delta$ satisfies the 
conditions in Proposition \ref{ramextenr}, so that $S$ is nonextendable.

We finally consider the case $\Delta.E=0$, where $E_1.\Delta >0$, so that
$E_1.M_2 \geq 3$ by \eqref{eq:star}.

Now $L_2 \sim (\beta-\alpha)E+E_1+M_2$ if $H$ is of type (I-B) in Section \ref{class} and $L_3 \sim
(\beta-\alpha-2)E+E_1+M_2$ if $H$ is of type (I-A) in Section \ref{class}. We claim that the removing 
conventions of Section \ref{class}, page \pageref{remconv} now imply that $E_1.M_2 \leq E.M_2 + 1$ and, if
$\beta=3$, that $E_1.M_2 \leq E.M_2$. In fact if the latter inequalities do not hold we have that $E.L_2 \leq
E_1.L_2$, $E.L_3 < E_1.L_3$ and $(E_1+M_2-E)^2 \geq 0$, contradicting the fact that 
$L_2$ and $L_3$ are of small type. To summarize, we must have
\begin{equation} 
\label{eq:star2}
E.M_2 \geq 2, \; \mbox{and furthermore} \; E.M_2 \geq 3\; \mbox{if} \; \beta=3. 
\end{equation}
This yields $H^2 =6\beta+2 \beta E.M_2+6E_1.M_2 + M_2^2 \geq 54$ in any case. Moreover, using
\eqref{eq:star2}, one easily verifies that $B:=E+2E_1+\Delta$ satisfies the conditions in
Proposition \ref{ramextenr}, so that $S$ is nonextendable.
\end{proof}

The main result of this section is the following:

\begin{prop} 
\label{sum:I-3}
If $H$ is of type {\rm (I)} with $\gamma =3$ and $M_2 >0$ then $S$ is nonextendable.
\end{prop}

\begin{proof}
Set $D_0=2N = 2(E+E_1)$, which is nef by Lemma \ref{nef:main} with $D_0^2=8$ and
$\phi(D_0)=2$. By Lemma \ref{nef:I-3} we can assume that $H-D_0$ is base-point free as well. 

By assumption we have $H.D_0=2(\beta + 3 + (E+E_1).M_2) \geq 16$ with equality only if $\beta = 3$ and 
$E.M_2 = E_1.M_2 = 1$. But in the latter case, since $M_2$ does not contain $E$ in its arithmetic genus $1$ 
decompositions, we have, by Lemma \ref{lemma0}, that $M_2^2=0$ and $H^2=30$, a contradiction. Hence 
$(2D_0-H).D_0 < 0$ and consequently $h^0(2D_0-H)=0$, so that $\Phi_{H_D, \omega_D}$ is surjective
by Theorem \ref{tendian}(c).

Let $D_1, D_2 \in |N|$ be two general members. By Lemma \ref{nef:I-3} we can assume
$h^1(H-D_0-D_i)=0$. Hence $\mu_{V_{D_1}, \omega_{D_1}}$ is surjective by (\ref{mu1}) since
$(H-D_0).D_1=\beta -1 + (E+E_1).M_2 \geq 5$ by our assumptions, and the map 
$\mu_{V_{D_2}, \omega_{D_2}(D_1)}= \mu_{\O_{D_2}(H-D_0), \omega_{D_2}(D_1)}$ is surjective 
by \cite[Cor.4.e.4]{gr} since $\deg \omega_{D_2}(D_1) = 2g(D_2)$ and $\deg \O_{D_2}(H-D_0) \geq
2g(D_2)+1$.

By Lemma \ref{multhelp}, $\mu_{V_D, \omega_D}$ is surjective and by Proposition
\ref{mainextenr}, $S$ is nonextendable.
\end{proof}

\section{Case (I) with $\gamma =2$ and $M_2 > 0$}
\label{caseI-2}

We have
\[ H \sim \beta E + 2 E_1 + M_2, \hs E.E_1=1, \hs M_2 >0, \hs M_2^2 \geq 0, \hs H^2 \geq 32 \mbox{ or } 
H^2 = 28 \]
and, as above, either $L_2$ or $L_3$ is of small type, whence so is $M_2$. 

Recall that $E_1.M_2 \leq E_1.M_2+\beta-\alpha =E_1.L_1 =\phi(L_1) \leq \phi(H)=2+E.M_2 \leq E_1.H
= \beta + E_1.M_2$. Moreover, since by construction $M_2$ neither contains 
$E$ nor $E_1$ in its arithmetic genus $1$ decompositions, we have, by Lemma \ref{lemma0}, that 
$(M_2 - E)^2 < 0$ and $(M_2 - E_1)^2 < 0$. Hence 
\begin{eqnarray}
\label{eq:I-R.2} 
\frac{1}{2}M_2^2 +1 \leq & E.M_2 & \leq E_1.M_2 + \beta-2, \hs \mbox{and} \\
\label{eq:I-R.1}  
\frac{1}{2}M_2^2 +1 \leq & E_1.M_2 & \leq E.M_2 + 2 - \beta + \alpha \leq
E.M_2 + 2.
\end{eqnarray}

In this section we will first prove, in Lemma \ref{lemma:inizioI}, nonextendability up to some explicit
decompositions of $H$ and we will then prove, in Proposition \ref{cor:I-2.1}, nonextendability if 
$\beta \geq 5$. 

\begin{lemma} 
\label{lemma:inizioI}
Let $H$ be of type {\rm (I)} with $\gamma=2$ and $M_2^2 \geq 2$. Then $S$ is nonextendable unless, 
possibly, we are in
one of the following cases (where all the $E_i$'s are effective and isotropic):
\begin{itemize}
\item[(a)] $M_2^2 =2$, $M_2 \sim E_2+E_3$, $E_2.E_3=1$, and either
\begin{itemize}
\item[(a-i)]  $\beta=2$, $(E.E_2, E.E_3, E_1.E_2, E_1.E_3)=(1,2,1,2),(1,2,2,1),(1,1,2,2),(2,2,2,2)$,
\linebreak $(1,2,2,2)$; or 
\item[(a-ii)] $\beta=3$, $(E.E_2, E.E_3, E_1.E_2, E_1.E_3)=(2,2,2,2),(2,2,1,2)$; or 
\item[(a-iii)] $\beta \geq 3$, $(E.E_2, E.E_3, E_1.E_2, E_1.E_3)=(1,1,1,1),(1,1,1,2),(1,1,2,2)$.
\end{itemize}
\item[(b)] $M_2^2=4$, $M_2 \sim E_2+E_3$, $E_2.E_3=2$, and either
\begin{itemize}
\item[(b-i)]  $\beta=2$, $(E.E_2, E.E_3, E_1.E_2, E_1.E_3)=(1,2,1,2),(1,2,2,1),(1,2,2,2),(1,2,1,3)$; or 
\item[(b-ii)] $\beta=3$, $(E.E_2, E.E_3, E_1.E_2, E_1.E_3)=(1,2,2,1),(1,2,1,3)$.
\end{itemize}
\item[(c)] $M_2^2=6$, $M_2 \sim E_2+E_3+E_4$, $E_2.E_3=E_2.E_4=E_3.E_4=1$, and 
\begin{itemize}
\item[] $\beta=2$, $(E.E_2, E.E_3, E.E_4, E_1.E_2, E_1.E_3, E_1.E_4)=(1,1,2,1,1,2)$.
\end{itemize}
\end{itemize}
\end{lemma}

\begin{proof} 
We write $M_2 \sim E_2 + \ldots + E_{k+1}$ as in Lemma \ref{STlemma2} with $k =
2$ or $3$. Moreover we can assume that $1 \leq E.E_2 \leq \ldots \leq E.E_{k+1}$, whence that 
$E.M_2 \geq kE.E_2$.

We first consider the case $\beta \geq 4$.

We note that $(M_2-2E_2)^2=-2$ if $M_2^2=2$ or $6$, $(M_2-2E_2)^2=-4$ if $M_2^2=4$ and $(M_2-2E_2)^2 \geq -6$ 
if $M_2^2=10$. In the latter case $E.M_2 \geq 6$ by \eqref{eq:I-R.2}, whence $E.(M_2-2E_2) \geq 2$. Using this
and setting $B:=E+E_1+E_2$ one easily verifies that $(H-2B)^2 = 2(\beta-2)E.(M_2-2E_2) +(M_2-2E_2)^2 \geq 0$ and
$E.(H-2B) >0$ (whence $H-2B \geq 0$ by Riemann-Roch), except for the cases
\begin{equation} 
\label{eq:10.1.1}
M_2^2=2,4 \; \mbox{and} \; E.E_2=E.E_3. 
\end{equation}
Moreover, except for these cases, using (\ref{eq:I-R.2}) and (\ref{eq:I-R.1}), one easily verifies that 
$H^2 \geq 54$, except for the case $\beta=4$, $M_2^2=2$ and $(E.M_2,E_1.M_2)=(3,2)$, where $H^2=50$. In this
case $(3B-H).H =4 < \phi(H)=5$, so that, if $3B-H>0$ it must be a nodal cycle. Therefore either $h^0(3B-H)=0$ or 
$h^0(3B+K_S-H)=0$, so in any case $B$ satisfies the conditions in Propositions \ref{ramextenr} or
\ref{ramextenr6} and $S$ is nonextendable. 

In the remaining cases \eqref{eq:10.1.1} we can without loss of generality assume $1 \leq E_1.E_2 \leq E_1.E_3$
and we set $B:=E+E_2$. Then $(H-2B)^2 = 4(\beta-2)+4E_1.(E_3-E_2) + (E_3-E_2)^2 \geq 4$ and $(H-2B).E = 2$.
Using (\ref{eq:I-R.2}) and (\ref{eq:I-R.1}), one gets $H^2 \geq 64$ if $M_2^2=4$, $H^2 \geq 74$ if $M_2^2=2$ and
$E.E_2=E.E_3 = 3$, and $B.H \geq 17$ if $M_2^2=2$ and $E.E_2=E.E_3 =2$. Moreover, in the latter case, we have
that again $H^2 \geq 64$ unless $\beta=4$ and $E_1.M_2 = 2, 3$, which gives $E_1.E_2=1$ and $B$ is nef by Lemma
\ref{lemmone}(c) since $E_2.H = 11 < 2\phi(H)=12$, whence $E_2$ is quasi-nef by Lemma \ref{qnef0}. Therefore
$B$ satisfies the conditions in Propositions \ref{ramextenr} or \ref{ramextenr4} and $S$ is nonextendable unless
$M_2^2 = 2$ and $E.E_2=E.E_3=1$. In the latter case, by (\ref{eq:I-R.1}) we have 
$2 \leq E_1.M_2 \leq 4$ with $E_1.M_2=4$ if and only if $\alpha=\beta$. In this last case $L_1 \sim 2E_1+M_2$,
whence $\phi(L_1)=E_1.M_2=4$ and we get that $4 \leq E_i.L_1 = 2E_1.E_i+1$ for $i=2,3$, so that
$E_1.E_2=E_1.E_3=2$. Therefore we get the cases in (a-iii) with $\beta \geq 4$.

We next treat the cases $\beta \leq 3$. Then we must be in case (I-B) of Section \ref{class}, whence $L_2$ is of
small type and either $L_2 \sim M_2$ or $L_2 \sim E+M_2$.

Suppose first that $L_2 \sim E+M_2$.

Then $\beta \geq \alpha + 1 \geq 3$, whence $\beta = 3$, $\alpha = 2$ and, since $L_2$ is of small type, by 
(\ref{eq:I-R.2}), we can only have $(M_2^2,E.M_2)=(2,2)$, $(2,4)$ or $(4,3)$.

If $(M_2^2,E.M_2)=(2,2)$, then $E.E_2=E.E_3=1$ and by (\ref{eq:I-R.1}) we have $2 \leq E_1.M_2 \leq 3$, 
yielding the first two cases in (a-iii).

If $(M_2^2,E.M_2)=(2,4)$, then $L_2^2=10$ and $\phi(L_2)=3$. As $E.E_i+1  = L_2.E_i \geq \phi(L_2)=3$ for 
$i = 2, 3$, we must have $E.E_2=E.E_3=2$. Now $L_1 \sim E+2E_1+M_2$ and $(1+E_1.M_2)^2=\phi(L_1)^2 \leq
L_1^2=14+4E_1.M_2$ and (\ref{eq:I-R.2}) yield $E_1.M_2=3$ or $4$. Therefore, by Lemma \ref{fi} and symmetry, we
get the two cases in (a-ii).

If $(M_2^2,E.M_2)=(4,3)$, then $E_1.M_2=3$ or $4$ by (\ref{eq:I-R.1}). Since $L_2^2=10$ and $\phi(L_2)=E.L_2=3$, 
there is by \cite[Cor.2.5.5]{cd} an isotropic effective $10$-sequence $\{f_1, \ldots, f_{10}\}$ such that $E=f_1$
and $3L_2 \sim f_1+ \ldots + f_{10}$. 

In the case $E_1.M_2=3$ we get $E_1.L_2=4$, whence we can assume, possibly after renumbering, that
$E_1.f_i=1$ for $1 \leq i \leq 8$ and $(E_1.f_9,E_1.f_{10})=(2,2)$ or $(1,3)$. In the 
latter case we have $(E_1+f_{10})^2=6$ and $\phi(E_1+f_{10})=2$, whence we can write $E_1+f_{10} \sim
A_1 + A_2 + A_3$ for some $A_i>0$ such that $A_i^2 = 0$, $A_i.A_j=1$ for $i \neq j$. But
$f_i.(E_1+f_{10})=2$ for all $1 \leq i \leq 9$, a contradiction. Hence  $(E_1.f_9,E_1.f_{10})=(2,2)$. 
One easily sees that there is an isotropic divisor $f_{19}>0$ such that
$f_{19}.f_1=f_{19}.f_9=2$ and $L_2 \sim f_1+f_9+f_{19}$. Therefore $E_1.f_{19}=1$. Setting $E'_2=f_9$
and $E'_3=f_{19}$  we get the first case in (b-ii).

In the case $E_1.M_2=4$ we get $E_1.L_2=5$, whence we can assume, possibly after renumbering, that
$E_1.f_i=1$ for $1 \leq i \leq 5$. As above there is an isotropic divisor $f_{12}>0$ 
such that $f_{12}.f_1=f_{12}.f_2=2$ and $L_2 \sim f_1+f_2+f_{12}$. Therefore $E_1.f_{12}=3$. 
Setting $E'_2=f_2$ and $E'_3=f_{12}$ we obtain the second  case in (b-ii).

Finally, we have left the case with $L_2 \sim M_2$, where $\beta=\alpha$. 

We have $L_1 \sim 2E_1+M_2$, whence $(E_1.M_2)^2 =\phi(L_1)^2 \leq L_1^2 =4E_1.M_2+M_2^2$, so that
(\ref{eq:I-R.1}) and \cite[Prop.1]{kl1} give $E_1.M_2 \leq 4$. In particular $M_2^2 \leq 6$ by (\ref{eq:I-R.1}).

If $\beta=\alpha=3$, by definition of $\alpha$, we must have $1+E_1.M_2=E_1.(L_1+E) > \phi(H)=2+E.M_2$, whence
$E_1.M_2=4$, $E.M_2=2$ and $M_2^2=2$ by (\ref{eq:I-R.2}). Then $E.E_2=E.E_3=1$ and for $i = 2, 3$ we
have $E_i.L_1 = 2E_i.E_1 + 1 \geq \phi(L_1) = E_1.M_2 = 4$, whence $E_1.E_2=E_1.E_3=2$ and we get the third 
case in (a-iii).

In the remaining cases we have $\beta=\alpha=2$.

If $M_2^2=2$ using again $\phi(L_1)^2 \leq L_1^2$, $E_i.L_1 \geq \phi(L_1)$, \eqref{eq:I-R.2} and
\eqref{eq:I-R.1} together with $H^2 \geq 32$ or $H^2=28$, we deduce the possibilities $(E.M_2,E_1.M_2)=(3,3)$,
$(2,4)$, $(3,4)$ or $(4,4)$. By symmetry one easily sees that one gets the cases in (a-i).

If $M_2^2=4$ we similarly get $(E.M_2,E_1.M_2)=(3,3)$, $(3,4)$ or $(4,4)$. From the first two cases, using 
Lemma \ref{fi} for the second, we obtain the cases in (b-i). If $(E.M_2,E_1.M_2)=(4,4)$, we now show that $H$
also has a ladder decomposition of type (III). It will follow that $S$ is nonextendable by Section
\ref{caseIII}.

We have $E.H = 6$, whence $(H - 3E)^2 = 8$ and $H - 3E>0$ by Lemma \ref{nefred}. If $\phi(H-3E)=1$ we can
write $H - 3E \sim 4A_1 + A_2$ with $A_i>0$, $A_i^2=0$ and $A_1.A_2=1$. Now $6 \leq H.A_1 = 3E.A_1+1$ gives 
$E.A_1 \geq 2$, whence the contradiction $6 = H.E = 4E.A_1 + E.A_2 \geq 8$. Therefore there is an $E_1'>0$ such
that $(E_1')^2 = 0$ and $E_1'.(H-3E)=2$. Since $(H-3E-2E_1')^2=0$, by Lemma \ref{nefred} we can write $H \sim
3E+2E_1'+E_2'$, with $E_2'>0$, $(E_2')^2 = 0$ and $E_1'.E_2'=2$. From $6 \leq H.E_1' = 3E.E_1'+2$ we get
$E.E_1' \geq 2$. Now from $6 = H.E = 2E.E_1' + E.E_2'$ we see that we cannot have $E.E_1' \geq 3$, for then 
$E.E_1' = 3$, $E.E_2' = 0$, but this gives $E_2' \eqv qE$ for some $q \geq 1$ by Lemma \ref{nefred}, whence the 
contradiction $2 = E_1'.E_2' = 3q$. Therefore $E.E_1' = 2$ and since $E_1'.L_1=E_1'.(H-3E)+E_1'.E=4=\phi(L_1)$ we
obtain a ladder decomposition of $H$ of type (III), as claimed.

If $M_2^2=6$, by (\ref{eq:I-R.2}) and (\ref{eq:I-R.1}) we get, as above, $E_1.M_2=E.M_2=4$, yielding by
symmetry the case in (c) in addition to the case  $(E.E_2, E.E_3, E.E_4, E_1.E_2, E_1.E_3,
E_1.E_4)=(1,1,2,1,2,1)$. In the latter case we note that $\phi(H)=E.H=E_1.H=6$ and
$\phi(H-2E_1)=\phi(2E+E_2+E_3+E_4)=E_3.(H-2E_1)=4$. Hence we can decompose $H$ with respect to $E_1$ and $E_3$,
which means that $H$ is also of type (III) and $S$ is nonextendable by Section \ref{caseIII}.
\end{proof}

For the proof of Proposition \ref{cor:I-2.1} we will need the following:

\begin{claim}
\label{I-2.3}
If $\beta \geq 5$ and either $M_2$ or $E_1+M_2$ is not quasi-nef, then
$S$ is nonextendable.
\end{claim}

\begin{proof}
We first claim that $E_1.M_2 \geq E.M_2$.

Indeed, if $L_3 \sim E+M_2$, then by the removing conventions in Section \ref{class}, page \pageref{remconv}, 
we must have $E_1.L_3 >E.L_3$, whence $E_1.M_2 \geq E.M_2$. In all other cases we have $L_1 \sim
2E_1+\varepsilon E + M_2$, with $0 \leq \varepsilon \leq 2$, so that $\beta \geq 5$ implies that $\alpha \geq 3$. Therefore
$\varepsilon + 1 + E_1.M_2 =E_1.(L_1+E) > \phi(H) = 2+E.M_2$, whence $E_1.M_2 > 1-\varepsilon+E.M_2 \geq 
E.M_2-1$. This proves our assertion. 

Assume now there is a $\Delta >0$ such that $\Delta^2=-2$ and $\Delta.M_2 \leq -2$. By Lemma \ref{A} we can 
write $M_2 \sim A + k\Delta$, where $A>0, A^2=M_2^2, k=-\Delta.M_2 = \Delta.A \geq 2$ and $A$ is primitive.

If $\Delta.E_1 <0$, then $\Delta.E >0$ since $H$ is ample and by Lemma \ref{nef:main2} we have $E_1 \eqv
E+\Delta$. But this yields $E_1.M_2 = (E+\Delta).M_2 < E.M_2$, a contradiction. Hence $\Delta.E_1 \geq 0$.

If $\Delta.E=0$, then since $H$ is ample, we get $\Delta.E_1 \geq 2$, whence $(E_1+\Delta)^2 \geq 2$. Since
$E.(E_1+\Delta)=1$, we can write $E_1+\Delta \sim (E_1.\Delta-1)E+E'$ for $E' >0$ satisfying $(E')^2=0$ and
$E.E'=1$. From $E.M_2 = \phi(H)-2 \geq \phi(L_1)-2 = E_1.M_2+\beta - \alpha -2 \geq E_1.M_2 + \Delta.M_2 =
(E_1.\Delta-1)E.M_2+ E'.M_2$, we get $E_1.\Delta=2$, $\Delta.M_2 = -2$, $\beta=\alpha$, $E_1.M_2=E.M_2+2$ and
$E'.M_2=0$, whence $M_2 \eqv E'$ by Lemma \ref{nefred}. It follows that $E_1+\Delta \eqv E +
M_2$, whence the absurdity $0 = (E_1+\Delta).\Delta = (E + M_2).\Delta < 0$. Hence $\Delta.E >0$.

Now define $B=\lfloor \frac{\beta}{2} \rfloor E+ E_1 + \Delta  = \lfloor \frac{\beta-2}{2} \rfloor E+ E_1 +
(E+\Delta)$. Then $B^2 \geq 6$, $\phi(B) \geq 2$ and $H-2B \sim \epsilon E + (M_2-2\Delta)$, with
$M_2-2\Delta > 0$, $(M_2-2\Delta)^2 \geq 0$ and $\epsilon =0$ or $1$. By Proposition \ref{ramextenr} we are
immediately done if $B^2 \geq 8$ and if $B^2=6$, we only need to prove that $H^2 \geq 54$.

This is satisfied. Indeed, if $B^2=6$, we must have $\beta=5$, $E.\Delta=1$ and $E_1.\Delta=0$. Note that
$E.M_2=E.A+k\Delta.E = E.A+k \geq 2$, and if $E.M_2=2$, then $k=2, E.A=0$ and $A \eqv E$, whence the
contradiction $2 = A.\Delta = E.\Delta = 1$.  Then $E.M_2 \geq 3$, whence also $E_1.M_2 \geq 3$, so that $H^2 =
(5E+2E_1+M_2)^2 \geq 62$.

Assume similarly that there is a $\Delta >0$ such that $\Delta^2=-2$ and $\Delta.(E_1+M_2) \leq -2$. By what we
have just proved and Lemma \ref{nef:main2}, we can assume that $\Delta.E_1=\Delta.M_2=-1$, but then we get  $E_1
\eqv E+\Delta$, whence $E_1.M_2 = (E+\Delta).M_2 < E.M_2$, a contradiction.
\end{proof}

\begin{prop}
\label{cor:I-2.1}
Let $H$ be of type {\rm (I)} with $\gamma=2$ and $M_2 >0$. Then $S$ is nonextendable if 
$\beta \geq 5$.
\end{prop}

\begin{proof}
By Claim \ref{I-2.3} we can assume that $E_1+M_2$ and $M_2$ are quasi-nef. If there exists a 
nodal curve $\Gamma$ such that $\Gamma.(E+E_1+M_2) <0$, then $\Gamma.E_1 =0$ and $\Gamma.M_2=-1$,
so that $\Gamma.E=0$. But this yields $\Gamma.H=-1$, a contradiction. Therefore $E+E_1+M_2$
is nef.

Set $D_0=kE+E_1+M_2$ with $k = \lfloor \frac{\beta-1}{2} \rfloor \geq 2$. Then $H-D_0
\sim (\beta-k)E+E_1$ and $H-D_0-2E = (\beta - k -2)E + E_1$ are base-component free by Lemma
\ref{nef:main}.

To prove the surjectivity of $\mu_{V_D, \omega_D}$ we apply Lemma \ref{2E}. We have
$h^1(D_0+K_S-2E)=h^1((k-2)E+E_1+M_2+K_S)=0$ by Theorem \ref{lemma:qnef} and
$h^2(D_0+K_S-4E)=h^0((4-k)E-E_1-M_2)=0$ by the nefness of $E$. 

Now $h^0(H-2D_0-2E)=0$ by the nefness of $E$ and the exact sequence
\[ 0  \hpil \O_S(H-2D_0-4E) \hpil \O_S(H-D_0-4E)  \hpil \O_D(H-D_0-4E) \hpil  0, \]
shows that $h^0(\O_D(H-D_0-4E)) \leq h^0(H-D_0-4E) + h^1(H-2D_0-4E)$. 

As $\beta-k-4 \geq -1$ we have $h^0(H-D_0-4E)=h^0((\beta-k-4)E+E_1)=\beta-k-3$ by Lemma
\ref{nef:main} and as $\beta-2k-4 <0$ we have $h^1(H-2D_0-4E) = h^1((\beta-2k-4)E-M_2)=0$ 
by Theorem \ref{lemma:qnef}. Therefore (\ref{green}) holds and $\mu_{V_D, \omega_D}$ is surjective by Lemma
\ref{2E}.

To end the proof we deal with the Gaussian map $\Phi_{H_D, \omega_D}$. By Lemma \ref{lemma:inizioI}
we can assume that $M_2^2 \leq 2$ and using Theorem \ref{lemma:qnef} it can be easily seen that
$h^0(M_2-E) \leq 1$. 

We have $D_0^2= 2k(1+E.M_2)+2E_1.M_2 +M_2^2 \geq 12$, unless $\beta = 5, 6$,  
$E.M_2=E_1.M_2=1$ and $M_2^2=0$, whence $D_0^2=10$. Since $2D_0-H \sim
(2k-\beta)E+M_2$, the map $\Phi_{H_D, \omega_D}$ is surjective by Theorem
\ref{tendian}(c)-(d) unless possibly if $(E.M_2,E_1.M_2,M_2^2)=(1,1,0)$ and $h^0(M_2-E)>0$ if
$\beta = 5$ or $h^0(M_2-2E)>0$ if $\beta = 6$. But in the latter case we have the contradiction $2 =
h^0(2E) \leq h^0(M_2) = 1$. Therefore $S$ is nonextendable by Proposition \ref{mainextenr} except
possibly for the case $(E.M_2,E_1.M_2,M_2^2,\beta)=(1,1,0,5)$ and $h^0(M_2-E)>0$.

This case will be treated in Lemma \ref{lemma:I12} below (where we will set $E_2 = M_2$).
\end{proof}

\begin{lemma} 
\label{lemma:I12}
If $H \sim 5E + 2E_1 + E_2$ with $E.E_1 = E.E_2 = E_1.E_2 = 1$ and $E_2 >E$, then $S$ is nonextendable.
\end{lemma}

To prove the lemma we will use the following two results:

\begin{claim} 
\label{F}
Set $E_0 = E$. Let $F > 0$ be a divisor such that $F^2 = 0$ and $F.E = F.E_1 = F.E_2 = 1$. If $F$ is 
not nef there exists a nodal curve $R$ such that $F \eqv E_i + R$ and $E_i.R = 1$ for some $i \in \{0, 1,
2\}$.
\end{claim}

\begin{proof}
Let $R$ be a nodal curve such that $k:= - F.R \geq 1$. By Lemma \ref{A} we can write $F \sim A + k R$
with $A > 0$ primitive and $A^2 = 0$. Since $H$ is ample there is an $i \in \{0, 1, 2\}$ such that
$E_i.R \geq 1$. Hence $1 = E_i.F = E_i.A + k E_i.R \geq k$, so that $k = 1$, $E_i.A = 0$, 
$E_i.R = 1$ and $A \eqv E_i$ by Lemma \ref{nefred}.
\end{proof}

\begin{claim}
\label{Fgiusto}
There is an isotropic effective $10$-sequence $\{F_1, \ldots, F_{10}\}$ such that $F_1 = E$, $F_2 = E_1$, 
$F_3 = E_2$. For $4 \leq i \leq 10$ set $F'_i = E + E_1 + E_2 - F_i$. Then $F'_i > 0$, $(F'_i)^2 = 0$
and $F'_i.E = F'_i.E_1 = F'_i.E_2 = 1$. Moreover, up to renumbering $F_4, \ldots, F_{10}$, we can
assume that:
\begin{itemize}
\item[(i)] $F_i$ is nef for $7 \leq i \leq 10$.
\item[(ii)] $E + F'_i$ is nef for $9 \leq i \leq 10$.
\item[(iii)] If $E_2 > E$ then $h^0(2F_{10} + E - E_2 + K_S) = 0$.

\end{itemize}
\end{claim}
\begin{proof}
First of all the $10$-sequence exists since by \cite[Cor.2.5.6]{cd} we can complete the isotropic
$3$-sequence $\{E, E_1, E_2\}$ to an isotropic effective $10$-sequence. 

To see (i) suppose that $F_4, \ldots, F_7$ are not nef. By Claim \ref{F} there is an 
$i \in \{0, 1, 2 \}$ and two indices $j, k \in \{4, \ldots, 7 \}$, $j \neq k$ such that 
$F_j \eqv E_i + R_j$ and $F_k \eqv E_i + R_k$. Therefore $R_j.R_k = (F_j - E_i).(F_k - E_i) = -1$, a
contradiction. Upon renumbering we can assume that $F_i$ is nef for $7 \leq i \leq 10$.

Now the definition of $F'_i$ easily gives that $(F'_i)^2 = 0$ and $F'_i.E = F'_i.E_1 = 
F'_i.E_2 = 1$. Also $F'_i > 0$ by Riemann-Roch since $F'_i.E = 1$ implies that $H^2(F'_i) = 0$. To see
(ii) suppose that $E + F'_7$, $E + F'_8$ and $E + F'_9$ are not nef. By Claim \ref{F} there is an $i
\in \{1, 2 \}$ and two indices $j, k \in \{7, 8, 9 \}$, $j \neq k$ such that $F_j' \eqv E_i + R_j$ and
$F_k' \eqv E_i + R_k$, giving a contradiction as above. Upon renumbering we can assume that $E + F'_i$ is nef for
$9 \leq i \leq 10$.

To see (iii) let $F$ be either $F_9$ or $F_{10}$ and suppose that $2F + E - E_2 + K_S 
\geq 0$. Let $\Gamma$ be a nodal component of $E_2 - E$. Since $2F + K_S \geq E_2 - E \geq
\Gamma$ and $h^0(2F + K_S) = 1$, we get that $\Gamma$ must be either a component of $F$ or of
$F + K_S$. Therefore $\Gamma$ is, for example, a component of both $F_9$ and $F_{10}$ and this is not
possible since $F_9.F_{10} = 1$ and they are both nef and primitive.
\end{proof}

\renewcommand{\proofname}{Proof of Lemma {\rm \ref{lemma:I12}}.}  
\begin{proof}
By Claim \ref{Fgiusto}(ii) we know that $E + F'_{10} = 2E + E_1 + E_2 - F_{10}$ is nef, whence, using
\cite[Prop.3.1.6 and Cor.3.1.4]{cd}, we can choose $F \eqv F_{10}$ so that, setting $F' = E + E_1 + E_2 -
F$, we have that $E + F' = 2E + E_1 + E_2 - F$ is a base-component free pencil. Let $D_0 = 3E + E_1 + F$.
Then $D_0^2 = 14$, $\phi(D_0) = 2$ and $D_0$ is nef by Lemma \ref{nef:main} and Claim
\ref{Fgiusto}(i). Now $H - D_0 = 2E + E_1 + E_2 - F = E + F'$ is a base-component free pencil. Also $2D_0 - H 
= 2F + E - E_2$ so that $E.(2D_0 - H) = 1$ and $F.(2D_0 - H) = 0$. If $h^0(2D_0 - H) \geq 2$ we can write $2D_0 -
H \sim G + M$ where $G$ is the base component and $M$ is base-component free. Since both $E$ and $F$ are nef we
get that $E.M \leq 1$ and $F.M = 0$. The latter implies that $M^2 = 0$ whence $M \sim 2hF$ for some $h \geq 1$ by
\cite[Prop.3.1.4]{cd}, but then we get the contradiction $E.M \geq 2$. Therefore $h^0(2D_0 - H) \leq 1$ and
$\Phi_{H_D, \omega_D}$ is surjective by Theorem \ref{tendian}(c)-(d). 

Now $E.(H - 2D_0) = - 1$ whence $h^0(H - 2D_0) = 0$. Also $(H - 2D_0)^2 = -2$ whence, by 
Riemann-Roch, $h^1(H - 2D_0) = h^2(H - 2D_0) = h^0(2D_0 - H + K_S) = h^0(2F + E - E_2 + K_S) = h^0(2F_{10} +
E - E_2 + K_S) = 0$ by Claim \ref{Fgiusto}(iii).

Therefore $\mu_{V_D, \omega_D}$ is surjective by \eqref{pencil} and $S$ is nonextendable by Proposition 
\ref{mainextenr}.
\end{proof}
\renewcommand{\proofname}{Proof}

\section{Remaining cases in (I) with $\gamma=2$ and $M_2 >0$}
\label{caseI-R}

As the cases left have $\beta \leq 4$ by Proposition \ref{cor:I-2.1}, we have, for the whole section,
\[ H \sim \beta E + 2E_1 + M_2, \: \mbox{with } \beta =2, \; 3 \; \mbox{or} \; 4, \]
and either $M_2^2=0$ or we are in one of the cases of Lemma \ref{lemma:inizioI}.

\subsection{The case $M_2^2=0$}

We write $M_2=E_2$ for a primitive $E_2 >0$ with $E_2^2=0$. 

\subsubsection{$\beta=2$}

From (\ref{eq:I-R.2}) and (\ref{eq:I-R.1}) we get $1 \leq E.E_2 \leq E_1.E_2 \leq
E.E_2+2$. Moreover, since $L_1  \sim 2E_1+E_2$, we get $(\phi(L_1))^2= (E_1.E_2)^2 \leq L_1^2 = 4E_1.E_2$, 
whence $E_1.E_2 \leq 3$ by \cite[Prop.1]{kl1}, as $E_2$ is primitive. Since $H^2 \geq 28$, we are left with the
cases $(E.E_2,E_1.E_2)=(2,3)$ or $(3,3)$, so that $S$ is nonextendable by Lemma \ref{subram}(iii-b).
 
\subsubsection{$\beta=3$}

From (\ref{eq:I-R.2}) and (\ref{eq:I-R.1}) we get $1 \leq E.E_2 \leq E_1.E_2+1 \leq 
E.E_2+\alpha$.

If $\alpha=2$ we get $E.E_2-1 \leq E_1.E_2 \leq E.E_2+1$. Moreover, since $L_2
\sim E+E_2$ is of small type, we must have $E.E_2 \leq 3$ or $E.E_2=5$. Furthermore,
since $L_1 \sim E+2E_1+E_2$, we get
$(\phi(L_1))^2=(1+E_1.E_2)^2 \leq L_1^2 = 4+ 4E_1.E_2 + 2E.E_2$. However, in
the case $(E.E_2, E_1.E_2)=(3,4)$, we find $(L_1^2, \phi(L_1))=(26,5)$, which
is impossible by \cite[Prop.1]{kl1}. This yields that $E.E_2 =2,3,5$ if $E_1.E_2=E.E_2-1$; $E.E_2 =1,2,3$ if
$E_1.E_2=E.E_2$; and $E.E_2 =1,2$ if $E_1.E_2=E.E_2+1$.

If $\alpha=3$ we must have, by \eqref{eq:class.1}, that $E_1.(H-3E)= \phi(H)$, whence $E_1.E_2=2+E.E_2$.
Moreover, since $L_1 \sim 2E_1+E_2$, we get $(\phi(L_1))^2=(E_1.E_2)^2 \leq L_1^2 = 4 E_1.E_2$, whence 
$E_1.E_2 \leq 3$ by \cite[Prop.1]{kl1} since $E_2$ is primitive. Hence $E_1.E_2=3$ and $E.E_2=1$.

To summarize, using $H^2 \geq 32$ or $H^2 = 28$, we have the following cases:
\begin{eqnarray} 
\nonumber  E_1.E_2=E.E_2-1, & E.E_2=2, \; 3 \; \mbox{or} \; 5, & g=15, \;  20 \;  \mbox{or} \; 30. 
\\ \label{eq:lista7casi} E_1.E_2=E.E_2, &  E.E_2=2 \; \mbox{or} \; 3, & g=17  \; \mbox{or} \; 22.\\
\nonumber E_1.E_2=3, & E.E_2=2,  & g=19.
\end{eqnarray}

We will now show, in Lemmas \ref{I-R.3}-\ref{I-R.6}, that $S$ is nonextendable in 
the five cases of genus $g \geq 17$.
The case with $g = 15$ is case (b1) \label{b1} in the proof of Proposition \ref{precisa}. 

\begin{lemma} 
\label{I-R.3}
Let $H \sim 3E + 2E_1 + E_2$ be as in \eqref{eq:lista7casi} with $(E.E_2,E_1.E_2,g)=(5,4,30)$. 
Then $S$ is nonextendable.
\end{lemma}

\begin{proof}
We have $H^2= 58$ and $\phi(H)=E.H=E_1.H=7$. Hence both $E$ and $E_1$ are nef by Lemma \ref{qnef0}.
Let now $H' = H-4E$. Then $(H')^2=2$ and consequently we can write
$H \sim 4E+A_1+A_2$ for $A_i >0$ primitive with $A_i^2=0$ and $A_1.A_2=1$. Since $E.H=E.A_1+ E.A_2=7$ we can
assume by symmetry that either (a) $(E.A_1,E.A_2)=(2,5)$ or (b)
$(E.A_1,E.A_2)=(3,4)$. Also since $E_1.H=7$ we have $E_1.(A_1+A_2)=3$,
whence we have the two possibilities $(E_1.A_1,E_1.A_2)=(2,1)$ or $(1,2)$.

In case (b) we get $A_1.H=13$, whence $(H-2(E+A_1))^2 = 2$. Since
$(H-2(E+A_1)).E = 1$, we have $H-2(E+A_1) >0$ by Riemann-Roch, whence
$S$ is nonextendable by Proposition \ref{ramextenr}.

In case (a) we get $A_1.H=9$. Now if $E_1.A_1=2$, we get
$(H-2(E+A_1+E_1))^2=6$, and as above $S$ is nonextendable by Proposition
\ref{ramextenr}.

Hence the only case left is (a) with $(E_1.A_1,E_1.A_2)=(1,2)$. 
Note that $E_1.(H-2E)=A_1.(H-2E)=5$, whence $L_1 \sim H-2E$ and $\phi(L_1)=A_1.L_1=5$.
Therefore we can continue the decomposition with respect to $A_1$ instead of $E_1$. Since
$H$ now is also of type (III), $S$ is nonextendable by Sections \ref{class} and \ref{caseIII}. 
\end{proof}

\begin{claim} 
\label{claimI11} 
Let $H \sim 3E + 2E_1 + E_2$ be as in \eqref{eq:lista7casi} with $(E.E_2,E_1.E_2,g)=(3,2,20)$ 
(respectively $(E.E_2,E_1.E_2,g)=(3,3,22)$). 
Then there exists an isotropic effective $5$-sequence $\{E, F_1, F_2, F_3, F_4\}$ (respectively an isotropic
effective $4$-sequence $\{E, F_1, F_2, F_3\}$ together with an isotropic divisor $F_4>0$ such that
$E.F_4=F_2.F_4=F_3.F_4=1$ and $F_1.F_4=2$) such that $H \sim 2E + 2F_1 + F_2+F_3+F_4$ and:
\begin{itemize}
\item[(a)] $F_1$ is nef and $F_i$ is quasi-nef for $i=2,3,4$;
\item[(b)] $|E+F_2|$ and $|F_1+F_3|$ are without base components;
\item[(c)] $|E+F_1+F_2+F_3|$ and $|E+F_1+F_4|$ are base-point free;
\item[(d)] $h^1(F_1+F_4-F_2) = h^2(F_1+F_4-F_2) =0$. 
\end{itemize}
\end{claim}

\begin{proof}
Since $(E+E_2)^2=6$ and both $E$ and  $E_2$ are primitive, we can write 
$E+E_2 \sim A_1+A_2+A_3$ with $A_i >0$, $A_i^2=0$ and $A_i.A_j=1$ for $i \neq j$. 
We easily find (possibly after renumbering) that $A_i.E=A_i.E_2=A_1.E_1=A_2.E_1=1$ for $i=1,2,3$ and 
$A_3.E_1=1$ if $g=20$ and $2$ if $g=22$. Moreover $A_i.H \leq 8 <2\phi(H)=10$, whence
all the $A_i$'s are quasi-nef by Lemma \ref{qnef0}.

Assume now there is a nodal curve $R_i$ with $R_i.A_i = -1$ for $(i,g) \neq (3,22)$.
Then we can as usual write $A_i \sim B_i + R_i$, with $B_i > 0$ primitive and isotropic. 
Since $A_i.H=6$ we deduce that $B_i \eqv E$ or $B_i \eqv E_1$,
where the latter case only occurs if $g=20$.

If $g=20$, then, since for $i \neq j$, we have $(E+R_i).(E+R_j) = 2 + R_i.R_j = 
(E_1+R_i).(E_1+R_j)$, we see that at most two of the $A_i$'s can be not nef, otherwise we would 
get $R_i.R_j=-1$, a contradiction.  Possibly after reordering the $A_i$'s and adding $K_S$ to two 
of them, we can therefore assume that $A_1$ is nef, and that either $A_2$ is nef or 
$A_2 \sim E+R+K_S$ for $R$ a nodal curve with $E.R=1$.  Now $E_1$ is nef, by Lemma \ref{qnef0}, as $E_1.H
=\phi(H)=5$, so that both $|E_1+A_1|$ and $|E+A_2|$ are without fixed components. Setting $F_1=E_1$, $F_2=A_2$,
$F_3=A_1$ and $F_4=A_3$ we therefore have the desired decomposition satisfying (a) and (b). It also follows by
construction that $E+F_1+F_2+F_3$ and $E+F_1+F_4$ are nef, the latter because $E$ and $F_1$ are, and $F_4$ is
either nef or $F_4 \eqv A+R'$ with $A=E$ or $A=E_1$, for $R'$ a nodal curve with $A.R'=1$. Therefore (c) also 
follows.

If $g=22$, we similarly find that we can assume that $A_1$ and $A_2$ are nef. Moreover
$A_1.L_1=A_1.(H-2E)=E_1.(H-2E)=4$, so if $E_1$ is not nef, we can substitute $E_1$ with $A_1$ and 
repeat the process. Therefore we can assume that $E_1$ is nef as well. Again both  $|E_1+A_1|$ and
$|E+A_2|$ are without fixed components, and  setting $F_1=E_1$, $F_2=A_2$, $F_3=A_1$ and $F_4=A_3$ we
therefore have the desired decomposition satisfying (a) and (b). Now $E+F_1+F_2+F_3$ is again nef by
construction. To see that $E+F_1+F_4$ is nef, assume, to get a contradiction, that there is a nodal curve
$\Gamma$ with $\Gamma.(E+F_1+F_4) <0$. Then $\Gamma.F_4=-1$ and $\Gamma.(E+F_1)=0$ by (a). The ampleness
of $H$ implies $\Gamma.(F_2+F_3) \geq 2$, whence the contradiction $(F_4 -\Gamma)^2=0$ and $(F_4
-\Gamma).(F_2+F_3) \leq 0$, recalling that $F_4 -\Gamma >0$ by Lemma \ref{A}. Therefore (c) is proved.

We now prove (d).

If $g=20$ then $(F_1+F_4-F_2)^2 = -2$ and $(F_1+F_4-F_2).H=5=\phi(H)$, whence $h^2(F_1+F_4-F_2)=0$ and
if $F_1+F_4-F_2>0$ it is a nodal cycle, so that either $h^0(F_1+F_4-F_2)=0$ or $h^0(F_1+F_4-F_2+K_S)=0$. 
Replacing $F_1$ with $F_1+K_S$ if necessary, we can arrange that $h^0(F_1+F_4-F_2)=0$,
whence also $h^1(F_1+F_4-F_2)=0$ by Riemann-Roch.

If $g=22$, then $(F_1+F_4-F_2)^2 = 0$ and $(F_1+F_4-F_2).H=8 < 2\phi(H)$, whence 
(d) follows by Lemma \ref{qnef0} and Theorem \ref{lemma:qnef}.
\end{proof}

\begin{lemma} 
\label{I-R.5}
Let $H \sim 3E + 2E_1 + E_2$ be as in \eqref{eq:lista7casi} with $(E.E_2,E_1.E_2,g)=(3,2,20)$ 
or $(3,3,22)$. Then $S$ is nonextendable.
\end{lemma}

\begin{proof} 
By Claim \ref{claimI11} we can choose $D_0=E+F_1+F_2+F_3$ with $D_0^2=12$ and both $D_0$ and 
$H-D_0 \sim E+F_1+F_4$ base-point free.
 
We have $h^0(2D_0-H)=h^0(F_2+F_3-F_4) \leq 1$ by Lemma \ref{qnef0}, as $(F_2+F_3-F_4).H \leq 6 <2\phi(H)$. 
Hence the map $\Phi_{H_D, \omega_D}$ is surjective by Theorem \ref{tendian}(c)-(d).

To show the surjectivity of $\mu_{V_D, \omega_D}$ we use Claim \ref{claimI11}(b)
and let $D_1 \in |E+F_2|$ and $D_2 \in |F_1+F_3|$ be general smooth curves and apply Lemma 
\ref{multhelp}. 

Now $H-D_0-D_1 \sim F_1+F_4-F_2$ whence $h^1(H-D_0-D_1) = 0$ by Claim \ref{claimI11}(d), so that
$\mu_{V_{D_1}, \omega_{D_1}}$ is surjective by (\ref{mu1}) since $(H-D_0).D_1 = (E+F_1+F_4).(E+F_2)=5$. 

Since $(H-D_0-D_2).H = (E+F_4-F_3).H \leq  7 < 2\phi(H)$ we have that $h^0(H-D_0-D_2) \leq 1$ by Lemma 
\ref{qnef0} and $\mu_{V_{D_2}, \omega_{D_2}(D_1)}$ is surjective by (\ref{mu2}).   
 
Therefore $\mu_{V_D, \omega_D}$ is surjective whence $S$ is nonextendable by Proposition 
\ref{mainextenr}.
\end{proof}

\begin{lemma} 
\label{I-R.6}
Let $H \sim 3E + 2E_1 + E_2$ be as in \eqref{eq:lista7casi} with $E.E_2=2$ and
$(E_1.E_2,g)=(2,17)$ or $(3,19)$. Then $S$ is nonextendable.
\end{lemma}

\begin{proof}
We first observe that it is enough to find an isotropic divisor $F>0$ such that $E.F=1$, $F.H=6$ 
if $g=17$ and $F.H=7$ if $g=19$ and $B:=E+F$ is nef. 

In fact the latter implies that $H \sim 2B+A$, with $A>0$ isotropic with $E.A=2$ and $F.A=4$ if 
$g=17$ and $F.A=5$ if $g=19$. As we assume that $H$ is not $2$-divisible in $\Num S$, $A$ is automatically 
primitive and it follows that $S$ is nonextendable by Lemma \ref{subram}(iii-b).

To find the desired $F$ we first consider the case $g=17$.

Set $Q=E+E_1+E_2$. Then $Q^2=10$ and $\phi(Q)=3$. By \cite[Cor.2.5.5]{cd} there is an 
isotropic effective $10$-sequence $\{f_1, \ldots, f_{10}\}$ with 
\[ 3Q \sim f_1 + \ldots + f_{10}. \] 

Since $E.Q = E_1.Q = 3$, we can assume that $f_1=E$ and $f_2=E_1$ and then $E_2.f_i=1$ for $i \geq 3$.

We now claim that $E+f_i$ is not nef for at most one $i \in \{3, \ldots ,10\}$.

Indeed, note that, for $i \geq 3$, we have $f_i.H = 6 < 2\phi(H) = 8$, whence each $f_i$ is quasi-nef by Lemma
\ref{qnef0}. Now assume that $R_i.(E+f_i) <0$  for some nodal curve $R_i$. Then $R_i.E=0$ and $R_i.f_i=-1$, so
that $f_i \sim \overline{f_i} +R_i$,  by Lemma \ref{A}, with $\overline{f_i}>0$ primitive and 
$\overline{f_i}^2=0$. Since $H$ is ample we  must have $R_i.E_j >0$ for $j=1$ or $2$.

If $R_i.E_2 >0$ then $E_2.f_i=1$ implies $\overline{f_i} \eqv E_2$ and $R_i.E_2=1$.
But then we get the contradiction $E.f_i= E.(E_2+R_i) =2$. 

Therefore $R_i.E_1 >0$, so that $\overline{f_i} \eqv E_1$ and $R_i.E_1=1$. 

Now suppose that also $E+f_j$ is not nef for $j \in \{3, \ldots ,10\}-\{i\}$. Then 
$R_i.R_j = (f_i - E_1).(f_j - E_1) = -1$, a contradiction. Therefore $E+f_i$ is not nef for at most one 
$i \in \{3, \ldots ,10\}$.

Now one easily verifies that any $F \in \{f_3, \ldots ,f_{10}\}$ such that $E+F$ is nef satisfies 
the desired numerical conditions. 

We next consider the case $g=19$.

Since $(E_1 + E_2)^2 = 6$ and $\phi(E_1 + E_2) = 2$ we can find an isotropic effective $3$-sequence 
$\{f_3, f_4, f_5\}$ such that $E_1 + E_2 \sim f_3+f_4+f_5$. Since $E.(E_1 + E_2) = E_1.(E_1 + E_2) = 
3$ we have $f_i.E=f_i.E_1=1$ for $i=3,4,5$, so that we have  an isotropic effective $5$-sequence $\{f_1, \ldots,
f_5\}$ with $f_1=E$ and $f_2=E_1$ such that $H \sim 3f_1 + f_2 + f_3 + f_4 + f_5$.
By \cite[Cor.2.5.6]{cd} we can complete the sequence to an isotropic effective $10$-sequence $\{f_1, \ldots, 
f_{10}\}$. Note that for $i \geq 6$ we have  $f_i.H = 7 < 2\phi(H)=8$, whence each $f_i$ is quasi-nef by Lemma
\ref{qnef0}.

Now the same arguments as above can be used to prove that $E+f_i$ is nef for at least one $i \in \{6, \ldots
,10\}$, whence any $F \in \{f_6, \ldots ,f_{10}\}$ such that $E+F$ is nef satisfies the desired numerical
conditions. 
\end{proof}

\subsubsection{$\beta=4$}

From (\ref{eq:I-R.2}) and (\ref{eq:I-R.1}) we get $1 \leq E.E_2 \leq E_1.E_2+2 \leq 
E.E_2+\alpha$.

If $\alpha=2$ we get $E.E_2-2 \leq E_1.E_2 \leq E.E_2$. Moreover, since $L_2
\sim 2E+E_2$ is not of small type, we get $(\phi(L_2))^2=(E.E_2)^2 \leq L_2^2 = 4E.E_2$,
whence $E.E_2 \leq 3$ by \cite[Prop.1]{kl1}. Therefore $(E.E_2, E_1.E_2) \in
\{(1,1),(2,1),(2,2),(3,1),(3,2),(3,3)\}$. The first case is case (b2) \label{b2} in the proof of Proposition 
\ref{precisa} and in the other cases $S$ is nonextendable by Lemmas \ref{nef:main} and \ref{subram}(iii-a).

If $\alpha=3$ or $4$ we must have $E_1.(H-\alpha E) = \phi(H)$ by \eqref{eq:class.1}, whence
$E_1.E_2=E.E_2+\alpha-2$. Moreover $L_1 \sim (4-\alpha)E+2E_1+E_2$ and using
$(\phi(L_1))^2 \leq L_1^2$, we get $E_1.E_2 \leq 4$. If equality holds then $(L_1^2, \phi(L_1)) = (26, 5)$ or
$(16, 4)$, both excluded by \cite[Prop.1]{kl1}, as $E_2$ is primitive. Therefore 
$(E.E_2, E_1.E_2)=(1,2), (1,3)$ or $(2,3)$ and $S$ is nonextendable by Lemmas \ref{nef:main} and
\ref{subram}(iii-a). 

\subsection{The case $M_2^2=2$}

We write $M_2=E_2+E_3$ for primitive $E_2 >0$ and $E_3 >0$ with $E_2^2=E_3^2 =0$ and
$E_2.E_3=1$, as in Lemma \ref{lemma:inizioI}(a). 

\subsubsection{$\beta=2$}
 
By Lemma \ref{lemma:inizioI} we have left to treat the cases (a-i), that is
\begin{equation} 
\label{eq:poss}
(E.E_2,E.E_3, E_1.E_2, E_1.E_3) =  (1,2,1,2), \; (1,2,2,1), \; (1,1,2,2), \; (2,2,2,2), \; (1,2,2,2).
\end{equation}

We first show that $S$ is nonextendable in the first case of \eqref{eq:poss}. 

Since $E_2.H=\phi(H)=5$ and $E_3.H=9 <2\phi(H)$ we have that $E_2$ is nef and $E_3$ is quasi-nef
by Lemma \ref{qnef0}. In particular we get that $h^1(E_2 + E_3)=h^1(E_2 + E_3+K_S)=0$ by Theorem
\ref{lemma:qnef} and $h^0(E_2 + E_3)=2$ by Riemann-Roch. Now $D_0 := E + E_1 + E_2 + E_3$ is nef by Lemma
\ref{lemmone}(b) with $\phi(D_0) = 3$ and $D_0^2 = 16$.  Also $H-D_0 \sim E+E_1$ is base-component free and 
$2D_0-H \sim E_2 + E_3$. Then $h^0(2D_0-H)=2$ and $h^1(H-2D_0)=0$, so that $\mu_{V_D, \omega_D}$ is surjective 
by \eqref{pencil} and $\Phi_{H_D, \omega_D}$ is surjective by Theorem \ref{tendian}(e), as $\gon(D)=6$ by
\cite[Cor.1]{kl1}, whence $\Cliff(D) = 4$, as $D$ has genus $9$ \cite[\S 5]{elms}. By Proposition
\ref{mainextenr}, $S$ is nonextendable. 

We next show that $S$ is nonextendable in the last four cases in \eqref{eq:poss}. 

By Lemmas \ref{qnef0} and \ref{lemmone}(b) we see that $E_2$ and $E_3$ are quasi-nef and $E+E_1+E_2$ and
$E+E_1+E_3$ are base-point free.

Now set $D_0= E+E_1+E_2$. Then $D_0^2 \geq 8$, $D_0$ is nef, $\phi(D_0) \geq 2$ and $H-D_0 \sim
E+E_1+E_3$ is base-point free. Moreover $h^0(2D_0-H)=0$ as $(2D_0-H).H=(E_2-E_3).H \leq 0$, so that
$\Phi_{H_D, \omega_D}$ is surjective by Theorem \ref{tendian}(c). Now, in all cases except for 
$(E.E_2,E.E_3, E_1.E_2, E_1.E_3) =(1,2,2,2)$, we have $(H-2D_0)^2=-2$ and $(H-2D_0).H=0$, so that 
$h^0(H-2D_0)=h^2(H-2D_0)=0$, whence $h^1(H-2D_0)=0$ by Riemann-Roch and $\mu_{V_D, \omega_D}$ is
surjective by \eqref{biraz} (noting that $(H-D_0)^2 = 10$ in the case $(2,2,2,2)$, while $H-D_0$ is not
$2$-divisible in $\Pic S$ as either $E.(H-D_0)=3$ or $E_1.(H-D_0)=3$ in the other two cases). By Proposition
\ref{mainextenr}, $S$ is nonextendable in those cases.

We now prove the surjectivity of $\mu_{V_D, \omega_D}$ in the case
$(E.E_2,E.E_3, E_1.E_2, E_1.E_3) =(1,2,2,2)$.

Note that $E_1+E_2$ is nef by Lemma \ref{lemmone}(e), whence base-point free, and that $E_1+E_3$ is quasi-nef.
To see the latter, let $\Delta > 0$ be such that $\Delta^2 = -2$ and $\Delta.E_1 + \Delta.E_3 \leq -2$. As $E_1$
is quasi-nef by Lemma \ref{nef:main2} and $E_3$ is quasi-nef we get, again by Lemma \ref{nef:main2}, that
$\Delta.E_1 = \Delta.E_3 = -1$ and $E_1 \eqv E + \Delta$, giving the contradiction $\Delta.E_3 = 0$. Hence
$E_1+E_3$ is quasi-nef.

To show the surjectivity of $\mu_{V_D, \omega_D}$ we let $D_1 = E$ and $D_2 \in
|E_1+E_2|$ be a general smooth curve and apply Lemma \ref{multhelp}. The map $\mu_{V_{D_1}, \omega_{D_1}}$ is
surjective by \eqref{mu1isot} since $h^1(H-D_0-D_1)=h^1(E_1+E_3)=0$ by Theorem \ref{lemma:qnef}.
Finally, $\mu_{V_{D_2}, \omega_{D_2}(D_1)}$ is surjective by \eqref{mu2}, using the fact that
$h^0(H-D_0-D_2)=h^0(E+E_3-E_2) \leq 1$ by Lemma \ref{qnef0}, as $(E+E_3-E_2).H=7 < 2\phi(H)$. Therefore
$\mu_{V_D, \omega_D}$ is surjective and $S$ is nonextendable by Proposition \ref{mainextenr}.
 
\subsubsection{$\beta = 3, 4$}
 
By Lemma \ref{lemma:inizioI} we have left to treat the cases (a-ii) and (a-iii), that is
\begin{eqnarray} 
\label{bb} \beta=3,  \hs (E.E_2, E.E_3, E_1.E_2, E_1.E_3, E_2.E_3) &=& (2,2,2,2,1), \\
\label{aa} \beta=3,  \hs (E.E_2, E.E_3, E_1.E_2, E_1.E_3, E_2.E_3) &=& (2,2,1,2,1), \\
\label{dd} \beta = 3, 4, \hs (E.E_2, E.E_3, E_1.E_2, E_1.E_3, E_2.E_3) &=& (1,1,2,2,1), \\
\label{ee} \beta = 3, 4, \hs (E.E_2, E.E_3, E_1.E_2, E_1.E_3, E_2.E_3) &=& (1,1,1,2,1), \\
\label{ff} \beta = 3, 4, \hs (E.E_2, E.E_3, E_1.E_2, E_1.E_3, E_2.E_3) &=& (1,1,1,1,1).
\end{eqnarray}
 
\begin{claim} 
\label{nef:I-R.1}
In the cases {\rm (\ref{bb})-(\ref{ff})} both $E_2$ and $E_3$ are quasi-nef.
\end{claim}
 
\begin{proof}
We first prove that $E_2$ is quasi-nef. Assume, to get a contradiction, that there exists a $\Delta
>0$ satisfying $\Delta^2 =-2$ and $\Delta.E_2 \leq -2$. Write $E_2 \sim A+k\Delta$, for $A >0$  primitive
with $A^2=0$ and $k=-\Delta.E_2 = \Delta.A \geq 2$. From $E_2.E_3=1$ it follows that 
$\Delta.E_3 \leq 0$. If $\Delta.E >0$, we get from $2 \geq E.E_2 = E.A + k E.\Delta$ that 
$E.E_2=k=2$, $E.\Delta=1$ and $E.A=0$, whence the contradiction $E \eqv A$. Hence $\Delta.E=0$ and 
the ampleness of $H$ gives $\Delta.E_1 \geq 2$, whence the contradiction $E_1.E_2=E_1.A + k E_1.\Delta
\geq 4$. Hence $E_2$ is quasi-nef. The same reasoning works for $E_3$.
\end{proof}

\begin{lemma} 
\label{casesaa-cc,dd-ffbeta>3}
In the cases {\rm (\ref{bb})}-{\rm (\ref{dd})} and in the cases {\rm (\ref{ee})-(\ref{ff})}
with $\beta = 4$ we have that $S$ is nonextendable.
\end{lemma}
 
\begin{proof}
Define $D_0= 2E+E_1+E_2$, which is nef by Lemma \ref{lemmone}(a) with $\phi(D_0) \geq 2$ and 
$D_0^2 \geq 12$ in cases (\ref{bb})-(\ref{dd}) and $D_0^2 = 10$ in cases (\ref{ee}) and (\ref{ff}).

Also $H-D_0 \sim (\beta-2)E+E_1+E_3$, whence $\phi(H-D_0) \geq 2$ and $H-D_0$ is base-point
free by Lemma \ref{lemmone}(b). 

We have $2D_0-H \sim (4-\beta)E+E_2-E_3$, whence $h^0(2D_0-H) \leq 1$ 
in the cases (\ref{bb})-(\ref{dd}), as $(2D_0-H).H \leq \phi(H)$, and $h^0(2D_0-H) =0$ in
cases (\ref{ee})-(\ref{ff}), as $(2D_0-H).H \leq 0$. It follows from Theorem \ref{tendian}(c)-(d) that
the map $\Phi_{H_D, \omega_D}$ is surjective.

We next note that $\mu_{V_D, \omega_D}$ is surjective by \eqref{biraz} if
$h^1(H-2D_0)=h^1(E_3-(4-\beta)E-E_2)=0$. 

Since $(E_3-E_2).H=0$ in cases (\ref{dd}) and (\ref{ff}) we have $h^0(E_3-E_2)=h^2(E_3-E_2)=0$, whence
$h^1(E_3-E_2)=0$ by Riemann-Roch. It follows that $\mu_{V_D, \omega_D}$ is surjective, whence $S$ is
nonextendable by Proposition \ref{mainextenr} in cases (\ref{dd}) and (\ref{ff}) with $\beta=4$. In the
remaining cases we can assume that
\begin{equation} 
\label{eq:h1nonzero}
h^1(E_3-(4-\beta)E-E_2) >0.
\end{equation}

We next show that $\mu_{V_D, \omega_D}$ is surjective in case (\ref{aa}). For this we use Lemmas 
\ref{multhelp}, \ref{nef:main} and \ref{lemmone}(c) and let $D_1 \in |E+E_1|$ and $D_2 \in |E+E_2|$ be general
smooth members.

By Claim \ref{nef:I-R.1} and Theorem \ref{lemma:qnef} we have that $h^1(H-D_0-D_1)=h^1(E_3)=0$, whence 
$\mu_{V_{D_1}, \omega_{D_1}}$ is surjective by \eqref{mu1}. Furthermore $\mu_{V_{D_2}, \omega_{D_2}(D_1)}$ is
surjective by \eqref{mu2}, where one uses that $h^0(H-D_0-D_2)=h^0(E_1+E_3-E_2) \leq 1$ by Lemma \ref{qnef0}
since $(E_1+E_3-E_2).H < 2\phi(H)$. Hence $\mu_{V_D, \omega_D}$ is surjective and $S$ is nonextendable by
Proposition \ref{mainextenr}.

Finally we treat the cases \eqref{bb}, \eqref{dd} (with $\beta=3$) and \eqref{ee} (with $\beta=4$). 
Since $(E_3-(4-\beta)E-E_2)^2=-2$ and $(E_3-(4-\beta)E-E_2).H=-\phi(H)$ in \eqref{bb} and \eqref{dd}
(respectively $2$ in \eqref{ee}), we see that Riemann-Roch and \eqref{eq:h1nonzero} imply
that $E+E_2-E_3+K_S$ is a nodal cycle in \eqref{bb} and \eqref{dd} and $E_3-E_2$ is a nodal cycle in \eqref{ee}.
With $\beta$ as above, it follows that 
\begin{equation} 
\label{eq:h1nonzero1}
h^i(E+E_2-E_3) =0  \; 
\mbox{in (\ref{bb}) and (\ref{dd}) and} \;  
h^i(E_3-E_2+K_S)=0 \; \mbox{in (\ref{ee})}, \; i=0,1,2.  
\end{equation}

We now choose a new $D_0:=(\beta-2)E+E_1+E_3$, which is nef with $\phi(D_0) \geq 2$
and with $H-D_0$ base-point free by Lemma \ref{lemmone}(a) and (b). 
Then $D_0^2 \geq 8$ with $h^0(2D_0-H)=h^0(E_3-E-E_2)=0$ in \eqref{bb} and \eqref{dd}
and $D_0^2=12$ with $h^0(2D_0-H)=h^0(E_3-E_2)=1$ in \eqref{ee}, whence
$\Phi_{H_D, \omega_D}$ is surjective by Theorem \ref{tendian}(c)-(d).

Now \eqref{eq:h1nonzero1} implies $h^1(H-2D_0)=0$, so that $\mu_{V_D, \omega_D}$ is surjective by
\eqref{biraz} and $S$ is nonextendable by Proposition \ref{mainextenr}.  
\end{proof}

We have left the cases (\ref{ee}) and (\ref{ff}) with $\beta=3$, which we treat in Lemmas \ref{lemmaI24} and
\ref{I15}.

\begin{lemma}
\label{lemmaI24}
Suppose $H \sim 3E + 2E_1 + E_2 +E_3$ with $E.E_2=E.E_3=E_1.E_2=E_2.E_3=1$, $E_1.E_3=2$ (the case {\rm 
(\ref{ee})} with $\beta=3$). Then $S$ is nonextendable.
\end{lemma}

\begin{proof}
Since $E_2.H=6$ one easily finds another ladder decomposition 
\begin{equation} 
\label{eq:dec2}
H \sim 3E+2E_2+E_1+E'_3, \; \mbox{with } E_2.E'_3=2,
\end{equation}
and all other intersections equal to one.

We first claim that either $E_1$ or $E_2$ is nef.

In fact $\phi(L_1) = E_1.L_1 = E_1.(E+2E_1+E_2+E_3) = 4 = E_2.L_1$. By Lemma \ref{nef:main2}, if neither $E_1$ 
nor $E_2$ are nef, there are two nodal curves $R_1$ and $R_2$ such that $R_i.E=1$ and $E_i \eqv E+R_i$, for
$i=1,2$. But then we get the absurdity $R_1.R_2= (E_1-E).(E_2-E)=-1$.  

By \eqref{eq:dec2} we can and will from now on assume that we have a ladder
decomposition $H \sim 3E + 2E_1 + E_2 +E_3$ with $E_2$ nef.

\begin{claim} 
\label{cl:11.15b}
Either $h^0(E+E_3-E_2+K_S)=0$, or $h^0(E+E_2-E_3)=0$, or $h^0(E_2+E_3-E+K_S)=0$.
\end{claim}

\begin{proof}
Let $\Delta_1:=E+E_3-E_2+K_S$, $\Delta_2:=E+E_2-E_3$ and $\Delta_3:=E_2+E_3-E+K_S$. Assume, to get a 
contradiction, that $\Delta_i \geq 0$ for all $i=1,2,3$. Since $\Delta_i^2=-2$ we get that 
$\Delta_i > 0$ for all $i=1,2,3$. 

We have $\Delta_2 \sim 2E+K_S -\Delta_1$. Since $\Delta_1.H =6$ and $E.H = 4$, we can neither have 
$\Delta_1 \leq E$ nor $\Delta_1 \leq E+K_S$. Therefore, as $E$ and $E+K_S$ have no common components,
we must have $\Delta_1 = \Delta_{11}+\Delta_{12}$ with $0< \Delta_{11} \leq E$ and $0< \Delta_{12} 
\leq E+K_S$ and $\Delta_{11}.\Delta_{12} =0$. Moreover we have $E.\Delta_{11} = E.\Delta_{12} = 0$,
whence $\Delta_{1i}^2 \leq 0$ for $i = 1, 2$. From $-2=\Delta_1^2=\Delta_{11}^2 +\Delta_{12}^2$ we
must have $\Delta_{1i}^2=0$ either for $i=1$ or for $i = 2$. By symmetry we can assume that
$\Delta_{11}^2=0$. Therefore $\Delta_{11} \eqv qE$ for some $q \geq 1$ by Lemma \ref{nefred}, but
$\Delta_{11} \leq E$, whence $\Delta_{11} = E$ and $\Delta_{12}^2=-2$. Moreover $\Delta_{12}.H=2$. 

Now since $E + \Delta_{12} \eqv \Delta_1 \eqv E+E_3-E_2$, we get $E_3 \eqv E_2 +
\Delta_{12}$ and $E_2.\Delta_{12} = 1$. Hence
$\Delta_3 \sim E_2+E_3-E+K_S \sim (E+E_3+K_S-\Delta_1) + E_3-E+K_S 
\sim 2E_3 - \Delta_1 \sim 2(E_2+ \Delta_{12})-\Delta_1 \sim 2E_2+\Delta_{12} -\Delta_{11}$, therefore
\begin{equation} 
\label{comp}
\Delta_{11} + \Delta_3 \in |2E_2+\Delta_{12}|.
\end{equation}
We claim that $|2E_2+\Delta_{12}|= |2E_2|+\Delta_{12}$. To see the latter observe that it certainly 
holds if $\Delta_{12}$ is irreducible, for then it is a nodal curve with $E_2.\Delta_{12} = 1$ (recall
that $|2E_2|$ is a genus one pencil). On the other hand if $\Delta_{12}$ is reducible then, using
$\Delta_{12}.H=2$ and the ampleness of $H$ we deduce that $\Delta_{12} = R_1 + R_2$ where $R_1, R_2$
are two nodal curves with $R_1.R_2 = 1$. Moreover the nefness of $E_2$ allows us to assume that
$E_2.R_1 = 1$ and $E_2.R_2 = 0$. But then $R_2.(2E_2+\Delta_{12}) = -1$ so that $R_2$ is a
base-component of $|2E_2+\Delta_{12}|$ and of course $R_1$ is a base-component of
$|2E_2+\Delta_{12}-R_2| = |2E_2+R_1|$ and the claim is proved.

Since $\Delta_{11}$ and $\Delta_{12}$ have no common components we deduce from (\ref{comp}) 
that each irreducible component of $E = \Delta_{11}$ must lie in some element of $|2E_2|$. The latter
cannot hold if $E$ is irreducible for then we would have that $2E_2 - E >0$ and $(2E_2-E).E_2 =-1$
would contradict the nefness of $E_2$. Therefore, as is well-known, we have that $E = R_1 + \ldots +
R_n$ is a cycle of nodal curves and we can assume, without loss of generality, that $E_2.R_1 = 1$ and
$E_2.R_i = 0$ for $2 \leq i \leq n$. As we said above, we have $2E_2 - R_1 > 0$. Now for $2 \leq i \leq
n - 1$ we get $R_i.(2E_2 - R_1 - \ldots - R_{i-1}) = -1$, whence $2E_2 - R_1 - \ldots - R_i > 0$.
Therefore $2E_2 - R_1 - \ldots - R_{n-1} > 0$ and since $R_n.(2E_2 - R_1 - \ldots - R_{n-1}) = -2$ we
deduce that $2E_2 - E >0$, again a contradiction.
\end{proof}

{\it Conclusion of the proof of Lemma {\rm \ref{lemmaI24}}}. We divide the proof into the three 
cases of Claim \ref{cl:11.15b}.

{\bf Case A: $h^0(E+E_3-E_2+K_S)=0$.} Set $D_0=2E+E_1+E_3$. Then $D_0^2=12$ and 
$\phi(D_0)=2$. Moreover $D_0$ is nef by Claim \ref{nef:I-R.1} and Lemma 
\ref{lemmone}(a) and $H-D_0 \sim E+E_1+E_2$ is nef since
$E+E_1$ and $E_2$ are (the first by Lemma \ref{nef:main}), so that $|H-D_0|$ is base-point free, since
$\phi(H-D_0)=E.(H-D_0)=2$.

We have $2D_0-H \sim E+E_3-E_2$ and since $(2D_0-H).H =6 < 2\phi(H)=8$, we have 
$h^0(2D_0-H) \leq 1$ by Lemma \ref{qnef0}, so that $\Phi_{H_D, \omega_D}$ is surjective by Theorem
\ref{tendian}(c)-(d).

Clearly $h^0(H-2D_0)=0$ and we also have $h^2(H-2D_0)=h^0(2D_0-H+K_S)=h^0(E+E_3-E_2+K_S)=0$ 
by assumption. Therefore $h^1(H-2D_0)=0$ by Riemann-Roch and $\mu_{V_D, \omega_D}$ is surjective by 
(\ref{biraz}). Hence $S$ is nonextendable by Proposition \ref{mainextenr}.

{\bf Case B: $h^0(E+E_2-E_3)=0$.} We set $D_0= E+E_1+E_3$, so that $D_0^2=8$, $\phi(D_0)=2$ and both $D_0$ and 
$H-D_0 \sim 2E+E_1+E_2$ are nef by Claim \ref{nef:I-R.1} and Lemma \ref{lemmone}(a) and (b), whence base-point 
free. Since $2D_0-H \sim E_3-E-E_2$ and $(E_3-E-E_2).H <0$ we have $h^0(2D_0-H) =0$, whence $\Phi_{H_D,
\omega_D}$ is  surjective by Theorem \ref{tendian}(c).

Now by hypothesis $h^0(H-2D_0)=0$ and we also have $h^0(2D_0-H+K_S) =h^0(E_3-E-E_2+K_S)=0$,
and by Riemann-Roch we get $h^1(H-2D_0)=0$ as well. Therefore $\mu_{V_D, \omega_D}$ is surjective 
by (\ref{biraz}). Hence $S$ is nonextendable by Proposition \ref{mainextenr}.

{\bf Case C: $h^0(E_2+E_3-E+K_S)=0$.} Set $D_0 = E+E_1+E_2+E_3$, which is nef (since $E+E_1+E_3$ is
nef by Claim \ref{nef:I-R.1} and Lemma \ref{lemmone}(b) and $E_2$ is nef by assumption) with $D_0^2=14$ and
$\phi(D_0)=3$. Moreover $H-D_0 \sim 2E+E_1$ is without fixed components.

We have $H-2D_0 \sim E-E_2-E_3$ and since $(H-2D_0).E=-2$ we have $h^0(E-E_2-E_3)=0$. By 
hypothesis we have $h^2(E-E_2-E_3)=0$, whence $h^1(H-2D_0)=0$ by Riemann-Roch. It follows that
$\mu_{V_D, \omega_D}$ is surjective by (\ref{biraz}).

Furthermore, since $2D_0-H \sim E_2+E_3-E$ and $h^0(E_2+E_3-E+K_S)=0$ we have $h^0(2D_0-H) 
\leq 1$, and  $\Phi_{H_D, \omega_D}$ is surjective by Theorem
\ref{tendian}(c)-(d). Hence $S$ is nonextendable by Proposition \ref{mainextenr}.
\end{proof}

\begin{lemma} 
\label{I15}
Suppose $H \sim 3E + 2E_1 + E_2 + E_3$ with $E.E_1 = E.E_2 = E.E_3 = E_1.E_2 = E_1.E_3 = E_2.E_3 = 1$ 
(the case {\rm (\ref{ff})} with $\beta=3$). Then $S$ is nonextendable.
\end{lemma}

\begin{proof} 
By Claim \ref{nef:I-R.1}, Lemma \ref{lemmone}(d) and symmetry, and adding $K_S$ to both $E_2$ and $E_3$ if 
necessary, we can assume that $|E + E_2|$ is base-component free. 

Now set $D_0 =2 E + 2E_1 + E_3$. Then $D_0^2 = 16$ and $\phi(D_0)=3$. Hence Lemmas \ref{nef:main} and 
\ref{lemmone}(b) give that $D_0$ is nef and $H - D_0 \sim E + E_2$ is base-component free. 
 
We have $H - 2D_0 \sim - (2E_1+E+E_3-E_2)$ and we now prove that $h^0(2D_0 - H) = 2$ and $h^1(H - 2D_0)=0$. To 
this end, by Theorem \ref{lemma:qnef} and Riemann-Roch, we just need to show that $B:=2E_1+E+E_3-E_2$ is
quasi-nef.  Let $\Delta > 0$ be such that $\Delta^2 = -2$ and $\Delta.B \leq -2$. By Lemma
\ref{A} we  can write $B \sim B_0 + k\Delta$ where $k = - \Delta.B \geq 2$, $B_0 > 0$ and $B_0^2 = B^2 = 2$. 
Now $2 = E.B = E.B_0 + k E.\Delta \geq 1 + 2E.\Delta$, therefore $E.\Delta = 0$. The ampleness of $H$
implies that $E_2.\Delta \geq 2$, giving the contradiction $4 = E_2.B = E_2.B_0 + k E_2.\Delta \geq 5$.
Therefore $B$ is quasi-nef.

Now let $D \in |D_0|$ be a general curve. By \cite[Cor.1]{kl1} we know that $\gon(D) = 2\phi(D_0) = 6$ whence
$\Cliff(D) = 4$, as $D$ has genus $9$ \cite[\S 5]{elms}. Therefore the map $\Phi_{H_D, \omega_D}$ is
surjective by  Theorem \ref{tendian}(e). Also $\mu_{V_D, \omega_D}$ is surjective by \eqref{pencil} and $S$ is
nonextendable by Proposition \ref{mainextenr}.
\end{proof}

\subsection{The case $M_2^2=4$}

We write $M_2=E_2+E_3$ for primitive $E_2 >0$ and $E_3 >0$ with $E_2^2=E_3^2 =0$ and $E_2.E_3=2$, as
in Lemma \ref{lemma:inizioI}(b).

\subsubsection{$\beta=2$}  
\label{beta2}
By Lemma \ref{lemma:inizioI} we have $(E.E_2,E.E_3)=(1,2)$ and the four cases $(E_1.E_2,E_1.E_3)=(1,2)$,
$(2,1)$, $(2,2)$ and $(1,3)$. Note that in all cases $E_2.H < 2\phi(H)=10$, whence $E_2$ is quasi-nef by Lemma
\ref{qnef0}.

If $(E_1.E_2,E_1.E_3)=(1,2)$ we claim that either $E+E_2$ or $E_1+E_2$ is nef. Indeed if
there is a nodal curve $\Gamma$ such that $\Gamma.(E+E_2)<0$ then $\Gamma.E_2=-1$ and $\Gamma.E=0$. 
By Lemma \ref{lemmone}(a) we have $\Gamma.E_1 >0$, so that $E_2 \eqv E_1+\Gamma$ and 
$E_1+E_2 \eqv 2E_1+\Gamma$ is nef. By symmetry the same arguments work if there is a nodal curve 
$\Gamma$ such that $\Gamma.(E_1+E_2)<0$ and the claim is proved.

By symmetry between $E$ and $E_1$ we can now assume that $E+E_2$ is nef. Setting $A:= H-2E-2E_2$ we have
$A^2=0$. As $E.A=3$ and $E_2.A=4$ we have that $A>0$ is primitive and $S$ is nonextendable by Lemmas
\ref{nef:main} and \ref{subram}(iii-b).

If $(E_1.E_2,E_1.E_3)=(2,1)$ one easily sees that $H \sim 2(E_1+E_2)+A$, with $A^2=0$, $E_1.A=1$ and $E_2.A=4$.
Then $A>0$ is primitive, $E_1+E_2$ is nef by Lemma \ref{lemmone}(e) and $S$ is nonextendable by Lemmas
\ref{nef:main} and \ref{subram}(ii). 

If $(E_1.E_2,E_1.E_3)=(1,3)$ we have $(E_1+E_3)^2=6$ and we can write
$E_1+E_3 \sim A_1+A_2+A_3$ with $A_i >0$, $A_i^2=0$ and $A_i.A_j=1$ for $i \neq j$. Then
$E.A_i=E_1.A_i=E_2.A_i=E_3.A_i=1$ and $A_i.H=6$.

We now claim that either $A_i$ is nef or $A_i \eqv E+\Gamma_i$ for a nodal curve $\Gamma_i$
with $\Gamma_i.E=1$. In particular, at least two of the $A_i$'s are nef.

As a matter of fact if there is a nodal curve $\Gamma$ with $\Gamma.A_i <0$, then since
$A_i.L_1=4=\phi(L_1)$ we must have $\Gamma.L_1 \leq 0$, whence $\Gamma.E >0$ by the ampleness of
$H$ and the first statement immediately follows.
If two of the $A_i$'s are not nef, say $A_1 \eqv E+\Gamma_1$ and $A_2 \eqv E+\Gamma_2$
then $1=A_1.A_2=(E+\Gamma_1).(E+\Gamma_2) = 2 + \Gamma_1.\Gamma_2$ yields the contradiction 
$\Gamma_1.\Gamma_2=-1$ and the claim is proved.

We can therefore assume that $A_1$ and $A_2$ are nef. Let $A=H-2A_1-2A_2$. Then $A^2=0$ and $E.A=1$, whence 
$A>0$ is primitive. As $A_1.A=A_2.A=4$ and $\phi(H)=5$, we have that $S$ is nonextendable by Lemma
\ref{subram}(iii-b). 

If $(E_1.E_2,E_1.E_3)=(2,2)$, note first that $E_1+E_2$ is nef by Lemma \ref{lemmone}(e).
Set $A:= H-2E_1-2E_2$. Then $A^2=0$ and $A.E=1$, so that $A>0$ is primitive. 
As $(E_1+E_2).A=6$, we have that $S$ is nonextendable by Lemma \ref{subram}(ii). 

\subsubsection{$\beta = 3$}

By Lemma \ref{lemma:inizioI} we have $(E.E_2,E.E_3)=(1,2)$ and $(E_1.E_2,E_1.E_3)=(1,3)$ 
or $(2,1)$.

We first show that $E_i$ is quasi-nef for $i=2,3$. We have $H.E_2 \leq 9 < 2\phi(H) = 10$, whence $E_2$ is
quasi-nef by Lemma \ref{qnef0}. Now let $\Delta >0$ be such that $\Delta^2 = -2$ and $\Delta.E_3 \leq -2$. 
By Lemma \ref{A} we can write $E_3 \sim A+ k\Delta$, for $A >0$ primitive with $A^2=0$, $k=-\Delta.E_3 = 
\Delta.A \geq 2$.

If $\Delta.E >0$, from $E.E_3 = E.A +k\Delta.E$ we get that $k=2$, $\Delta.E=1$ and $E.A=0$, whence the
contradiction $E \eqv A$. Hence $\Delta.E=0$.

We get the same contradiction if $\Delta.E_2 >0$. Hence, by the ampleness of $H$ we must have $\Delta.E_1 
\geq 2$, but this gives the contradiction $E_1.E_3 = E_1.A +k\Delta.E_1 \geq 4$. Hence also $E_3$ is quasi-nef.

We now treat the case $(E_1.E_2,E_1.E_3)=(1,3)$.

Let $D_0=2E+E_1+E_2$. Then $D_0^2 =10$ and $\phi(D_0) =2$. Moreover $D_0$ and $H-D_0 \sim E+E_1+E_3$ are 
base-point free by Lemma \ref{lemmone}(a)-(b).
  
Moreover $2D_0-H \sim E+E_2-E_3$, and since $(2D_0-H).E =-1$, we have $h^0(2D_0-H) =0$ and it follows from 
Theorem \ref{tendian}(c) that the map $\Phi_{H_D, \omega_D}$ is surjective.

After possibly adding $K_S$ to both $E_2$ and $E_3$, we can assume, by Lemmas \ref{nef:main2}
and \ref{lemmone}(c), that the general members of both $|E+E_1|$ and $|E+E_2|$ are smooth
irreducible curves. Let $D_1 \in |E+E_1|$ and $D_2 \in |E+E_2|$ be two such curves.

By Theorem \ref{lemma:qnef} we have $h^1(H-D_0-D_1)=h^1(E_3)=0$, whence $\mu_{V_{D_1}, \omega_{D_1}}$ is 
surjective by (\ref{mu1}).

We now claim that $h^0(E_1+E_3-E_2) \leq 2$. Indeed, assume that $h^0(E_1+E_3-E_2) \geq 3$.
Then $|E_1+E_3-E_2| = |M| + G$, with $G$ the base-component and $|M|$ base-component free with $h^0(M)
\geq 3$. If $M^2=0$, then $M \sim lP$, for an elliptic pencil $P$ and an integer $l \geq 2$. But then
$14=(E_1+E_3-E_2).H= (lP+G).H \geq lP.H \geq 4\phi(H)=20$, a contradiction. Hence $M^2 \geq 4$, but
since $M.H \leq (E_1+E_3-E_2).H=14$, this contradicts the Hodge index theorem.

Therefore we have shown that $h^0(E_1+E_3-E_2) \leq 2$ and $\mu_{V_{D_2}, \omega_{D_2}(D_1)}$ is
surjective by (\ref{mu2}). By Lemma \ref{multhelp}, $\mu_{V_D, \omega_D}$ is surjective and by Proposition
\ref{mainextenr}, $S$ is nonextendable.

Next we treat the case $(E_1.E_2. E_1.E_3)=(2,1)$.

Let $D_0=2E+E_1+E_3$. Then $D_0^2 = 14$, $\phi(D_0) =3$ and $D_0$ and $H-D_0 \sim E+E_1+E_2$ are base-point 
free by Lemma \ref{lemmone}(a)-(b).
  
Moreover $2D_0-H \sim E+E_3-E_2$, and since $E+E_3$ is nef by Lemma \ref{lemmone}(c) and $(2D_0-H).(E+E_3) =
(E+E_3-E_2).(E+E_3)=1$, we get that $h^0(2D_0-H) \leq 1$. It follows from Theorem \ref{tendian}(c)-(d) that
the map $\Phi_{H_D, \omega_D}$ is surjective.

Let $D_1 \in |E+E_1|$ and $D_2 \in |E+E_3|$ be two general members.

By Theorem \ref{lemma:qnef} we have that $h^1(H-D_0-D_1)=h^1(E_2)=0$, whence $\mu_{V_{D_1}, \omega_{D_1}}=
\mu_{\O_{D_1}(H-D_0), \omega_{D_1}}$. Since $\omega_{D_1}$ is a base-point free pencil we get that
$\mu_{\O_{D_1(H-D_0)}, \omega_{D_1}}$ is surjective by the base-point free pencil trick because
$\deg(\O_{D_1}(H-D_0-D_1+K_S))=3$, whence
$h^1(\O_{D_1}(H-D_0-D_1+K_S))=0$.

We have $(E_1+E_2-E_3).H =5=\phi(H)$, whence $h^0(E_1+E_2-E_3) \leq 1$ and $\mu_{V_{D_2},
\omega_{D_2}(D_1)}$ is surjective by (\ref{mu2}).

By Lemma \ref{multhelp}, $\mu_{V_D, \omega_D}$ is surjective and, by Proposition \ref{mainextenr}, $S$ is
nonextendable.

\subsection{The case $M_2^2=6$}

By Lemma \ref{lemma:inizioI} we have $\beta=2$ and $M_2=E_2+E_3+E_4$ for primitive $E_i >0$ with $E_i^2=0$,
$E_i.E_j=1$ for $i \neq j$ and $(E.E_2, E.E_3, E.E_4, E_1.E_2, E_1.E_3, E_1.E_4)=(1,1,2,1,1,2)$.
We note that $E_1$, $E_2$ and $E_3$ are nef by Lemma \ref{qnef0} and $E_4$ is quasi-nef by the same lemma.

By the ampleness of $H$ it follows that $D_0:=E+E_1+E_2+E_3+E_4$ is nef with $D_0^2= 24$, $\phi(D_0) = 4$
and $H-D_0 \sim E+E_1$ is base-component free. Since $H-2D_0 \sim -(E_2+E_3+E_4)$
we have $h^1(H-2D_0)=0$ by Theorem \ref{lemma:qnef} and $h^0(2D_0-H)=4$ by Riemann-Roch. Then $\mu_{V_D,
\omega_D}$ is surjective by \eqref{pencil} and $\Phi_{H_D, \omega_D}$ is surjective by Theorem \ref{tendian}(e),
since $\gon(D)=8$ by \cite[Cor.1]{kl1}, whence $\Cliff D=6$ by \cite[\S 5]{elms}, as $g(D)=13$. Hence $S$ is
nonextendable by Proposition \ref{mainextenr}. 

\section{Case (II)}
\label{caseII}

We have
\[ H \eqv \beta E + \gamma E_1 + \delta E_2 + M_3, \hs
E.E_1=E.E_2=E_1.E_2=1, \hs \beta, \gamma, \delta \in \{2,3 \}, \]
$32 \leq H^2 \leq 52$ or $H^2 = 28$.  

Since $M_3$ does not contain $E$, $E_1$ or $E_2$ in its arithmetic genus $1$
decompositions, we have:
\begin{equation} 
\label{eq:II3}
\mbox{ If } M_3 >0 \mbox{ then } E.M_3 \geq \frac{1}{2}M_3^2 +1, \ E_1.M_3 \geq \frac{1}{2}M_3^2 +1 
\mbox{ and } E_2.M_3 \geq \frac{1}{2}M_3^2 +1.
\end{equation}

\begin{claim} 
\label{nef:II1}
$E+E_1+E_2$ is nef.
\end{claim}

\begin{proof}
Let $\Gamma$ be a nodal curve such that $\Gamma.(E+E_1+E_2) <0$. By Lemma \ref{nef:main} we must have
$\Gamma.E_2 <0$. We can then write $E_2=A+k\Gamma$, for $A >0$ primitive with $A^2=0$,
$k=-\Gamma.E_2 = \Gamma.A \geq 1$.

Since $\phi(L_2)=E_2.L_2 = A.L_2 + k\Gamma.L_2 \leq A.L_2$, we must
have $\Gamma.L_2 \leq 0$, whence either $\Gamma.E >0$ or $\Gamma.E_1 >0$,
since $H$ is ample. If $\Gamma.E >0$, then $1=E.E_2= E.A + k\Gamma.E$, whence $k=1$ and $A 
\eqv E$, which means $E_2 \eqv E + \Gamma$. But then $E_1.E_2=1$ yields $\Gamma.E_1=0$,
whence $\Gamma.(E+E_1+E_2)=0$, a contradiction. We get the same contradiction if $\Gamma.E_1 >0$.
\end{proof}

Set $B = E+E_1+E_2$. Then $B^2 = 6$ and $(3B-H).B=18-2(\beta+\gamma+\delta)-(E+E_1+E_2).M_3$. If  
\begin{equation} 
\label{eq:II2}
2(\beta+\gamma+\delta) + (E+E_1+E_2).M_3 \geq 17,
\end{equation}
then $(3B-H).B \leq 1$, whence if $3B-H>0$, it is a nodal cycle by Claim \ref{nef:II1}. Thus \eqref{eq:II2} 
implies that either $h^0(3B-H)=0$ or $h^0(3B+K_S-H)=0$ and $S$ is nonextendable by Proposition \ref{ramextenr6}.

We now deal with \eqref{eq:II2}.

Assume first that $M_3 >0$. Then, in view of (\ref{eq:II3}), the condition
(\ref{eq:II2}) is satisfied unless $\beta = \gamma = \delta =2$ and $(E+E_1+E_2).M_3=3,4$, which means that
$M_3^2=0$, whence $S$ is nonextendable by Claim \ref{nef:II1} and Lemma \ref{subram}(ii). 

Assume now that $M_3 =0$. Then the condition (\ref{eq:II2}) is satisfied unless $6 \leq \beta
+ \gamma + \delta \leq 8$. Since $E.H=\gamma+\delta$ and $E_1.H= \beta +  \delta$, we get $\gamma \leq \beta$.
At the same time, since  $E_1.L_1= \beta-\alpha +  \delta$ and $E_2.L_1=  \beta-\alpha + \gamma$, we get $\gamma
\geq \delta$.  Recalling that we assume that $H$ is not 2-divisible in $\Num S$, we end up with the cases
$(\beta,\gamma,\delta)=(3,2,2)$ or $(3,3,2)$.

The first case has $g=17$ and is case (a3) \label{a3} in the proof of Proposition \ref{precisa}. 
In the second case, set $D_0=2E+E_1+E_2$, which is nef by Claim \ref{nef:II1} and satisfies $D_0^2=10$ and
$\phi(D_0) = 2$. Note that $E_1$ is nef by Lemma \ref{qnef0} since $E_1.H=E.H=\phi(H)$. Now $H-D_0 \eqv
E+2E_1+E_2$ is nef by Claim \ref{nef:II1} with $(H-D_0)^2 = 10$ and $\phi(H-D_0) = 2$, whence $|H-D_0|$ is
base-point free. We have $H-2D_0 \eqv E_1-E$, so that $(H-2D_0)^2=-2$ and $(H-2D_0).H=0$. Therefore, by
Riemann-Roch, $h^i(H-2D_0)=h^i(H-2D_0+K_S)=0$ for all $i=0,1,2$. By Theorem \ref{tendian}(c) 
$\Phi_{H_D, \omega_D}$ is surjective. Moreover $\mu_{V_D, \omega_D}$ is surjective by (\ref{biraz}). Therefore
$S$ is nonextendable by Proposition \ref{mainextenr}.

\section{Case (III)}
\label{caseIII}

We have
\[ H \eqv \beta E + \gamma E_1 + M_2, \hs E.E_1=2, \hs \beta, \gamma \in \{2,3 \}, \]
$32 \leq H^2 \leq 62$ or $H^2 = 28$ and $L_2$ is of small type. In particular
\begin{equation} 
\label{eq:IIIfi}
\phi(H)=E.H=2\gamma+E.M_2 \leq 7
\end{equation}
and
\begin{equation} 
\label{eq:semIII1}
\mbox{either $M_2>0$ or $\beta=\gamma=3$.}
\end{equation}

Since $M_2$ contains neither $E$ nor $E_1$ in its arithmetic genus $1$ decompositions, we have:
\begin{equation} 
\label{eq:III1'}
\mbox{ If } M_2 > 0 \mbox{ then } E.M_2 \geq \frac{1}{2}M_2^2 +1 \mbox{ and } E_1.M_2 \geq \frac{1}{2}M_2^2 + 1.
\end{equation}

By Proposition \ref{ramextenr4} and Lemma \ref{nef:main} we have that $S$ is nonextendable if $(E+E_1).H \geq
17$, therefore in the following we can assume
\begin{equation} 
\label{eq:III2}
(E+E_1).H = 2(\beta +\gamma) +(E+E_1).M_2 \leq 16.
\end{equation}

We now divide the rest of the treatment into the cases $\beta=2$ and $\beta=3$.

\subsection{The case $\beta=2$}
We have $M_2 >0$ by \eqref{eq:semIII1} and $E.M_2 \geq 1$ by \eqref{eq:III1'}.

If $\gamma=3$, then $E.M_2=1$ and $\phi(H)=7$ by \eqref{eq:IIIfi}, so that $M_2^2=0$ 
by \eqref{eq:III1'}. As $L_2 \eqv E_1+M_2$ the removing conventions of Section \ref{class}, page
\pageref{remconv}, require $E_1.L_2 < E.L_2 $. Hence $E_1.M_2 \leq 2$, giving the contradiction
$49 = \phi(H)^2 \leq H^2  \leq 40$.

Therefore $\gamma=2$, so that $E.M_2 \leq 3$ by \eqref{eq:IIIfi}, whence $M_2^2 \leq 4$ by \eqref{eq:III1'}.
Moreover $(E+E_1).M_2 \leq 8$ by \eqref{eq:III2}, whence
\begin{equation} 
\label{eq:semIII4} 
\phi(H)^2 = (4+E.M_2)^2 \leq H^2 =16+M_2^2 + 4(E+E_1).M_2 \leq 48 +M_2^2.
\end{equation}
Combining with \cite[Prop.1]{kl1} we get $E.M_2 \leq 2$, whence $M_2^2 \leq 2$ by \eqref{eq:III1'}. 

We now treat the two cases $M_2^2=0$ and $M_2^2=2$ separately.

If $M_2^2=2$, then $E.M_2=2$ by \eqref{eq:III1'} and since $(E_1.M_2)^2 = \phi(L_1)^2 \leq L_1^2=4E_1.M_2+2$,
we must have $E_1.M_2 \leq 4$. Writing $M_2 \sim E_2+E_3$ for isotropic $E_2>0$ and $E_3>0$ with $E_2.E_3=1$, 
we must have $E.E_2=E.E_3=1$. As $E_i.H \geq \phi(H)=E.H=6$ for $i=2,3$, we find the only possibility
$E_1.E_2=E_1.E_3=2$.

Since $H.E_2=H.E_3 = 7 < 2\phi(H)$ we have that $E_2$ and $E_3$ are quasi-nef by Lemma \ref{qnef0},
whence $E+E_1+E_2$ and $E+E_1+E_3$ are nef by Lemma \ref{lemmone}(b).

Set $D_0=E+E_1+ E_2$. Then $D_0^2=10$, $\phi(D_0)=3$ and both $D_0$ and $H-D_0$ are base-point free.

Now $H-2D_0 \eqv E_3-E_2$ and $(H-2D_0)^2=-2$ with $(H-2D_0).H=0$. Therefore 
$h^i(2D_0-H)=h^i(2D_0-H+K_S)=0$ for $i = 0, 1, 2$. By Theorem \ref{tendian}(c) $\Phi_{H_D, \omega_D}$ is 
surjective. Moreover $\mu_{V_D, \omega_D}$ is surjective by (\ref{biraz}). By Proposition \ref{mainextenr} we
get that $S$ is nonextendable.

Finally, if $M_2^2=0$, then $S$ is nonextendable by Lemmas \ref{nef:main} and \ref{subram}(ii) 
unless $(E+E_1).M_2 \leq 3$. 

In the latter case, by \eqref{eq:semIII4}, we get $E.M_2=1$, whence $M_2^2=0$ by 
\eqref{eq:III1'} and $E_1.M_2 = 2$.

Set $E_2:=M_2$ and grant for the moment the following:

\begin{claim} 
\label{cl:III2i}
There is an isotropic effective $10$-sequence $\{f_1, \ldots, f_{10}\}$, with $f_1=E$, $f_{10}=E_2$, all $f_i$ 
nef for $i \leq 9$, and, for each $i=1, \ldots,9$, there is an effective decomposition $H \sim
2f_i+2g_i+h_i$, where $g_i>0$ and $h_i>0$ are primitive, isotropic with $f_i.g_i=g_i.h_i=2$ and
$f_i.h_i=1$.  Furthermore, $g_i+h_i$ is not nef for at most one $i \in \{1, \ldots, 9\}$.
\end{claim}

By the claim we can assume that $H \sim 2 E + 2 E_1 + E_2$ with $E_1+E_2$ nef. We have $(E_1+E_2-E)^2=-2$.
Since $1=(E_1+E_2).(E_1+E_2-E) < \phi(E_1 + E_2) = 2$ we have that if $E_1+E_2-E>0$ it is a nodal cycle, whence
either $h^0(E_1+E_2-E) =0$ or $h^0(E_1+E_2-E+K_S) =0$. By replacing $E$ with $E+K_S$ if necessary, we can assume 
that $h^0(E_1+E_2-E)=0$. As $h^2(E_1+E_2-E) = h^0(E-E_1-E_2+K_S)=0$ by the nefness of $E$, we find from
Riemann-Roch that $h^1(E_1+E_2-E)=0$ as well. 

Set $D_0 = 2E+E_1$, so that $D_0$ is nef by Lemma \ref{nef:main} with $D_0^2=8$ and $\phi(D_0)=2$.  
Moreover $H-D_0 = E_1+E_2$ is nef by assumption, with $\phi(H-D_0)=2$, whence base-point free.

We have $2D_0-H = 2E-E_2$, and since $(2D_0-H).E =-1$, we have $h^0(2D_0-H)=0$
and by Theorem \ref{tendian}(c) we get that $\Phi_{H_D, \omega_D}$ is surjective.

Now let $D_1 =E$ and $D_2 \in |E+E_1|$ be a general smooth irreducible curve. Since
$h^1(H-D_0-D_1)=h^1(E_1+E_2-E)=0$,  we have that $\mu_{V_{D_1}, \omega_{D_1}}$ is surjective by (\ref{mu1isot}).

Now $h^0(H-D_0-D_2)=h^0(E_2-E) \leq h^0(E_1+E_2-E)=0$, whence $h^0(E_2- 2E) = h^0(E_2-E) 
= 0$, so that $h^1(E_2-2E)=1$ by Riemann-Roch. Therefore 
$h^0(\O_{D_2}(H-D_0-D_1)) = h^0(\O_{D_2}(E_1+E_2-E)) \leq h^0(E_1+E_2-E) +h^1(E_2-2E)=1$
and $\mu_{V_{D_2}, \omega_{D_2}(D_1)}$ is surjective by (\ref{mu2}). By Lemma \ref{multhelp},
$\mu_{V_D, \omega_D}$ is surjective and by Proposition \ref{mainextenr}, $S$ is nonextendable.

We have left to prove the claim.

\renewcommand{\proofname}{Proof of Claim {\rm \ref{cl:III2i}}}

\begin{proof}
Let $Q=E+E_1+E_2$. Then $Q^2=10$ and $\phi(Q)=3$, therefore, by \cite[Cor.2.5.5]{cd}, there is an 
isotropic effective $10$-sequence $\{f_1, \ldots, f_{10}\}$ such that $3Q \sim f_1+ \ldots + f_{10}$. Since
$E.Q=E_2.Q=3$ we can without loss of generality assume that $f_1=E$ and $f_{10}=E_2$. Now $Q \sim
f_1+f_{10}+E_1$, whence $f_i.E_1=1$ for all $i \in \{2, \ldots, 9\}$. It follows that
$f_i.H=\phi(H)=5$ for all $i \leq 9$, whence all $f_i$ are nef for $i \leq 9$ by Lemma \ref{qnef0}.

Now for $i \leq 9$ we have $(H-2f_i)^2=8$. If $\phi(H-2f_i)=1$, then $H-2f_i= 4F_1+F_2$ 
for $F_i>0$, $F_i^2=0$ and $F_1.F_2=1$, but then $f_i.H=5$ implies $f_i.F_1=1$, so that $F_1.H=3$, a
contradiction. Therefore $\phi(H-2f_i)=2$, so that $H-2f_i= 2g_i+h_i$ for isotropic $g_i>0$ 
and $h_i>0$ with $g_i.h_i=2$. Moreover $g_i$ is primitive since it computes $\phi(H-2f_i)$. Now $5 \leq
g_i.H = 2+2f_i.g_i$ implies $f_i.g_i \geq 2$, and $f_i.H=5$ implies $f_i.g_i=2$ and $f_i.h_i=1$, so
that $h_i$ is primitive. Moreover $H.g_i = H.h_i = 6 < 2\phi(H)$, whence $g_i$ and $h_i$ are
quasi-nef by Lemma \ref{qnef0}.

Assume that $g_i+h_i$ is not nef for some $i \leq 9$ and let $R$ be a nodal curve with $R.(g_i+h_i) <0$. 

If $R.g_i <0$ then $R.g_i = -1$ and $R.(2g_i+h_i) \leq -2$, whence $R.f_i \geq 2$ by the ampleness of $H$. 
By Lemma \ref{A} we can write $g_i \sim A + R$, with $A>0$ primitive such that $A^2=0$ and $A.R = 1$. From
$2=f_i.g_i=f_i.A + f_i.R$ we get $f_i.R=2$ and $f_i \eqv A$, a contradiction.

Therefore $R.g_i = 0$, $R.h_i = -1$ and as above we can write $h_i \sim A + R$, with $A>0$ primitive such 
that $A^2=0$ and $A.R = 1$. Now $R.f_i >0$ by the ampleness of $H$, and again by $1=f_i.h_i=f_i.A + f_i.R$ 
we get $f_i.R=1$ and $f_i \eqv A$, so that $h_i \eqv f_i+R$ with $R.f_i=1$.

It follows that if $g_i+h_i$ and $g_j+h_j$ are not nef for two distinct $i,j \leq 9$, say 
for $i=1$ and $j=2$ for simplicity, then $h_1 \eqv f_1+R_1$ and $h_2 \eqv f_2+R_2$ where $R_1$ and $R_2$
are nodal curves such that $R_1.f_1=R_2.f_2=1$. Then
\[ H \eqv 3f_1+2g_1+R_1 \eqv 3f_2+2g_2+R_2. \] 
Now the nefness of $f_2$ and
\[ 5 = f_2.H = 3 + 2g_1.f_2 + R_1.f_2 \]
imply that $g_1.f_2=1$ and $R_1.f_2=0$. As $(R_1+R_2).H =2 < \phi(H)$, we get $R_1.R_2 \leq 1$. Hence
\[ 1 = R_1.H = 3R_1.f_2+2R_1.g_2 +R_1.R_2 = 2R_1.g_2 +R_1.R_2, \]
so that $R_1.g_2=0$ and $R_1.R_2=1$. Similarly $R_2.g_1=0$, whence we get the absurdity
\[ 6 = g_1.H= 3g_1.f_2 + 2g_1.g_2 + g_1.R_2=3 + 2g_1.g_2. \]
Therefore $g_i+h_i$ is not nef for at most one $i \leq 9$ and the claim is proved.
\end{proof}
\renewcommand{\proofname}{Proof}

\subsection{The case $\beta=3$} 
Replacing $E$ with $E+K_S$ if necessary, we can assume that $H \sim 3 E + \gamma E_1 + M_2$. We first claim that
\begin{equation} 
\label{eq:semIII2} 
(\gamma-1)E_1+M_2 \; \mbox{and} \; (\gamma-2)E_1+M_2 \; \mbox{are quasi-nef}. 
\end{equation}

Let $\varepsilon = 0, 1$ and let $\Delta >0$ be such that $\Delta^2=-2$ and 
$\Delta.((\gamma-1-\varepsilon)E_1+M_2) \leq -2$. 

If $\Delta.E_1 <0$, then $\Delta.E \geq 2$ by the ampleness of $H$. By Lemma \ref{A} we can write $E_1 \sim A +
k \Delta$ with $A>0$ primitive, $A^2=0$ and $k = -E_1.\Delta = A.\Delta \geq 1$. But then $2 =
E.E_1 = E.A + k E.\Delta$ implies the contradiction $k=1$, $E.\Delta = 2$ and $E \eqv A$. 

Hence $\Delta.E_1 \geq 0$, so that $M_2 >0$ and $l:= - \Delta.M_2 \geq 2$. By Lemma \ref{A} we can write $M_2
\sim A_2 + l \Delta$ with $A_2 > 0$ primitive, $A_2^2 = M_2^2$ and $\Delta.A_2 = l$.

If $\Delta.E=0$, then $\Delta.E_1 \geq 2$ by ampleness of $H$, whence $E_1.M_2 = E_1.(A_2 + l\Delta) \geq 4$, so
that $\gamma=2$ by \eqref{eq:III2}, which moreover implies $E_1.M_2 \leq 5$, so that $l=E_1.\Delta=2$. As
$(E_1+\Delta)^2=2$, we must have
\[ 2\phi(L_1) \leq (E_1+\Delta).L_1 = \phi(L_1) + \Delta.((3-\alpha)E+2E_1+M_2) = \phi(L_1)+2, \]
and we get the contradiction $4 \leq E_1.M_2 \leq E_1.L_1=\phi(L_1) \leq 2$.

Therefore $\Delta.E >0$, so that $E.M_2 = E.(A_2 + l\Delta) \geq 3$, whence $E.M_2=3$, $\gamma = 2$ and
$\phi(H)=7$ by \eqref{eq:IIIfi}, whence $M_2^2 \leq 4$ by \eqref{eq:III1'}. By \eqref{eq:III2} we must have
$E_1.M_2 \leq 3$, but as $H^2=42+4E_1.M_2+M_2^2 \geq 54$ by \cite[Prop.1]{kl1}, using \eqref{eq:III1'}, we
must have $E_1.M_2=3$. Since $E_1.(H-2E) = 5 \leq \phi(H) = 7$ we have $\alpha = 2$, $L_1 \sim E+2E_1+M_2$ and
$L_2 \sim E+M_2$. Since the latter is of small type and $M_2^2 \leq 4$, we must have either $M_2^2=0$ or
$M_2^2=4$. In the latter case we get $L_2^2 = 10$ and $\phi(L_2) = 3$. Now $(E+\Delta)^2 \geq 0$
and $(E+\Delta).M_2 \leq 1$, whence $\phi(M_2) = 1$ and we can write $M_2 \sim 2F_1 + F_2$ for some
$F_i>0$ with $F_i^2=0$ and $F_1.F_2=1$. Therefore $3 = \phi(L_2) \leq F_1.L_2 = F_1.E + 1$, so that $F_1.E \geq
2$, giving the contradiction $3 = E.M_2 = 2F_1.E + F_2.E \geq 4$. Hence $M_2^2=0$, $L_1^2=26$ and
$\phi(L_1)=E_1.L_1=5$, contradicting \cite[Prop.1]{kl1}. Therefore \eqref{eq:semIII2} is proved. 

To show that $S$ is nonextendable set $D_0=2E+E_1$, which is nef by Lemma \ref{nef:main} with $D_0^2=8$ and 
$\phi(D_0)=2$. Moreover $H-D_0 \sim E+(\gamma-1)E_1+M_2$ is easily seen to be nef by \eqref{eq:semIII2}. Since
$\phi(H-D_0) \geq 2$ we have that $H-D_0$ is base-point free.
 
We have $h^0(2D_0-H)=h^0(E+(2-\gamma)E_1-M_2) =0$ by the nefness of $E$ and \eqref{eq:semIII1},
whence the map $\Phi_{H_D, \omega_D}$ is surjective by Theorem \ref{tendian}(c). 

If $M_2 >0$ and $(\gamma,E.M_2,E_1.M_2)=(2,1,1)$, then $M_2^2=0$ by \eqref{eq:III1'}, $(H-2D_0)^2=-2$ and
$(H-2D_0).H=0$, whence $h^1(H-2D_0) = 0$ by Riemann-Roch, so that $\mu_{V_D, \omega_D}$ is surjective by
\eqref{biraz}, as $E_1$ is primitive.

In the remaining cases, to show the surjectivity of $\mu_{V_D, \omega_D}$ we apply Lemma \ref{multhelp} with
$D_1 = E + K_S$ and $D_2$ general in $|E+E_1+K_S|$.

Since $h^1(H-D_0-D_1) = h^1((\gamma-1)E_1+M_2+K_S)=0$ by \eqref{eq:semIII2} and Theorem \ref{lemma:qnef}, we 
have that $\mu_{V_{D_1}, \omega_{D_1}}$ is surjective by (\ref{mu1isot}).

From \eqref{eq:semIII2} and Theorem \ref{lemma:qnef} we also have $h^1(H-D_0-D_2)=h^1((\gamma-2)E_1+M_2+K_S)=0$,
whence $\mu_{V_{D_2}, \omega_{D_2}(D_1)}= \mu_{\O_{D_2}(H-D_0), \omega_{D_2}(D_1)}$. This is surjective by
\cite[Cor.4.e.4]{gr} if $M_2 >0$, since we assume $(\gamma,E.M_2,E_1.M_2) \neq (2,1,1)$.

Finally, if $M_2=0$, then $\gamma=3$ by \eqref{eq:semIII1}, whence $H.E_1 = 6 = \phi(H)$, so that $E_1$
is nef by Lemma \ref{qnef0} and $h^0(H-D_0-D_2)=h^0(E_1+K_S)=1$. We get 
$h^0(\O_{D_{2}}(H-D_0-D_{1})) \leq h^0( H-D_0-D_1) + h^1(H-2D_0) = 2$, since $h^0(H-D_0-D_1)=h^0(2E_1+K_S)=1$
and $h^1(H-2D_0) = h^1(E_1-E) = 1$ by Riemann-Roch and the nefness of $E$ and $E_1$. Hence $\mu_{V_{D_2},
\omega_{D_2}(D_1)}$ is surjective by \eqref{mu2}.

Therefore $\mu_{V_D, \omega_D}$ is surjective in all cases and $S$ is nonextendable by
Proposition \ref{mainextenr}.

\section{Case (S)}
\label{caseS}

We have $H \sim \alpha E + L_1$ with $L_1^2 >0$ by Lemma \ref{STlemma3} and $L_1$ of small type by hypothesis.
Also we assume that $H$ is not numerically $2$-divisible in $\Num S$ and $H^2 \geq 32$ or $H^2 = 28$.

If $\alpha=2$ we get $H^2 = 4E.L_1+L_1^2 = 4\phi(H) + L_1^2$, whence
$(\phi(H))^2 \leq 4\phi(H)+ L_1^2$ and Lemma \ref{STlemma2} yield $\phi(H) \leq 5$, incompatible with the
hypotheses on $H^2$. Therefore $\alpha \geq 3$ and we can write $L_1 \sim F_1 + \ldots + F_k$ as in Lemma
\ref{STlemma2} with $k = 2$ or $3$ and $E.F_1 \geq \ldots \geq E.F_k$.

If $E.F_k =0$ then $E \eqv F_k$ and $3 \leq \phi(H) = E.H = E.L_1 = F_k.L_1$, so that
$E.L_1 = F_k.L_1 = 3$, $L_1^2 = 10$ and we can write $L_1 \sim E + E_1 + E_2$ with $(E.E_1, E.E_2, E_1.E_2) 
= (1, 2, 2)$.

If $E.F_k > 0$, by definition of $\alpha$ we must have
\[ \phi(H) + 1 \leq F_k.(L_1+E) \leq F_k.L_1 + \frac{1}{k}E.L_1 = F_k.L_1 + \frac{1}{k}\phi(H), \]
whence $F_k.L_1 = 3$ or $4$, $L_1^2 = 10$, $k = 3$ and $\phi(H) = 3$ or $4$.
Hence we can decompose $L_1 \sim E + E_1 + E_2$ with $(E.E_1, E.E_2, E_1.E_2) 
= (1, 2, 2)$ if $\phi(H) = 3$ and $(E.E_1, E.E_2, E_1.E_2)  = (2, 2, 1)$ if $\phi(H) = 4$.
Therefore, setting $\beta = \alpha+1$, we get the following cases: 
\begin{eqnarray}
\label{s1} H \sim \beta E + E_1+E_2, \hs \beta \geq 4, \hs E.E_1=1, E.E_2=E_1.E_2=2, \\
\label{s2} H \sim \beta E + E_1+E_2, \hs \beta \geq 4, \hs E.E_1=E.E_2=2, E_1.E_2=1.
\end{eqnarray}

\begin{claim} 
\label{cl:semplice}
{\rm (i)} In the cases \eqref{s1} and \eqref{s2} we have that $E+E_2$ is nef and $E_2$ is quasi-nef.

{\rm (ii)} In case \eqref{s1} both $nE+E_2-E_1$ and $nE+E_2-E_1+K_S$ are effective and quasi-nef for
all $n \geq 2$, and moreover they are primitive and isotropic for $n=2$. 
\end{claim}

\begin{proof}
Assume that $\Delta >0$ satisfies $\Delta^2=-2$ and $\Delta.E_2 =-k$ for some $k>0$. By Lemma \ref{A} we can
write $E_2 = A+k\Delta$, for $A >0$ primitive with $A^2=0$ and $A.\Delta=k$. If $\Delta.E=0$ the ampleness
of $H$ yields $\Delta.E_1 \geq 2$, and, from $E_1.E_2= E_1.A + kE_1.\Delta$, we get $E_1.E_2=2$,
$k=1$ and $E_1.A=0$, whence the contradiction $E_1 \eqv A$. Therefore $\Delta.E >0$ and it follows
that if $\Delta.(E+E_2) <0$, then $\Delta.E_2 \leq -2$. Hence we can assume $k \geq 2$ and we get
from $2=E.E_2= E.A + kE.\Delta$ that $k=2$, $E.\Delta=1$ and $E.A=0$, whence the contradiction $E
\eqv A$. This proves (i).

As for (ii), note that $(2E+E_2-E_1)^2=0$ and $(E+E_2).(2E+E_2-E_1)=3 < 2\phi(E+E_2)=4$, so that
$h^0(2E+E_2-E_1)=h^0(2E+E_2-E_1+K_S)=1$ by Lemma \ref{qnef0}, whence also $h^1(2E+E_2-E_1)
=h^1(2E+E_2-E_1+K_S)=0$ by Riemann-Roch. Since $E.(2E+E_2-E_1)=1$, the statement follows for $n=2$ by Theorem
\ref{lemma:qnef}, and consequently for all $n \geq 2$ again by the same theorem.
\end{proof}

\begin{lemma} 
Let $H$ be as in \eqref{s1} or \eqref{s2}. Then $S$ is nonextendable.
\end{lemma}

\begin{proof}
We first treat the case \eqref{s1} with $\beta=4$.

In this case we set $D_0=3E+E_2$, which is nef by Claim \ref{cl:semplice}(i) with $D_0^2 = 12$.
Then $H-D_0 \sim E+E_1$ is a base-component free pencil by Lemma \ref{nef:main2} and $H-2D_0 \sim -2E+E_1-E_2$. 
By Claim \ref{cl:semplice}(ii) we have $h^0(2D_0-H)=1$, so that the map $\Phi_{H_D, \omega_D}$ is surjective by
Theorem \ref{tendian}(d), and $h^1(H-2D_0)=0$ so that $\mu_{V_D, \omega_D}$ is surjective by
\eqref{biraz}. By Proposition \ref{mainextenr} we find that $S$ is nonextendable.

In the general case we set $D_0=kE+E_2$ with $k = \lfloor \frac{\beta}{2} \rfloor \geq 2$. 
Then $D_0^2 = 4k \geq 8$ and $D_0$ is nef by  Claim \ref{cl:semplice}(i) with $\phi(D_0)=2$. We have $H-D_0 \sim
(\beta-k)E+E_1$, whence by Lemma \ref{nef:main2} we deduce that $H-D_0$ is base-component free.

Since $2D_0-H \sim (2k-\beta)E+E_2-E_1 \leq E_2-E_1$ we have $h^0(2D_0-H)=0$ as
$(E+E_2).(E_2-E_1)=-1$ in case \eqref{s1} and $H.(E_2-E_1)=0$ in \eqref{s2}. Hence $\Phi_{H_D, \omega_D}$ is
surjective by Theorem \ref{tendian}(c).

Now if $\beta$ is even and we are in case \eqref{s2} we have $h^0(H-2D_0)=h^2(H-2D_0)=0$ as
$H.(H-2D_0)=H.(E_2-E_1)=0$. It follows that $h^1(H-2D_0)=0$ and consequently 
$\mu_{V_D, \omega_D}$ is surjective by \eqref{biraz} since $E_2$ is primitive. Hence $S$ is
nonextendable by Proposition \ref{mainextenr}.

We can therefore assume that $\beta$ is odd in case \eqref{s2}. In particular $\beta \geq 5$ for the
rest of the proof and, by Proposition \ref{mainextenr}, we just need to prove the surjectivity of $\mu_{V_D,
\omega_D}$, for which we will use Lemma \ref{2E}. 

We have $D_0+K_S-2E \sim (k-2)E+E_2+K_S$, whence $h^1(D_0+K_S-2E)=0$ by Claim \ref{cl:semplice}(i) and Theorem
\ref{lemma:qnef}. Moreover $h^2(D_0+K_S-4E)=h^0((4-k)E-E_2)=0$ by the nefness of $E$.

Since $H-D_0-2E \sim (\beta-k-2)E+E_1$ and $\beta-k-2 \geq 1$ as $\beta \geq 5$, we have that 
$|H-D_0-2E|$ is base-component free by Lemma \ref{nef:main2}. Since $(E+E_2).(-E+E_1-E_2) < 0$ we have that
$h^0(H-2D_0-2E)=h^0((\beta-2k-2)E+E_1-E_2) \leq h^0(-E+E_1-E_2)=0$, whence (\ref{green}) is
equivalent to
\begin{equation} 
\label{eq:S.3} 
h^0(\O_D(H-D_0-4E)) \leq (\beta-k-2)E.E_1-1.
\end{equation}
In the case \eqref{s2} with $\beta=5$ we have that $\deg \O_D(H-D_0-4E)= (-E+E_1).(2E+E_2)=3$ and $D$ is
nontrigonal by \cite[Cor.1]{kl1}, whence $h^0(\O_D(H-D_0-4E)) \leq 1$ and (\ref{eq:S.3}) is satisfied. 

Hence we can assume, for the rest of the proof, that $\beta \geq 5$ in case \eqref{s1} and  
$\beta \geq 7$ (and odd) in case \eqref{s2}. This implies $\beta-k-4 \geq -1$ in case \eqref{s1} and
$\geq 0$ in case \eqref{s2}, so that we have $h^0((\beta-k-4)E+E_1)=(\beta-k-4)E.E_1+1$ by Lemma
\ref{nef:main2} and Riemann-Roch. Hence
\begin{eqnarray*}
h^0(\O_D(H-D_0-4E)) & \leq & h^0(H-D_0-4E) + h^1(H-2D_0-4E)) \leq \\ & \leq & (\beta-k-4)E.E_1+
1+h^1(K_S+(2k+4-\beta)E+E_2-E_1),
\end{eqnarray*}
and to prove \eqref{eq:S.3} it remains to show
\begin{equation} 
\label{eq:S.3''} 
h^1(K_S+(2k+4-\beta)E+E_2-E_1) \leq 2E.E_1-2.
\end{equation}

In case \eqref{s1} we have $2k+4-\beta=3$ or $4$, and \eqref{eq:S.3''} follows from Claim
\ref{cl:semplice}(ii).

In case \eqref{s2} we have $2k+4-\beta=3$, and as $(3E+E_2-E_1)^2=-2$ and
$h^2(K_S+3E+E_2-E_1)=h^0(E_1-3E-E_2)=0$, we have that \eqref{eq:S.3''} is equivalent to
$h^0(K_S+3E-E_1+E_2) \leq 2$. If, by contradiction, $h^0(K_S+3E-E_1+E_2) \geq 3$, then we can write
$|K_S+3E-E_1+E_2|=|M|+ \Delta$ for $\Delta $ fixed and $h^0(M) \geq 3$. Since $E.(K_S+3E-E_1+E_2)=0$
and $E$ is nef, we must have $E.M=E.\Delta=0$, whence $M \sim 2lE$ for an integer $l \geq 2$ and
$E_2.\Delta \geq 0$ by the nefness of $E+E_2$. Now $5 = E_2.(K_S+3E-E_1+E_2) \geq 4l \geq 8$, a
contradiction. Hence \eqref{eq:S.3''} is proved.
\end{proof}

\section{Main theorem and surfaces of genus 15 and 17} 
\label{precisazioni}

We have shown, throughout Sections \ref{caseD}-\ref{caseS}, that every Enriques surface $S \subset
\PP^r$ of genus $g \geq 18$ is nonextendable, thus proving our main theorem.

Moreover the theorem can be made more precise in the cases $g=15$ and $g=17$:

\begin{prop}
\label{precisa}
Let $S \subset \PP^r$ be a smooth Enriques surface, let $H$ be its hyperplane bundle let $E>0$ such that 
$E.H = \phi(H)$ and suppose that either $H^2=32$ or $H^2=28$. Then $S$ is nonextendable if $H$ satisfies: 

{\rm (a)} $H^2=32$ and either $\phi(H) \neq 4$ or $\phi(H) = 4$ and neither $H$ nor $H-E$ are $2$-divisible in 
$\Pic S$.

{\rm (b)} $H^2=28$ and either $\phi(H) = 5$ or $(\phi(H), \phi(H-3E))=(4,2)$ or $(\phi(H), \phi(H-4E))=(3,2)$.    
\end{prop} 

\begin{proof}
We have shown, throughout Sections \ref{caseD}-\ref{caseS}, that $S$ is nonextendable unless it has one of the
following ladder decompositions:
\begin{itemize}
\item[(a1)]  $H \sim 4E+4E_1$, $E.E_1=1$, $H^2=32$ (page \pageref{a1});
\item[(a2)] $H \sim 4E+2E_1$, $E.E_1=2$, $H^2=32$ (page \pageref{a2});
\item[(a3)] $H \sim 3E+2E_1+2E_2$, $E.E_1=E.E_2=E_1.E_2=1$, $H^2=32$ (page \pageref{a3}).
\item[(b1)]  $H \sim 3E+2E_1+E_2$, $E.E_1=E_1.E_2=1$, $E.E_2=2$, $H^2=28$ (page \pageref{b1});
\item[(b2)] $H \sim 4E+2E_1+E_2$, $E.E_1=E.E_2=E_1.E_2=1$, $H^2=28$ (page \pageref{b2}).
\end{itemize}

Now in the cases (a1)-(a3) we have $\phi(H)=4$ and we see that $H$ is $2$-divisible in $\Pic S$ in the
cases (a1) and (a2) and $H-E$ is $2$-divisible in $\Pic S$ in case (a3). In case (b1) we have $(\phi(H),
\phi(H-3E))=(4,1)$ whereas in case (b2) we have $(\phi(H), \phi(H-4E))=(3,1)$.
\end{proof}

\section{A new Enriques-Fano threefold}
\label{new3-fold}

We know by the articles of Bayle \cite[Thm.A]{ba} and Sano \cite[Thm.1.1]{sa} that for every $g$ such
that $6 \leq g \leq 10$ or $g = 13$ there is an Enriques-Fano threefold in  $\PP^g$. As mentioned in
the introduction there has been some belief that the examples found by Bayle and Sano
exhaust the complete list of Enriques-Fano threefolds. We will see in this section
that this is not so (see also \cite[Prop.3.2 and Rmk.3.3]{pr2}).

We now prove a more precise version of Proposition \ref{newthreefold}.

\begin{prop} 
\label{prop:newthreefold2}
There exists an Enriques-Fano threefold $X \sub \PP^9$ of genus $9$ with the following properties:
\begin{itemize}
\item[(a)] $X$ does not have a $\QQ$-smoothing. In particular, it does not lie in the closure
of the component of the Hilbert scheme made of Fano-Conte-Murre-Bayle-Sano's examples.
\item[(b)] Let $\mu : \widetilde{X} \khpil X$ be the normalization. Then $\widetilde{X}$ has canonical but not
terminal singularities, it does not have a $\QQ$-smoothing and $(\widetilde{X}, \mu^*\O_X(1))$ does not belong
to the list of Fano-Conte-Murre-Bayle-Sano.
\item[(c)] On the general smooth Enriques surface $S \in |\O_X(1)|$, we have $\O_S(1) \cong
\O_S(2E+2E_1+E_2)$, where $E$, $E_1$ and $E_2$ are smooth irreducible elliptic curves with
$E.E_1=E.E_2=E_1.E_2=1$.
\end{itemize}
\end{prop} 

\begin{proof} 
Let $\overline{X} \subset \PP^{13}$ be the well-known Enriques-Fano threefold of genus $13$. By 
\cite{fa, cm} we have that $\overline{X} \subset \PP^{13}$ is the image of the blow-up of $\PP^3$ along
the edges of a tetrahedron, via the linear system of sextics double along the edges.

This description of $\overline{X}$ allows to identify the linear system embedding
the general hyperplane section $\overline{S} = \overline{X} \cap \overline{H} \subset \PP^{12}$. Let 
$P_1, \ldots, P_4$ be four independent points in $\PP^3$, let $l_{ij}$ be the line joining $P_i$ and
$P_j$ and denote by $\widetilde{\PP}^3$ the blow-up of $\PP^3$ along the $l_{ij}$'s with
exceptional divisors $E_{ij}$ and by $\widetilde{H}$ the pull-back of a plane in $\PP^3$.
Let $\widetilde{L} = 6 \widetilde{H} - 2 \sum\limits_{1 \leq i < j \leq 4} E_{ij}$. Therefore 
$\overline{S}$ is just a general element $\widetilde{S} \in |\widetilde{L}|$, embedded with
$\widetilde{L}_{|\widetilde{S}}$. Now let $\widetilde{l}_{ij}$ be the inverse image of $l_{ij}$ on
$\widetilde{S}$. Then by \cite[Ch.4, \S 6, page 634]{gh}, for each pair of disjoint lines $l_{ij},
l_{kl}$ on $\widetilde{S}$ there is a genus one pencil $|2 \widetilde{H}_{|\widetilde{S}} -
\widetilde{l}_{ik} - \widetilde{l}_{il} - \widetilde{l}_{jk} - \widetilde{l}_{jl}| = |2
\widetilde{l}_{ij}|$.  Therefore 
$\widetilde{L}_{|\widetilde{S}} \sim 2 \widetilde{l}_{12} + 2 \widetilde{l}_{13} + 
2 \widetilde{l}_{14}$
and we have decomposed the hyperplane bundle of $\overline{S} \subset \PP^{12}$ as $2E + 2E_1 + 2
E_2$ where $E := \widetilde{l}_{12}, E_1 := \widetilde{l}_{13}, E_2 := \widetilde{l}_{14}$
are half-pencils and $E.E_1 = E.E_2 = E_1.E_2 = 1$. Also $E, E_1$ and $E_2$ are smooth and irreducible.

To find a new Enriques-Fano threefold $X$ of genus 9 we consider the linear span $M \cong \PP^3$ of
$E_2$, the projection $\pi_M : \PP^{13} - - \to \PP^{9}$ and let $X = \pi_M (\overline{X}) \subset
\PP^{9}$. 

Let $\psi:\widetilde{\overline{X}} \khpil \overline{X}$ be the blow up of $\overline{X}$ along $E_2$ 
with exceptional divisor $F$ and set $\H = (\psi^* \O_{\overline{X}}(1))(-F)$ and let
$\widetilde{\overline{S}} \in |\H| \cong |\fake_{E_2 / \overline{X}}(1)|$ be the smooth Enriques
surface isomorphic to $\overline{S}$. 

Then, by 
\[ 0 \hpil \O_{\widetilde{\overline{X}}} \hpil \H \hpil \H_{|\widetilde{\overline{S}}} \hpil 0, \]
and the fact that $\O_{\widetilde{\overline{S}}}(\H) \cong \O_{\overline{S}}(2E+2E_1+E_2)$ is very ample
(which one can easily verify using the fact that $2E + 2E_1 + 2 E_2$ is very ample and
\cite[Cor.2, page 283]{cd}) and $h^1(\O_{\widetilde{\overline{X}}})=h^1(\O_{\overline{X}}) = 0$, we see that
$|\H|$ is base-point free and thus defines a morphism $\varphi_{\H}$ such that $X=
\varphi_{\H}(\widetilde{\overline{X}}) \sub \PP^9$. Note that $\H^3 = (2E+2E_1+E_2)^2 = 16$, whence
$X$ is a threefold.

Let us see now that $X$ is not a cone over its general hyperplane section $S :=
\psi(\widetilde{\overline{S}})$. 

Consider the four planes $H_1, ..., H_4$ in $\PP^3$ defined by the faces of the
tetrahedron. As any sextic hypersurface in $\PP^3$ that is double on the edges of the tetrahedron
and goes through another point of $H_i$ must contain $H_i$, we see that these four planes are
contracted to four singular points $Q_1, \ldots, Q_4 \in \overline{X}$. Moreover their linear span 
$<Q_1, \ldots, Q_4>$ in $\PP^{13}$ has dimension 3, since the hyperplanes containing $Q_1, \ldots,
Q_4$ correspond to sextics in $\PP^3$ containing $H_1, \ldots, H_4$. Now suppose that $X$ is a
cone with vertex $V$. Then $Q_1, \ldots, Q_4$ project to $V$, whence $\dim <M, Q_1, \ldots,
Q_4> \leq 4$ and $\dim M \cap <Q_1, \ldots, Q_4> \geq 2$. On the other hand we know that
$M = <E_2> \subset \overline{H}$, where $\overline{H}$ is a general hyperplane. Therefore we have that $Q_i
\not\in \overline{H}, 1 \leq i \leq 4$, whence $\dim \overline{H} \cap <Q_1, \ldots, Q_4> = \dim M \cap 
<Q_1, \ldots, Q_4> = 2$, so that $\overline{H} \cap <Q_1, \ldots, Q_4> = M \cap <Q_1, \ldots, Q_4>$. 
Now choose the projection from $M' = <E_1> \subset \overline{H}$. If also $\pi_{M'}(\overline{X})$ is 
a cone then, by the same argument above, we get $\overline{H} \cap <Q_1, \ldots, Q_4> = M' \cap 
<Q_1, \ldots, Q_4>$ and therefore $\dim M \cap M' \geq 2$. But this is absurd since 
$\dim M \cap M' = 6 - \dim <E_1 \cup E_2> = 6 - h^0(\O_{\overline{S}}(2E+E_1+E_2)) = 0.$
Therefore, $X$ is an Enriques-Fano threefold satisfying (c).

Now let $X'$ be the only threefold in $\PP^9$ appearing in Bayle-Sano's list, namely an
embedding, by a line bundle $L'$, of a quotient by an involution of a smooth complete
intersection $Z$ of two quadrics in $\PP^5$. Let $S'$ be a general hyperplane section of $X'$. We claim 
that the hyperplane bundle $L'_{|S'}$ is 2-divisible in $\Num S'$. As $2E+2E_1+E_2$ is not 2-divisible 
in $\Num S$, this shows in particular that $X$ does not belong to the list of Bayle-Sano.

By \cite[\S 3, page 11]{ba}, if we let $\pi: Z \to X'$ be the quotient map, we have that 
$- K_Z = \pi^{\ast}(L')$ and the K3 cover $\pi_{|S''} : S'' \to S'$ is an anticanonical surface in $Z$, 
that is a smooth complete intersection $S''$ of three quadrics in $\PP^5$. Therefore, if $H_Z$ is the 
line bundle giving the embedding of $Z$ in $\PP^5$, we have $- K_Z = 2 H_Z$, whence, setting 
$p = \pi_{|S''}$, $H_{S''} = (H_Z)_{|S''}$, we deduce that $p^{\ast}(L'_{|S'}) \cong 
(\pi^{\ast}L')_{|S''} = 2 H_{S''}$. Suppose now that $L'_{|S'}$ is not 2-divisible in $\Num S'$. Then
$(L'_{|S'})^2 = 16$ and by \cite[Prop.1]{kl1} we have that $\phi(L'_{|S'}) = 3$ and it is easily seen
that there are three isotropic effective divisors $E, E_1, E_2$ such that either (i) $L'_{|S'} \sim  2E
+ 2E_1 + E_2$ with $E.E_1 = E.E_2 = E_1.E_2 = 1$ or (ii) $L'_{|S'} \sim  2E + E_1 + E_2$ with $E.E_1 =
1$, $E.E_2 = E_1.E_2 = 2$. In case (i) we get that $p^{\ast}(E_2) \sim 2D$, for some $D \in \Pic S''$.
Since $(p^{\ast}(E_2))^2 = 0$, we have $D^2 = 0$ and, as we are on a K3 surface, either $D$ or $-D$ is
effective. Also $4 H_{S''}.D = p^{\ast}(L'_{|S'}).p^{\ast}(E_2) = 8$, therefore $H_{S''}.D = 2$ and $D$
is a conic of arithmetic genus $1$, a contradiction. In case (ii) we get that $p^{\ast}(E_1 + E_2) \sim
2D'$, for some $D' \in \Pic S''$ with $(D')^2 = 2$ and $H_{S''}.D' = 5$. But now $|D'|$ cuts out a
$g^2_5$ on the general element $C \in |H_{S''}|$ and this is a contradiction since $C$ is a smooth
complete intersection of three quadrics in $\PP^4$. Therefore $L'_{|S'}$ is 2-divisible in $\Num S'$.

Now assume that $X$ has a $\QQ$-smoothing, that is (\cite{mi}, \cite{re1}) a small deformation 
$\X \hpil \Delta$ over the 1-parameter unit disk, such that, if we denote a fiber by $X_t$, we have that
$X_0 = X$ and $X_t$ has only cyclic quotient terminal singularities. Let $L = \O_X(1)$. We have that
$H^1(N_{S/X_0}) = H^1(\O_S(1)) = 0$, whence the Enriques surface $S$ deforms with any deformation of
$X_0$. Therefore  we can assume, after restricting $\Delta$ if necessary, that there is an $\L \in \Pic \X$ such
that $h^0(\L) > 0$ and $\L_{|X} = L$ (this also follows from the proof of \cite[Thm.5]{ho}, since
$H^1(T_{{\PP^9}_{|X}}) = 0$). Taking a general element of $|\L|$ we therefore obtain a family $\Ss \hpil
\Delta$ of surfaces whose fibers $S_t$ belong to $|L_t|$, where $L_t := \L_{|X_t}$ and $S_0 = S \in |L|$ is
general, whence a smooth Enriques surface with hyperplane bundle $H_0 := L_{|S_0} \sim 2E+2E_1+E_2$ of
type (i) above. Therefore, after restricting $\Delta$ if necessary, we can also assume that the general
fiber $S_t$ is a smooth Enriques surface ample in $X_t$, so that $(X_t, S_t)$ belongs to the list of
Bayle \cite[Thm.B]{ba} and is therefore a threefold like $X' \subset \PP^9$.

Let $H_t = (L_t)_{|S_t}$. As we saw above, we have $H_t \eqv 2A_t$, for some $A_t \in \Pic S_t$. This
must then also hold at the limit, so that $H_0 \sim 2E+2E_1+E_2 \eqv 2A_0$, for some $A_0 \in \Pic
S_0$. But then $E_2$ would be $2$-divisible in $\Num S_0$, a contradiction.

We have therefore shown that $X$ does not have a $\QQ$-smoothing. In particular it does not lie in the 
closure of the component of the Hilbert scheme consisting of Enriques-Fano threefolds with only cyclic
quotient terminal singularities (the fact that such threefolds do fill up a component of the Hilbert
scheme is a simple consequence of the fact that one can globalize, on a family, the construction of the
canonical cover \cite[Proof of Thm.4.2]{mi}, \cite[5.3]{km}). Hence (a) is proved.

To see (b) note that $\widetilde{\overline{X}}$ is terminal (because $\overline{X}$ is), whence the
morphism $\varphi_{\H}$ factorizes through $\widetilde{X}$. Since $\widetilde{X}$ is $\QQ$-Gorenstein by
\cite{ch}, an easy calculation, using a common resolution of singularities of $\widetilde{\overline{X}}$ and
$\widetilde{X}$ and the facts that $-K_{\widetilde{X}} \eqv \mu^*\O_X(1)$ and $-K_{\overline{X}} \eqv
\O_{\overline{X}}(1)$, shows that $\widetilde{X}$ is canonical. 

Finally, the same proof as above shows that $(\widetilde{X}, \mu^* \O_X(1))$ does not belong to the list
of Fano-Conte-Murre-Bayle-Sano and that $\widetilde{X}$ does not have a $\QQ$-smoothing. Hence by
\cite[MainThm.2]{mi} $\widetilde{X}$ cannot be terminal. This proves (b). 
\end{proof}
\renewcommand{\proofname}{Proof}

\begin{rem} 
\label{remko}
{\rm Since $\widetilde{X}$ has canonical but not terminal singularities, the morphism $\varphi_{\H}$ in the
proof of Proposition \ref{prop:newthreefold2} must in fact contract divisors, for otherwise $\widetilde{X}$
would be terminal (as $\widetilde{\overline{X}}$ is). This contraction makes $\widetilde{X}$ acquire new
singularities. It would be interesting to understand how these singularities affect the non existence of a
$\QQ$-smoothing. Moreover we observe that $\widetilde{X}$ is a $\QQ$-Fano threefold of Fano index $1$ with
canonical singularities not having a $\QQ$-smoothing, thus showing that Minagawa's theorem \cite[MainThm.2]{mi}
cannot be extended to the canonical case.}
\end{rem}

\begin{rem} 
\label{rem2}
{\rm Somehow Proposition \ref{prop:newthreefold2}(c) shows the spirit of the method of classification we 
introduce in this paper. The question of existence of threefolds is reduced to the geometry of decompositions
of the hyperplane bundle of the surface sections. In fact, in the case of Enriques surfaces, we can write down
all ``decomposition types'' of hyperplane bundles of genus $g \leq 17$. In each case one can try to either show
nonextendability or to find a threefold with that particular hyperplane section, whence either get a new one or
one belonging to the list of Bayle-Sano. For instance, Prokhorov's new Enriques-Fano threefold of genus $17$ must
belong to one of the three cases (a1)-(a3) of Proposition \ref{precisa}. Once one proves existence one can use
the same construction method as in the proof of Proposition \ref{prop:newthreefold2} and project down to find new
Enriques-Fano threefolds. We also observe that our method shows that, in several ``decomposition types'' of
hyperplane bundles of genus $g \leq 17$, we can prove that the Enriques-Fano threefold is not itself hyperplane
section of some fourfold and that its general Enriques surface section must contain rational curves.}
\end{rem}


\begin{thebibliography}{[CCML]} 

\bibitem[AS]{as} E.~Arbarello, E.~Sernesi. 
\textit{Petri's approach to the study of the ideal associated to a  special divisor}. 
Invent. Math. \textbf{49}, (1978) 99--119. 

\bibitem[Ba]{ba} L.~Bayle.
\textit{Classification des vari{\'e}t{\'e}s complexes projectives de dimension trois dont une
section hyperplane g{\'e}n{\'e}rale est une surface d'Enriques}.
J. Reine Angew. Math. \textbf{449}, (1994) 9--63. 

\bibitem[BCP]{bcp} I.~Bauer, F.~Catanese, R.~Pignatelli. 
\textit{Complex surfaces of general type: some recent progress}.
ArXiv alg-geom math.AG/0602477.

\bibitem[Be]{be} A.~Beauville. 
\textit{Fano threefolds and K3 surfaces}. 
Proceedings of the Fano Conference, Dipartimento di Matematica, Universit\`a di Torino, (2004)
175--184.

\bibitem[BEL]{bel} A.~Bertram, L.~Ein, R.~Lazarsfeld. 
\textit{Surjectivity of Gaussian maps for line bundles of large degree on curves}. 
Algebraic geometry (Chicago, IL, 1989), 15--25, Lecture Notes in Math. \textbf{1479}. Springer,
Berlin, 1991.

\bibitem[Bo]{bo} E.~Bombieri. 
\textit{Canonical models of surfaces of general type}. 
Inst. Hautes \'Etudes Sci. Publ. Math. \textbf{42}, (1973) 171--219.

\bibitem[BPV]{bpv} W.~Barth, C.~Peters, A.~van de Ven. 
\textit{Compact complex surfaces}. 
Ergebnisse der Mathematik und ihrer Grenzgebiete \textbf{4}.
Springer-Verlag, Berlin-New York, 1984.

\bibitem[BS]{bs} M.~C.~Beltrametti, A.~J.~Sommese. 
\textit{The adjunction theory of complex projective varieties}. 
de Gruyter Expositions in Mathematics \textbf{16}.
Walter de Gruyter \& Co., Berlin, 1995.

\bibitem[CD]{cd} F.~R.~Cossec, I.~V.~Dolgachev. 
\textit{Enriques Surfaces I}.
Progress in Mathematics \textbf{76}. Birkh\"auser Boston, MA, 1989.

\bibitem[Ch]{ch} I.~A.~Cheltsov. 
\textit{Singularity of three-dimensional manifolds possessing an ample effective divisor---a
smooth surface of Kodaira dimension zero}.
Mat. Zametki \textbf{59}, (1996) 618--626, 640; translation in Math. Notes  \textbf{59}, (1996) 445--450.

\bibitem[CHM]{chm} C.~Ciliberto, J.~Harris, R.~Miranda. 
\textit{On the surjectivity of the Wahl map}.
Duke Math. J. \textbf{57}, (1988) 829--858.

\bibitem[CLM1]{clm1} C.~Ciliberto, A.~F.~Lopez, R.~Miranda. 
\textit{Projective degenerations of $K3$ surfaces, Gaussian maps, and Fano threefolds}.
Invent. Math. \textbf{114}, (1993) 641--667.

\bibitem[CLM2]{clm2}  C.~Ciliberto, A.~F.~Lopez, R.~Miranda.
\textit{Classification of varieties with canonical curve section via Gaussian maps on canonical
curves}.
Amer. J. Math. \textbf{120}, (1998) 1--21.

\bibitem[CM]{cm} A.~Conte, J.~P.~Murre. 
\textit{Algebraic varieties of dimension three whose hyperplane sections are Enriques surfaces}.
Ann. Scuola Norm. Sup. Pisa Cl. Sci. \textbf{12}, (1985) 43--80.

\bibitem[Co]{cos} F.~Cossec. 
\textit{On the Picard group of Enriques surfaces}. 
Math. Ann. \textbf{271}, (1985) 577--600.

\bibitem[ELMS]{elms} D.~Eisenbud, H.~Lange, G.~Martens, F-O.~Schreyer. 
\textit{The Clifford dimension of a projective curve}.
Compositio Math. \textbf{72}, (1989) 173--204.

\bibitem[Fa]{fa} G.~Fano.
\textit{Sulle variet{\`a} algebriche a tre dimensioni le cui sezioni iperpiane sono 
superficie di genere zero e bigenere uno}. 
Memorie Soc. dei XL \textbf{24}, (1938) 41--66.

\bibitem[GH]{gh} P.~Griffiths, J.~Harris. 
\textit{Principles of algebraic geometry}.
Wiley Classics Library. John Wiley \& Sons, Inc., New York, 1994.

\bibitem[GLM]{glm} L.~Giraldo, A.~F.~Lopez, R.~Mu\~{n}oz. 
\textit{On the existence of Enriques-Fano threefolds of index greater than one}.
J. Algebraic Geom. \textbf{13}, (2004) 143--166. 

\bibitem[Gr]{gr} M.~Green. 
\textit{Koszul cohomology and the geometry of projective varieties}. 
J. Differ. Geom. \textbf{19}, (1984) 125--171.

\bibitem [H]{ho} E.~Horikawa. 
\textit{On deformations of rational maps}.
J. Fac. Sci. Univ. Tokyo Sect. IA Math. \textbf{23} (1976) 581--600.

\bibitem [I1]{isk1} V.~A.~Iskovskih. 
\textit{Fano threefolds. I}.
Izv. Akad. Nauk SSSR Ser. Mat. \textbf{41}, (1977) 516--562, 717.

\bibitem [I2]{isk2} V.~A.~Iskovskih. 
\textit{Fano threefolds. II}.
Izv. Akad. Nauk SSSR Ser. Mat. \textbf{42}, (1978) 506--549.

\bibitem[JPR]{jpr} P.~Jahnke, T.~Peternell, I.~Radloff.
\textit{Threefolds with big and nef anticanonical bundles I}.
Math. Ann. \textbf{333}, (2005) 569--631.

\bibitem[KL1]{klvan} A.~L.~Knutsen, A.~F.~Lopez. 
\textit{A sharp vanishing theorem for line bundles on K3 or Enriques surfaces}. 
Preprint 2005.

\bibitem[KL2]{kl1} A.~L.~Knutsen, A.~F.~Lopez. 
\textit{Brill-Noether theory of curves on Enriques surfaces I}. 
Preprint 2006.

\bibitem[KL3]{klgm} A.~L.~Knutsen, A.~F.~Lopez. 
\textit{Surjectivity of Gaussian maps for curves on Enriques surfaces}. 
Preprint 2005.

\bibitem[KM]{km} J.~Koll\'ar, S.~Mori.
\textit{Birational geometry of algebraic varieties}. 
With the collaboration of C. H. Clemens and A. Corti. 
Cambridge Tracts in Mathematics \textbf{134}. Cambridge University Press, Cambridge, 1998.

\bibitem[Lv]{lv} S.~L'vovsky. 
\textit{Extensions of projective varieties and deformations. I}.
Michigan Math. J. \textbf{39}, (1992) 41-51.

\bibitem[Mi]{mi} T.~Minagawa. 
\textit{Deformations of $\QQ$-Calabi-Yau $3$-folds and $\QQ$-Fano $3$-folds of Fano index $1$}.  
J. Math. Sci. Univ. Tokyo \textbf{6}, (1999) 397--414.

\bibitem[MM]{mm} S.~Mori, S.~Mukai. 
\textit{Classification of Fano $3$-folds with $B_2 \geq 2$}.
Manuscripta Math. \textbf{36}, (1981/82) 147--162. 
\textit{Erratum: "Classification of Fano 3-folds with $B_2 \geq 2$"}.  
Manuscripta Math. \textbf{110}, (2003) 407.

\bibitem[Muk1]{muk1} S.~Mukai. 
\textit{New developments in the theory of Fano threefolds: vector bundle method and moduli
problems}.  
Sugaku Expositions \textbf{15}, (2002) 125--150.

\bibitem[Muk2]{muk2} S.~Mukai. 
\textit{Biregular classification of Fano $3$-folds and Fano manifolds of coindex $3$}.  
Proc. Nat. Acad. Sci. U.S.A. \textbf{86}, (1989) 3000--3002.

\bibitem[P1]{pr1} Yu.~G.~Prokhorov.
\textit{The degree of Fano threefolds with canonical Gorenstein singularities}.
Mat. Sb. \textbf{196}, (2005) 81--122; translation in Sb. Math. \textbf{196}, (2005) 77--114.

\bibitem[P2]{pr2} Yu.~G.~Prokhorov.
\textit{On Fano-Enriques threefolds}.
Preprint arXiv:math.AG/0604468.

\bibitem[R1]{re1} M.~Reid. 
\textit{Young person's guide to canonical singularities}.
Algebraic geometry, Bowdoin, 1985 (Brunswick, Maine, 1985).
Proc. Sympos. Pure Math. \textbf{46}, Part 1, 345--414. 
Amer. Math. Soc., Providence, RI, 1987. 

\bibitem[R2]{re2} M.~Reid. 
\textit{Update on 3-folds}.
Proceedings of the International Congress of Mathematicians, Vol. II (Beijing, 2002),  513--524.
Higher Ed. Press, Beijing, 2002.

\bibitem[Sa]{sa} T.~Sano. 
\textit{On classification of non-Gorenstein $\QQ$-Fano $3$-folds of Fano index $1$}. 
J. Math. Soc. Japan \textbf{47}, (1995) 369--380.

\bibitem[Sc]{sc} G.~Scorza. 
\textit{Sopra una certa classe di variet\`a razionali}. 
Rend. Circ. Mat. Palermo \textbf{28}, (1909) 400--401.

\bibitem[Sh1]{sh1} V.~V.~Shokurov. 
\textit{Smoothness of a general anticanonical divisor on a Fano variety}. 
Izv. Akad. Nauk SSSR Ser. Mat. \textbf{43}, (1979) 430--441. 

\bibitem[Sh2]{sh2} V.~V.~Shokurov. 
\textit{The existence of a line on Fano varieties}.  
Izv. Akad. Nauk SSSR Ser. Mat. \textbf{43} (1979) 922--964, 968.

\bibitem[Wa]{wa} J.~Wahl. 
\textit{Introduction to Gaussian maps on an algebraic curve}.
Complex Projective Geometry, Trieste-Bergen 1989, London Math. Soc. Lecture Notes Ser. 
\textbf{179}. Cambridge Univ.\ Press, Cambridge 1992, 304-323.  

\bibitem[Za]{za} F.~L.~Zak. 
\textit{Some properties of dual varieties and their applications in projective geometry}.
Algebraic geometry (Chicago, IL, 1989), 273--280. Lecture Notes in Math. \textbf{1479}. Springer,
Berlin, 1991.
 
\end{thebibliography}
\end{document}